\documentclass[a4paper,11pt]{article}

\usepackage[english]{babel}
\usepackage{amssymb}
\usepackage{amsmath,amsfonts,latexsym}
\usepackage[all]{xy}

\font\ggordo=cmbx10 scaled \magstep2

\font\nineit=cmti9

\font\smc=cmcsc10

\newtheorem{proposition}{Proposition}
\newtheorem{definition}{Definition}
\newtheorem{theorem}{Theorem}

\newtheorem{corollary}{Corollary}
\newtheorem{lemma}{Lemma}
\newtheorem{remark}{Remark}

{\vspace{3.3mm}
\noindent{\bf #1}\it}%
{\vspace{3.3mm}}

\newcommand{\slps} {straight--line programs}
\newcommand{\cts} {correct test sequence}

\def\ifm#1#2{\relax \ifmmode#1\else#2\fi}

\newcommand{\klk}    {\ifm {,\ldots,} {$,\ldots,$}}

\newcommand{\plp}    {\ifm {+\cdots+} {$+\ldots+$}}

\newcommand{\cqfd}{\hfill\vbox{\hrule height 5pt width 5pt }\bigskip}

\newcommand{\md}{\mathcal{D}}
\newcommand{\mau}{\mathcal{U}}
\newcommand{\mo}{\mathcal{O}}
\newcommand{\mc}{\mathcal{C}}
\newcommand{\ma}{\mathcal{A}}
\newcommand{\mb}{\mathcal{B}}

\newcommand{\omo}{\overline{\mathcal{O}}}
\newcommand{\omd}{\overline{\mathcal{D}}}

\newcommand{\hmd}{\widehat{\mathcal{D}}}

\newcommand{\N}{\ifm {\mathbb N} {$\mathbb N$}}
\newcommand{\C}{\ifm {\mathbb C} {$\mathbb C$}}
\newcommand{\Z}{\ifm {\mathbb Z} {$\mathbb Z$}}
\newcommand{\Zpos}{\ifm {\Z_{\ge 0}} {$\Z_{\ge 0}$}}
\newcommand{\Q}{\ifm {\mathbb Q} {$\mathbb Q$}}
\newcommand{\Pe}{\ifm {\mathbb P} {$\mathbb P$}}
\newcommand{\R}{\ifm {\mathbb R} {$\mathbb R$}}
\newcommand{\A}{\ifm {\mathbb A} {$\mathbb A$}}
\newcommand{\ok}{\ifm {\overline{k}} {$\overline{k}$}}
\newcommand{\thetah}{\ifm {\widehat{\theta}} {$\widehat{\theta}$}}

\newcommand{\om}[2]   {{#1}_1 \klk {#1}_{#2}}

\newcommand{\xon}    {\ifm {\om X n} {$\om X n$}}

\newcommand{\yot}    {\ifm {\om Y t} {$\om Y t$}}

\newcommand{\gammam}    {\ifm {\om \gamma m} {$\om \gamma m$}}

\newcommand{\Ring}[2]  {\ifm {{#1}[#2]} {${#1}[#2]$}}

\newcommand{\okyot}   {\Ring {\ok} {\yot}}

\newcommand{\iseq}[1] {(#1_i)_{i\in\N}}
\newcommand{\nseq}[1] {(#1_n)_{n\in\N}}

\begin{document}

\title{\ggordo The Hardness of Polynomial Equation Solving\footnote{
Research was partially supported by the following Argentinian,
French and Spanish grants~: UBACyT X198, PIP CONICET 2461, UNGS
30/3003, CNRS FRE 2341 MEDICIS, DGCYT BFM2000--0349, HF-1999-055 and
ALA 01--E3/02.}}

 \author{ D. Castro$^{1}$,
 M. Giusti$^{2}$, J. Heintz$^{1,\,3,\,4}$,
 G. Matera$^{4,\,5}$, L.M. Pardo$^{1,\,2}$ }

%\date{\today}

 \maketitle

\addtocounter{footnote}{1}\footnotetext{\nineit Depto. de
Matem\'aticas, Estad\'{\i}stica y Computaci\'on, Facultad de
Ciencias, Universidad de Cantabria, E-39071 Santander, Spain.}
\addtocounter{footnote}{1}\footnotetext{\nineit Laboratoire GAGE,
Ecole Polytechnique, F-91128 Palaiseau Cedex, France.}
\addtocounter{footnote}{1}\footnotetext{\nineit Depto. de
Matem\'aticas, Facultad de Ciencias Exactas y Naturales,
Universidad de Buenos Aires, Ciudad Universitaria, Pabell\'on I
(1428) Buenos Aires, Argentina.} \addtocounter{footnote}{1}
\footnotetext{\nineit  Member of the National Council of Science
and Technology (CONICET), Argentina.} \addtocounter{footnote}{1}
\footnotetext{\nineit  Instituto de Desarrollo Humano, Universidad
Nacional de General Sarmiento, Campus Universitario, Jos\'e M.
Guti\'errez 1150 (1613) Los Polvorines, Pcia. de Buenos Aires,
Argentina.}

\begin{center}
{\em Dedicated to Michel Demazure}
\end{center}
\bigskip

\begin{quotation}

{\it ``.......Il est fr\'equent, devant un probl\`eme concret, de
trouver un th\'eor\`eme qui ``s'applique presque"....... . Le
r\^ole des contre--exemples est justement de d\'elimiter le
possible, et ce n'est pas par perversit\'e (ou en tout cas pas
totalement) que les textes math\'ematiques exhibent des monstres"}
{\smc M. Demazure}, 1987

\end{quotation}

\begin{abstract}
Elimination theory is at the origin of algebraic geometry in the
19-th century and deals with algorithmic solving of multivariate
polynomial equation systems over the complex numbers, or, more
generally, over an arbitrary algebraically closed field. In this
paper we investigate the {\it intrinsic} sequential time
complexity of {\it universal} elimination procedures for {\it
arbitrary continuous} data structures encoding input and output
objects of elimination theory (i.e.  polynomial equation systems)
and admitting the representation of certain limit objects.

Our main result is the following: let be given such a data
structure and together with this data structure a universal
elimination algorithm, say $\mathcal{P}$, solving arbitrary
parametric polynomial equation systems. Suppose that the algorithm
$\mathcal{P}$ avoids ``unnecessary'' branchings and that
$\mathcal{P}$ admits the efficient computation of certain natural
limit objects (as e.g. the Zariski closure of a given
constructible algebraic set or the parametric greatest common
divisor of two given algebraic families of univariate
polynomials). Then $\mathcal{P}$ cannot be a polynomial time
algorithm.

The paper contains different variants of this result and discusses
their practical implications.

\end{abstract}

\smallskip


\tableofcontents

\section{Introduction. Basic notions.}
\label{introduccion}
Complexity theory deals with the efficiency of answering
mathematical questions about mathematical objects. In this
context, mathematical objects happen usually to posses a
{\em unique} encoding
in a previously fixed data structure (e.g. integers are encoded
by their bit
representation, polynomials by their coefficients, etc.). Once a
data structure is fixed, standard complexity theory searches for
an efficient algorithm answering the mathematical questions under
consideration and tries to certify the optimality of this
algorithm.

However, things become somewhat more complicated in the particular
case of geometric elimination theory (polynomial equation solving
in algebraically or real closed fields). Complexity theory for
geometric elimination requires {\em simultaneous} optimization
of data structures {\em and} algorithms. In order to illustrate
this statement, let us consider the following first order formula,
say $\Phi$, belonging to the language of the elementary theory of
algebraically
 closed fields of
 characteristic zero:
$$\!\!(\exists X_1)\cdots(\exists X_n)\big(X_1-T-1=0\,\wedge\,
X_1^2-X_2=0\,\wedge\cdots \wedge\, X_{n-1}^2-X_n=0\,\wedge\,
Y=X_n^2\big).$$ The formula $\Phi$ contains two free variables,
namely $T$ and $Y$. Moreover $\Phi$ is logically equivalent to the
following quantifier--free formula, which we denote by $\Psi$:
$$Y-\sum_{i=0}^{2^n}
\begin{pmatrix}
2^n \\ i
\end{pmatrix} T^i=0.$$
If we choose as our data structure the standard dense or sparse
encoding of polynomials by their coefficients, then $\Phi$
has length $O(n)$, whereas the length of $\Psi$ exceeds $2^n$.
However, if we encode polynomials by {\em arithmetic circuits} (or \slps),
then $\Phi$ and $\Psi$ happen both to be of length $O(n)$,
since the polynomial $\displaystyle\sum_{i=0}^{2^n}\begin{pmatrix}
2^n \\ i \end{pmatrix} T^i=(1+T)^{2^n}$ can be evaluated in $n+1$
steps, using iterated squaring.

For the dense (or sparse) representation of polynomials
{\em superexponential} (sequential) time is necessary (and
sufficient) in order to eliminate a {\em single} quantifier
block (see e.g. \cite{CaGaHe89}, \cite{DiFiGiSe91},
\cite{GrVo88}, \cite{HeRoSo89A}, \cite{Canny88}), whereas the elimination
of an {\em arbitrary} number of quantifier blocks requires
{\em doubly exponential} time in this data structure
(see \cite{Heintz83}, \cite{Weispfenning88}, \cite{DaHe88},
\cite{FiGaMo90}, \cite{FiGaMo90a}, \cite{MoPa93a} for lower
and upper complexity bounds and \cite{Ierardi89},
\cite{HeRoSo90A}, \cite{Renegar92} for upper complexity bounds
only).

The existing superexponential (mostly Gr\"obner basis) algorithms
for the elimination of a single quantifier block are often
asymptotically optimal for the dense and sparse encoding of
polynomials. Nevertheless their complexity makes them infeasible
for real world sized problems (however not so for impressive
software demos). Moreover, a simple minded Gr\"obner basis
approach to the elimination of a single block of quantifiers may
lead to a doubly exponential complexity.
\medskip

This situation suggests that these elimination algorithms require
alternative data structures if one wishes to improve their
complexity behaviour substantially. This observation led in the
past to the idea of using arithmetic circuits for the
representation of the polynomials occurring in the basic
elimination procedures of algebraic and semialgebraic geometry
(see \cite{HeSc82}, \cite{HeSi81} and \cite{Kaltofen88} for an
early application of this idea). This change of data structure
allowed in a first attempt to reduce the complexity of the
elimination of a single block of quantifiers from
superexponential to single exponential time (\cite{GiHe93},
\cite{GiHeSa93}, \cite{FiGiSm95}, \cite{KrPa94}, \cite{KrPa96},
\cite{Matera99}). However, the corresponding algorithms required
the dense representation of the input polynomials and returned a
circuit encoding of the output polynomials. Therefore these
algorithms were unadapted to successive elimination of several
quantifier blocks (see \cite{PuSa98} for more details) and unable
to profit from a possible special geometric feature of the input
system.

In a second attempt (\cite{GiHeMoPa95}, \cite{Pardo95},
\cite{GiHeMoMoPa98}, \cite{GiHaHeMoMoPa97}, \linebreak
\cite{GiHeMoPa97}), this problem could be settled by means of a
new elimination procedure which transforms a given circuit
representation of the input polynomials into a circuit
representation of the output polynomials. The time complexity of
this new procedure is roughly the circuit size of the input
polynomials multiplied by a {\em polynomial} function of a certain
geometric invariant of the input system, called its {\em degree}.
Let us observe that the degree is always bounded by the
B\'ezout--number of the input system and happens often to be
considerably smaller.

For worst case input systems, the new algorithm becomes
polynomial in the B\'ezout--number of the system, and this
was the first time that this complexity goal could be reached
without hiding an exponential extra factor (compare \cite{MoPa97},
\cite{Rojas00}).

Afterwards the new algorithm and its data structure was extended
and refined in \cite{HeKrPuSaWa00}, \cite{GiSc99},
\cite{HaMoPaSo00}, \cite{HeMaWa01}, \cite{GiLeSa01},
\cite{Schost00}, \linebreak \cite{Lecerf01}, and in
\cite{BaGiHeMb97}, \cite{BaGiHeMb01} it was adapted to the problem
of polynomial equation solving over the {\em reals}. A successful
implementation of the full algorithm (\cite{Lecerf00}) is based on
\cite{GiLeSa01}. A partial implementation of the algorithm
(including basic subroutines) is described in \cite{BrHeMaWa02}
(see also \cite{CaHeLlMa00}). So far the account of successive
improvements of data structures and algorithms for {\em symbolic}
elimination. The complexity aspect of {\em numeric} elimination
was treated in a series of papers (\cite{ShSm93}, \cite{ShSm93b},
\cite{ShSm93c}, \cite{ShSm96}, \cite{ShSm94}, \cite{CuSm99}; see
also \cite{BlCuShSm98}). In \cite{CaHaMoPa01} and
\cite{CaMoPaSa02} the bit complexity aspect of the above mentioned
symbolic and numeric algorithms was analyzed and compared. Taking
bit complexity and the bit representation of rational numbers into
account, it turns out that a suitable numerical adaptation of the
above mentioned new symbolic elimination algorithm has the best
complexity performance between all known numerical elimination
algorithms. Therefore we shall limit our attention in this paper
to {\em symbolic} elimination procedures.
\medskip

Let us now briefly sketch the known lower bound results for the
complexity of arithmetic circuit based procedures for the
elimination of a single quantifier block. Any such
elimination algorithm which is {\em geometrically robust} in
the sense of \cite{HeMaPaWa98} requires necessarily exponential
time on infinitely many inputs. Geometric robustness is a very
mild condition that is satisfied by all known (symbolic)
elimination procedures.

Moreover, suppose that there is given an algorithm for the
elimination of a single quantifier block and suppose that this
algorithm is in a suitable sense ``universal", avoiding
``unnecessary branchings" and able to compute Zariski closures of
constructible sets and ``parametric" greatest common divisors of
algebraic families of univariate polynomials. Then necessarily
this algorithm has to be robust and hence of non--polynomial time
complexity \cite{GiHe01}. In particular, any ``reasonable" and
``sufficiently general" procedure for the elimination of a
single quantifier block produces geometrically robust
arithmetic circuit representations of suitable elimination
polynomials and therefore outputs of non--polynomial size in
worst case.
\medskip

In this paper we are going to argue that the non--polynomial
complexity character of the known symbolic geometric elimination
procedures is not a special feature of a particular data structure
(like the dense, sparse or arithmetic circuit encoding of
polynomials), but rather a consequence of the information encoded
by the respective data structure (see Theorem
\ref{theorem:section5.3} below).

%-----------------------------------------------------------------
%-----------------------------------------------------------------
%-----------------------------------------------------------------
%-----------------------------------------------------------------

\subsection{Data structures for geometric objects.}
\label{section1.1} Informally, we understand by a {\em data
structure} a class, say $\mathcal{D}$, of ``simple" mathematical
objects which encode another class, say $\mathcal{O}$, of
``complicated" ones. An element $D\in\mathcal{D}$ which encodes a
mathematical object $O\in\mathcal{O}$ is called a {\em code} of
$\mathcal{O}$. The data structure $\mathcal{D}$ is supposed to be
embedded in a context, where we may process its elements, in order
to answer a certain catalogue of well defined {\em questions}
about the mathematical objects belonging to $\mathcal{O}$ (or
about the object class $\mathcal{O}$ itself). Of course, the
choice of the data structure $\mathcal{D}$ depends strongly on the
kind of potential questions we are going to ask and on the time we
are willing to wait for the answers.

The mathematical objects we are going to consider in this paper
will always be polynomial functions or algebraic varieties and
their codes will always belong to suitable affine ambient spaces.
The {\em size} of a code is measured by the dimension of its
ambient space.

The (optimal) encoding of {\em discrete} (e.g. finite) sets of
mathematical objects is a well known subject in theoretical
computer science and the main theme of Kolmogorov complexity
theory (\cite{LiVi93}; see also \cite{Borel48}).
\medskip

This paper addresses the problem of optimal encoding of {\em
continuous} classes of mathematical objects. The continuous case
differs in many aspects from the discrete one and merits
particular attention. Any object class $\mathcal{O}$ we are
considering in this paper will possess a {\em natural topology}
and may be thought to be embedded in a (huge) affine or projective
ambient space. The given topology of $\mathcal{O}$ becomes always
induced by the Zariski (or strong) topology of its ambient space.
In this paper, the {\em closure} $\overline{\mathcal{O}}$ of the
object class $\mathcal{O}$ in its ambient space will generally
have a natural interpretation as a class of objects of the same
nature as $\mathcal{O}$. In this sense we shall interpret a given
element of $\overline{\mathcal{O}}\setminus\mathcal{O}$ as a {\em
limit} (or {\em degenerate}) object of $\mathcal{O}$. We shall
always suppose that data structures, object classes and graphs of
encodings form {\em constructible} subsets of their respective
ambient spaces.

If $\mathcal{O}$ is for example a class of equidimensional closed
subvarieties of fixed dimension and degree of a suitable
projective space, then the topology and the ambient space of
$\mathcal{O}$ may be given by the Chow coordinates of the objects
of $\mathcal{O}$. Or if $\mathcal{O}$ is the class of polynomial
functions of bounded arithmetic circuit complexity $L$, then
$\mathcal{O}$ is contained in a finite dimensional linear subspace
of the corresponding polynomial ring and has therefore a natural
topology. The limit objects of $\mathcal{O}$ are then those
polynomials which have approximative complexity at most $L$ (see
\cite{Alder84}, \cite{BuClSh97} and Section \ref{circuitcorrect} for details).

Let be given a data structure $\mathcal{D}$ encoding an object
class $\mathcal{O}$. By assumption, $\mathcal{D}$ is embedded in a
suitable affine or projective ambient space from which
$\mathcal{D}$ inherits a natural topology. We shall always assume
that $\mathcal{D}$ encodes $\mathcal{O}$ {\em continuously} or
{\em holomorphically} (see Section \ref{holocont} for precise,
mathematical definitions). However, in order to capture the
important case of the arithmetic circuit representation of
polynomials, we shall not insist on the injectivity of the given
encoding. More precisely, we say that $\md$ encodes $\mo$ {\em
injectively} or {\em unambiguously} if for any object
$O\in\mathcal{O}$ the data structure $\mathcal{D}$ contains a
single element encoding $O$ (otherwise we call the encoding {\em
ambiguous}).
\medskip

A fundamental problem addressed in this paper is the
following:

\smallskip \noindent given an ``efficient" (i.e. short) data
structure $\mathcal{D}$ encoding the object class $\mathcal{O}$,
how may we find another data structure $\overline{\mathcal{D}}$
encoding the object class $\overline{\mathcal{O}}$? How does the
{\em size} of $\overline{\mathcal{D}}$ (i.e. the size of its
codes) depend on the size of $\mathcal{D}$? \smallskip

In Theorem \ref{theorem:encode1}, Corollary \ref{encode3},
Corollary \ref{coro:encode2} and in Section \ref{section4real}
below we shall see that the solution of this problem depends
strongly on the type of questions about the object class $\omo$
which the data structure $\overline{\mathcal{D}}$ allows to
answer. \medskip

This leads us to the subject of the questions we wish to
be answered by a given data structure $\mathcal{D}$
encoding a
given object class $\mathcal{O}$. We shall always
require that any element $D\in\mathcal{D}$ encoding an object
$O\in\mathcal{O}$ contains enough information in order to
distinguish $O$ from other elements of the object class
$\mathcal{O}$. Thus a typical question we wish to be answered by
the data structure $\mathcal{D}$ is the following: \smallskip

\noindent let $D$ and $D'$ be two elements of $\mathcal{D}$
encoding two objects $O$ and $O'$ of $\mathcal{O}$. Are $O$ and
$O'$ identical?
\smallskip

In other words, we require to be able to deduce whether $O=O'$
holds by means of processing the codes $D$ and $D'$. We call this
problem the {\em identity question} associated to the data
structure $\mathcal{D}$. A common way to solve this identity
question consists of the transformation of the (supposedly
ambiguous) data structure $\mathcal{D}$ in a new one, which
encodes the objects of $\mathcal{O}$ injectively.\medskip

Another typical question arises in the following context:

\noindent suppose additionally that the object class $\mathcal{O}$
consists of (total) functions which can be evaluated on a
continuous (or discrete) domain $R$. Suppose furthermore that we
have free access to any element of $R$. Let $O$ be a given element
of the object class $\mathcal{O}$ and let $D\in\mathcal{D}$ be an
arbitrary code of $O$.

The question we wish to be answered by the data structure
$\mathcal{D}$ about the object $O$ is the following:

\smallskip

\noindent for any given argument value $r\in R$, what is the
function value $O(r)$?
\smallskip

In other words, we require to be able to compute the function
value $O(r)$ by means of processing the code $D$ and the argument
value $r$. We call this problem the {\em value question}
associated to the data structure $\mathcal{D}$. Of course, for a
class $\mathcal{O}$ of polynomial functions whose number of
variables and degree was previously bounded, the value question
for a continuous domain $R$ can be reduced by means of
interpolation techniques to the value question for a
discrete domain, namely to the task of determining, for
any monomial $M$ and any polynomial function $O\in\mathcal{O}$,
the coefficient of $M$ in the polynomial $O$.

%-----------------------------------------------------------------
%-----------------------------------------------------------------
%-----------------------------------------------------------------
%-----------------------------------------------------------------

\subsection{The r\^ole of data structures in elimination theory.}
\label{section1.2} In algebraic geometry, polynomial equation
systems are the ``simple" mathematical objects which encode the
real objects of interest: algebraic varieties or schemes. Except
for the particular case of hypersurfaces, there is no {\em a
priori} privileged {\em canonical} equation system that defines a
given algebraic variety. It depends on the questions we are going
to ask about the given variety, whether we shall feel the need to
transform a given equation system into a new, better suited one
for answering our questions.

As far as possible, we wish just to modify the syntactical form of
our equations, without changing their meaning, represented by the
underlying variety or scheme. Let us explain this in two different
situations.

Very often the new equation system we are looking for is {\em
uniquely determined} by the underlying variety or scheme and the
syntactical requirements the new system has to satisfy. We meet
this situation in the particular case of the (reduced) Gr\"obner
basis of a given ideal (representing a scheme) for a previously
fixed monomial order. The monomial order we shall choose depends
on the kind of questions we are going to ask about the given
scheme: we choose an (e.g. lexicographical) elimination order if
we wish to ``solve" the given equation system (i.e. uncouple its
variables) or we choose a graded order if we wish to compute the
Hilbert polynomial (the dimension and the degree) of the given
(projective) scheme, etc. If we want to analyze a given scheme or
variety by means of deformations, suitable (i.e. flat) equation
systems, like Gr\"obner bases, are even mandatory (\cite{BaMu93}).
Although Gr\"obner bases are able to answer all typical questions
about the variety or scheme under consideration, they are not well
suited for the less ambitious task of polynomial equation solving
(this constitutes the main elimination problem the paper is
focusing on).

Since Gr\"obner bases are able to answer too many questions about
the scheme or variety  they define, they may become difficult to
encode: a complete intersection ideal given by low degree binomial
equations, may have a Gr\"obner basis of doubly exponential degree
for a suitable elimination order, whereas it is possible to solve
the corresponding elimination problem in singly exponential time
using only polynomials of singly exponential degree (see
\cite{DiFiGiSe91}, \cite{KrPa96}, \cite{HaMoPaSo00}).\medskip

Let us consider another case of this general situation:

\noindent for the particular task of polynomial equation solving
it suffices to replace the original algebraic variety (which is
supposed to be equidimensional) by a birationally equivalent
hypersurface in a suitable ambient space. This hypersurface and
its minimal equation may be produced by means of generic linear
projections (see e.g. \cite{Kronecker82}, \cite{ChGr83},
\cite{GiMo89}, \cite{Canny88}, \cite{CaGaHe89}, \cite{DiFiGiSe91},
\cite{GiHe91}, \cite{GiHe93}, \cite{KrPa96}, \cite{GiHeMoMoPa98},
\cite{GiHaHeMoMoPa97}) or by means of dual varieties (see
\cite{GeKaZe94} and the references cited there). The minimal
equation encodes the necessary information about the dimension and
degree of the original algebraic variety and about a suitable set
of independent variables. However, as a consequence of B\'ezout's
Theorem, the degree of the canonical output equation may increase
exponentially with respect to the degree of the given input
equations if we apply this strategy of elimination. We call an
elimination procedure {\em Kronecker--like} if in terms of
suitable data structures, the procedure computes from the
representation of equations of the given variety a representation
of the minimal equation of the corresponding hypersurface (see
Section \ref{section4.1} for more details).
\medskip

In either case of this general situation, we need an {\em input
object} (a polynomial equation system describing an algebraic
variety or scheme) and an {\em output object} that describes the
same variety or scheme (or a birationally equivalent one) and
satisfies some additional syntactical requirements (allowing e.g.
the uncoupling of the variables of the original system). The
corresponding elimination problem maps input objects to output
objects. Since an output object may have degree exponential in the
degree of the corresponding input object we are led to ask about
short encodings of high degree polynomials in few variables. In
this context let us mention the main outcome of
\cite{GiHeMoMoPa98}, namely the observation that using the
arithmetic circuit representation, elimination polynomials (i.e.
the output objects of Kronecker--like procedures) have always size
polynomial in their degree, whereas the size of their sparse (or
dense) representation may become exponential in this quantity. But
unfortunately, elimination polynomials may have exponential
degrees. This inhibits the elimination procedure of
\cite{GiHeMoMoPa98} and \cite{GiHaHeMoMoPa97} to become polynomial
in the input length, at least in worst case.\medskip

We consider therefore in more generality the following task:

\smallskip \noindent {\em let be given an elimination problem and a data
structure $\mathcal{D}$ encoding the input objects. Find a data
structure $\mathcal{D}^*$ encoding the corresponding output
objects and an elimination algorithm $\mathcal{P}$ which maps
input codes belonging to $\mathcal{D}$ to output codes belonging
to $\mathcal{D}^*$ and solves the given elimination
problem}.\smallskip

In this terminology, the main problem this paper tries to
solve can be formulated as follows:

\smallskip \noindent {\em is it possible to find in the given situation
a data structure $\mathcal{D}^*$ and a continuous algorithm
$\mathcal{P}$ (in the sense specified in Sections \ref{section2}
and \ref{section4}) such that the size of each output code
belonging to $\mathcal{D}^*$ is only polynomial in the size
of the corresponding input code belonging to $\mathcal{D}$? Under
which circumstances do such a data structure $\mathcal{D}^*$ and
such an algorithm $\mathcal{P}$ exist and under which
circumstances do they not?}
\smallskip

As mentioned before, the solution of this problem
 depends strongly on the
questions about the output objects
 we wish to be answered by the output data
structure
 $\mathcal{D}^*$.

In view of the methodological progress made in
\cite{GiHeMoMoPa98}, \cite{GiHaHeMoMoPa97}, \cite{GiHeMoPa97},
\cite{HeKrPuSaWa00}, \cite{GiLeSa01} and \cite{HeMaWa01} (leading
to a substantial improvement of previously known complexity
bounds) and motivated by our interest in lower complexity bounds,
we limit our attention to basic and relatively simple elimination
problems of the following type:
\begin{itemize}
\item[({\em i}\,)]
Let be given a zero--dimensional algebraic variety $V$ by an input
equation system in $n$ variables $X_1\klk X_n$ and let be given a
supplementary input polynomial $F$ in these $n$ variables and
possibly some additional parameters $U_1\klk U_r$. Let $X:=(\xon)$
and $U:=(U_1\klk U_r)$ and let us suppose that the input equation
system and $F$ have short encodings in a previously fixed input
data structure $\mathcal{D}$. The problem is to find an output
data structure $\mathcal{D}^*$ and a Kronecker--like elimination
procedure $\mathcal{P}$ such that $\mathcal{P}$ associates to each
input code of $\mathcal{D}$ representing a specialization $u$ of
the parameters $U$ of $F$, an output code of $\mathcal{D}^*$
representing the canonical elimination polynomial of  $F(u,X)$
with respect to the given variety $V$ (see Sections
\ref{section4.1} and \ref{section4.2} for definitions and an
example).
\item[({\em ii}\,)]
Let be given a class $\mathcal{O}$ of mathematical objects and a
data structure $\mathcal{D}$ encoding the object class
$\mathcal{O}$. Suppose that $\mathcal{O}$ and $\mathcal{D}$
satisfy all general assumptions we made before on this kind of
mathematical entities. Find a data structure
$\mathcal{D}^*$ and a procedure $\mathcal{P}$ such that
$\mathcal{D}^*$ encodes the topological closure
$\overline{\mathcal{O}}$ of the object class $\mathcal{O}$ and
such that $\mathcal{P}$ maps any element of $\mathcal{D}$ encoding
a given object $O$ of $\mathcal{O}$ to an element of
$\mathcal{D}^*$ encoding the same object $O$.
\end{itemize}

In case of problems of type ({\em ii}\,), a typical example of
such a procedure $\mathcal{P}$ for arithmetic circuit represented
rational functions of bounded degree $d$ is the ``Vermeidung von
Divisionen" algorithm of \cite{Strassen73} (see also
\cite{KrPa96}). In this case the limit objects are polynomials of
degree at most $2d+1$ (see \cite{Alder84} for details). Another
example of such a procedure is the transformation (by means of
``tensoring") of approximative algorithms for matrix
multiplication into exact ones \cite{BuClSh97}.\medskip

In elimination theory one meets very natural non--closed input
object classes with limit objects not encoded by the given input
data structure. However, such a limit object may possess
a well--defined output object. In this case one may require that
the given output data structure is able to encode this output
object. This is the typical context where a problem of type ({\em
ii}\,) arises in elimination theory.\medskip

Another context, related to approximation and interpolation
theory is the following:\smallskip

\noindent let $\omega: \mathcal{D}\to \mathcal{O}$ be a given
encoding of an object class of $t$--variate polynomial functions
over $\C$ having degree bounded by an a priori constant $\Delta$.
Consider $\mathcal{O}$ as a metric space equipped with the
corresponding (strong) topology. Suppose that the encoding
$\omega$ is holomorphic, allowing for each code $D\in
\mathcal{D}$ and any argument $r\in \C^t$ the computation of the
value $\omega(D)(r)$ using a fixed number of arithmetic
operations in $\C$ (see Section \ref{holocont} for details). Let
$O\in \overline{\mathcal{O}}\setminus \mathcal{O}$ be a limit
object of $\mathcal{O}$, let $(D_i)_{i\in \N}$ be a sequence of
codes of $\mathcal{D}$ such that the sequence
$(\omega(D_i))_{i\in \N}$ converges to the limit object $O$ and
let $r\in \C^t$ be a given argument. As one easily sees, $O$ is
again a $t$--variate polynomial over $\C$ of degree at most
$\Delta$ and therefore the value $O(r)$ is well defined.
Moreover, the sequence of complex numbers $(\omega(D_i)(r))_{i\in
\N}$ converges to the value $O(r)$. However, the
convergence rate of $(\omega(D_i)(r))_{i\in \N}$ will typically
depend on the argument $r$. Our goal is to compute the value
$O(r)$ using only a {\em fixed} number of arithmetic operations
and limit processes in $\C$ for sequences which do not depend on
the argument $r$. We reach this goal if we are able to solve
in this context problem $(ii)$ by a data structure which answers
the value question.
\medskip

All known algorithms solving problem ({\em i}\,) or, limited to
the context of classical elimination theory, problem ({\em ii}\,),
possess branching--free versions of the same order of complexity.
We shall therefore consider only {\em branching--free} algorithms
for the solution of these two elimination problems. In this case,
we shall always assume that our output codes depend {\em
holomorphically} (or at least continuously) on our input
codes (see Section \ref{models} for details).\medskip

An elimination algorithm is called {\em universal} if it solves
for appropriate input and output data structures any standard
elimination problem on arbitrary inputs consisting of boolean
combinations of {\em parameter dependent} polynomial equations. A
universal elimination algorithm is called {\em
branching--parsimonious} if it avoids branchings for the solution
of suitable instances of problems of type ({\em i}\,) and ({\em
ii}\,).
\medskip

This paper is organized as follows:

In Section \ref{section2} we introduce the language and tools from
algebraic geometry and algebraic complexity theory we are going to
use in this paper. In Section \ref{section3} we discuss different
types of encodings of object classes: holomorphic, robust and
continuous ones. We prove our first main result, namely Theorem
\ref{theorem:encode1}, saying that any holomorphic (ambiguous)
encoding may be replaced by a continuous and unambiguous one of
similar size. We retake the subject of this section in an appendix
of this paper, namely in Section \ref{appendix}, generalizing
Theorem \ref{theorem:encode1} to Corollary \ref{coro:encode2} and
estimating the VC--dimension of a given, holomorphically encoded
object class in terms of the size of its encoding.

In Section \ref{section4real} we introduce the main technique we
are going to apply in this paper in order to prove lower bounds
for robust encodings of specific object classes. We exemplify this
technique by two fundamental examples.

In Section \ref{section4} we apply the tools developed in the
preceding sections to elimination theory. We introduce the notion
of a robust elimination procedure for flat families of
zero--dimensional elimination problems and show that any robust
elimination procedure requires necessarily exponential
(sequential) time on infinitely many inputs (Theorem
\ref{lemma:section4.2}). This result is then used in order to
prove the second main result of this paper, namely Theorem
\ref{theorem:section5.3}, which may be paraphrased as follows:
\smallskip

\noindent Suppose that there is given a universal,
branching--parsimonious elimination procedure $\mathcal{P}$ which
is also able to solve in the context of elimination theory
suitable problems of the above type ({\em ii}\,). In particular,
we suppose that the procedure $\mathcal{P}$ is able to eliminate
quantifiers in parametric existential first order formulas of the
language of the elementary theory of algebraically closed fields
of characteristic zero and that $\mathcal{P}$ is able to compute
equations for the Zariski closure of any given constructible set
and the generically square--free parametric greatest common
divisor of any given algebraic family of univariate polynomials
(see Sections \ref{models}, \ref{section5.2} and \ref{section5.3}
for precise, mathematical definitions). Then, the elimination
procedure $\mathcal{P}$ cannot be of polynomial (sequential) time
complexity.
\medskip

In conclusion, a universal, branching--parsimonious
procedure
for the elimination of a single existential quantifier block
which is able to solve suitable problems of
type ({\em ii}\,) cannot be polynomial. Let us remark that all
known universal elimination procedures satisfy this requirement
since they are based on subroutines (in particular greatest common
divisor computations) which behave well under specialization.

All these results are formulated in an {\em exact}
computation model which allows to represent all known symbolic
and seminumeric elimination procedures (based on the sparse or
dense or the arithmetic circuit representation of polynomials).

%-----------------------------------------------------------------
%-----------------------------------------------------------------
%-----------------------------------------------------------------
%-----------------------------------------------------------------
%-----------------------------------------------------------------
%-----------------------------------------------------------------
%-----------------------------------------------------------------
%-----------------------------------------------------------------

\section{Notions and notations.}
\label{section2} \subsection{Language and tools from algebraic
geometry.} Let $k$ be an infinite, perfect field which we think to
be ``effective" with respect to arithmetic operations as
addition/subtraction, multiplication/division and extraction of
$p$--th roots in case $k$ has positive characteristic $p$. Let
$\overline{k}$ be an algebraically closed field containing
$k$ (in the sequel we shall call such a field {\em an algebraic
closure of $k$}.
\medskip

Most of the statements and arguments of this paper will be
independent of the characteristic of $k$. Therefore the reader may
assume without loss of generality that $k$ is of characteristic
zero. For the sake of simplicity we shall assume in this case
$k:=\Q$ and $\overline{k}=\C$. We denote by $\N$ the set of
natural numbers and by $\Zpos$ the set of nonnegative integers.
\medskip

Fix $n\in\Zpos$ and let $X_0\klk X_n$ be indeterminates
over $k$. We denote by $\A^n:=\A^n(\overline{k})$ the
$n$--dimensional affine space and by $\Pe^n:=\Pe^n(\overline{k})$
the $n$--dimensional projective space over $\overline{k}$. The
spaces $\A^n$ and $\Pe^n$ are thought to be endowed with their
respective Zariski topologies over $k$ and with their respective
sheaves of $k$--rational functions with values in $\overline{k}$.
Thus the points of $\A^n$ are elements $(x_1\klk x_n)$ of
$\overline{k}$ and the points of $\Pe^n$ are (non uniquely)
represented by nonzero elements $(x_0\klk x_n)$ of
$\overline{k}^{n+1}$ and denoted by $(x_0:\dots: x_n)$. The
indeterminates $X_1\klk X_n$ are considered as the coordinate
functions of the affine space $\A^n$. The coordinate ring (of
polynomial functions) of $\A^n$ is identified with the polynomial
ring $k[X_1\klk X_n]$. Similarly we consider the (graded)
polynomial ring $k[X_0\klk X_n]$ as the projective coordinate ring
of $\Pe^n$. Consequently we represent rational functions of
$\Pe^n$ as quotients of homogeneous polynomials of equal degree
belonging to $k[X_0\klk X_n]$. Let $F_1\klk F_s$ be polynomials
which belong to $k[X_1\klk X_n]$ or are homogeneous and belong to
$k[X_0\klk X_n]$. We denote by $\{F_1=0\klk F_s=0\}$ or $V(F_1\klk
F_s)$ the algebraic set of common zeroes of the polynomials
$F_1\klk F_s$ in $\A^n$ and $\Pe^n$ respectively. We consider the
set $V:=\{F_1=0\klk F_s=0\}$ as (Zariski--)closed (affine or
projective) subvariety of its ambient space $\A^n$ or $\Pe^n$ and
call $V$ the affine or projective variety defined by the
polynomials $F_1\klk F_s$. We think the variety $V$ to be equipped
with the induced Zariski topology and its sheaf of rational
functions. The irreducible components of $V$ are defined with
respect to its Zariski topology over $k$ . We call $V$ irreducible
if $V$ contains a single irreducible component and equidimensional
if all its irreducible components have the same dimension. The
dimension $dim\,V$ of the variety $V$ is defined as the maximal
dimension of all its irreducible components. If $V$ is
equidimensional we define its (geometric) degree as the number of
points arising when we intersect $V$ with $dim\,V$ many generic
(affine) linear hyperplanes of its ambient space $\A^n$ or
$\Pe^n$. For an arbitrary closed variety $V$ with irreducible
components $\mathcal{C}_1\klk\mathcal{C}_t$ we define its degree
as $\deg V:=\deg \mathcal{C}_1\plp\deg\mathcal{C}_t$. With this
definition of degree the intersection of two closed subvarieties
$V$ and $W$ of the same ambient space satisfies the B\'ezout
inequality $$\deg V\cap W\le \deg V\deg W$$ (see \cite{Heintz83},
\cite{Fulton84}, \cite{Vogel84}).

We denote by $k[V]$ the affine or (graded) projective coordinate
ring of the variety $V$. If $V$ is irreducible we denote by $k(V)$
its field of rational functions. In case that $V$ is a closed
subvariety of the affine space $\A^n$ we consider the elements of
$k[V]$ as $\overline{k}$--valued functions mapping $V$ into
$\overline{k}$. The restrictions of the projections $X_1\klk X_n$
to $V$ generate the coordinate ring $k[V]$ over $k$ and are called
the coordinate functions of $V$. The data of $n$ coordinate
functions of $V$ fixes an embedding of $V$ into the affine space
$\A^n$. Morphisms between affine and projective varieties are
induced by polynomial maps between their ambient spaces which are
supposed to be homogeneous if the source and target variety is
projective.

Replacing the ground field $k$ by its algebraic closure
$\overline{k}$, we may apply all this terminology again. In this
sense we shall speak about the Zariski topologies and coordinate
rings over $\overline{k}$ and sheaves of $\overline{k}$--rational
functions. In this more general context varieties are defined by
polynomials with coefficients in $\overline{k}$. If we want to
stress that a particular variety $V$ is defined by polynomials
with coefficients in the ground field $k$, we shall say that $V$
is $k$--{\em definable} or $k$--{\em constructible}. The same
terminology is applied to any set determined by a (finite) boolean
combination of $k$--definable closed subvarieties of $\A^n$ or
$\Pe^n$. By a {\em constructible} set we mean simply a
$\overline{k}$--constructible one. Constructible and
$k$--constructible sets are always thought to be equipped with
their corresponding Zariski topology. In case of $k:=\Q$ and
$\overline{k}:=\C$ we shall sometimes also consider the {\em
euclidean} (i.e. ``{\em strong}") topology of $\A^n$ and $\Pe^n$
and their constructible subsets.
\medskip

The rest of our terminology and notation of algebraic
geometry and commutative algebra is standard and can be found in
\cite{Lang58}, \cite{Shafarevich84}, \cite[Chapter I]{Mumford88},
and in \cite{Lang93}, \cite{AtMa69}, \cite{Matsumura80}.

%-----------------------------------------------------------------
%-----------------------------------------------------------------
%-----------------------------------------------------------------
%-----------------------------------------------------------------

\subsection{Algorithmic models and complexity measures.}
\label{models}
The algorithmic problems we are going to consider
in this paper will depend on {\em continuous parameters} and
therefore the corresponding input data structures have to contain
entries for these parameters. We call them {\em problem} or {\em
input parameters}.

Once such a parametric problem is given, the specialization of the
parameters representing input objects are called (admissible) {\em
problem} or {\em input instances}. Thus the problem parameters may
in principle be algebraically dependent. An algorithm solving the
given problem operates on the corresponding input data structure
and produces for each admissible input instance an {\em output
instance} which belongs to a previously chosen output data
structure. We shall always require that output instances depend
{\em rationally} on the input parameters. Since we limit in this
paper our attention to branching--free algorithms, particular
admissible input instances may not produce well defined output
instances. In order to surmount this difficulty, we shall in the
sequel admit certain limit processes which we modelize using the
notion of places from valuation theory. These places will mimic
the process of limit determination and calculation by means of de
l'H\^opital's rule.

The chosen output data structure must enable us to answer certain
previously fixed questions about the output objects of our
algorithmic problem.

Let us consider the case that these output objects are polynomial
functions and that we wish to answer the value question for these
functions (we shall say that we want to ``compute" or ``evaluate"
them). For the sake of definiteness let us suppose that there is
given an algorithmic problem depending on $r$ parameters and that
this problem is expressible in the elementary language of
algebraically closed fields over the ground field $k$. Let
$U_1\klk U_r$ be indeterminates representing the input parameters
of the given problem. Let $S\subset\A^r$ be the Zariski closure of
the set of admissible input instances and suppose that $S$ is
irreducible. Since our algorithmic problem is elementarily
expressible over $k$, we conclude that $S$ is $k$--definable. Let
$m$ be the size of the output data structure we are going to use
for the solution of our problem. In the sense of this paper, a
{\em (branching--free) continuous algorithm} computing for each
admissible input instance the code of the corresponding output
object, is given by certain rational functions $\theta_1\klk
\theta_m$ of $k(S)$ such that the rational map
$\theta=(\theta_1\klk\theta_m)$ is well--defined for any
admissible input instance $u\in S$ and such that $\theta$ maps $u$
to the corresponding output instance $\theta(u)$ (observe that the
admissible input instances form a Zariski dense subset of $S$).
Suppose now that our output objects are polynomial functions in
the variables $Y_1\klk Y_t$. We call our algorithm {\em
essentially division--free} if these polynomial functions belong
to the polynomial ring $k[\theta_1\klk \theta_m][Y_1\klk Y_t]$.
Thus essentially division--free algorithms do not contain
divisions which involve any of the arguments $Y_1\klk Y_t$ of the
output objects. Nevertheless such an algorithm is allowed to
contain divisions involving exclusively elements of $k(U)$ which
represent rational functions of $k(S)$. Once the value $\theta(u)$
is determined for an admissible input instance $u\in S$, the
output objects may be evaluated in any point $y\in\A^t$ without
using additional divisions. If moreover the parameter functions
$\theta_1\klk\theta_m$ belong to the coordinate ring $k[S]$, we
shall say that our algorithm if {\em totally division--free}.
Unfortunately, the limitation to totally division--free algorithms
would be too restrictive for an appropriate complexity analysis of
geometric elimination problems. On the other hand, the notion of
essentially division--free algorithm modelizes in a fairly
realistic manner the intuitive meaning of algebraic (symbolic)
tools in situations which admit branching--free procedures. In
particular, it captures all today known parametric elimination
procedures for these situations.\medskip

Suppose now that our algorithmic problem is well defined for
any element of $S$. Thus the set of admissible input instances is
the Zariski closed set $S$. Suppose furthermore that there are
given rational functions $\theta_1\klk \theta_m\in k(U)$ and a
constructible Zariski dense subset $S_0$ of $S$, such that the
rational map $\theta=(\theta_1\klk \theta_m)$ is defined in any
point of $S_0$ and such that $\theta$ represents an essentially
division--free algorithm which solves our algorithmic problem for
any input instance belonging to $S_0$ correctly. We shall say that
the given algorithm can be ({\em uniquely}) {\em extended to the
limit data structure $S$ of $S_0$} (and to the corresponding limit
input objects) if the following condition is satisfied:
\smallskip

\noindent {\em for any input instance $u\in S$ and any place
$\varphi:\overline{k}(S)\to\overline{k}\cup\{\infty\}$ whose
valuation ring contains the local ring of the variety $S$ at the
point $u$, the values $\varphi(\theta_1)\klk\varphi(\theta_m)$
are} finite {\em and} uniquely determined {\em by the input
instance $u$}.
\smallskip

\noindent We observe that this condition implies that
$\theta_1\klk\theta_m$ belong to the integral closure of $k[S]$ in
$k(S)$.

Intuitively speaking, we admit certain (algebraic) limit processes
in the spirit of de l'H\^opital's rule in order to extend the
given algorithms from $S_0$ to the limit data structure $S$. These
limit processes are necessary because in elimination theory one
often faces situations where parameters become algebraically
dependent elements of domains which are not factorial. Greatest
common divisor computations for polynomials with coefficients in
these domains lead then to essential divisions of elements of
these domains (i.e. to divisions whose results do not anymore
belong to the given domain). These kind of situations can be found
in \cite{GiHe01} and Section \ref{section5.3}.\medskip

In the context of this paper we shall not care about the
representation of the rational map $\theta$. However, in
concrete situations, it is reasonable to think that the
rational functions $\om \theta m$ are represented by numerator and
denominator polynomials belonging to $k[\om U r]$, and that these
polynomials are holomorphically encoded by a suitable data
structure (see Section \ref{holocont} for the notion of
holomorphic encoding).\medskip

Let us finally exemplify the abstract notion of an
essentially division--free algorithm in the context of arithmetic
circuits (see \cite{BuClSh97} for details).

An {\em essentially division--free arithmetic circuit} is an
algorithmic device that can be represented by a labeled directed
acyclic graph ({\em dag}) as follows:

\noindent the circuit depends on certain input nodes, labeled by
indeterminates over the ground field $k$. These indeterminates are
thought to be subdivided in two disjoints sets, representing the
{\em parameters} and the {\em variables} of the given circuit. For
the sake of definiteness, let $U_1\klk U_r$ be the parameters and
$Y_1\klk Y_t$ the variables of the circuit. Let $K:=k(U_1\klk
U_r)$. We call $K$ the {\em parameter field} of the circuit. The
circuit nodes of indegree zero which are not inputs are labeled by
elements of $k$, which are called the {\em scalars} of the circuit
(here ``indegree" means the number of incoming edges of the
corresponding node). Internal nodes are labeled by arithmetic
operations (addition, subtraction, multiplication and division).
We require that the internal nodes of the circuit represent
polynomials in the variables $Y_1\klk Y_t$. We call these
polynomials the intermediate results of the given circuit. The
coefficients of these polynomials belong to the parameter field
$K$. In order to achieve this requirement, we allow in an
essentially division-free circuit only divisions which involve
elements of $K$. Thus essentially division--free circuits do not
contain divisions involving intermediate results which depend on
the variables $Y_1\klk Y_t$. A circuit which contains only
divisions by nonzero elements of $k$ is called {\em totally
division-free}.

Finally we suppose that the given circuit contains one or more
nodes which are labeled as output nodes. The results of these
nodes are called {\em outputs} of the circuit. Output nodes may
occur labeled additionally by sign marks of the form ``$=0$" or
``$\not=0$" or may remain unlabeled. Thus the given circuit
represents by means of the output nodes which are
labeled by sign marks a system of parametric polynomial equations
and inequations. This system determines in its turn for each
admissible parameter instance a locally closed set (i.e. an
embedded affine variety) with respect to the Zariski topology of
the affine space $\A^t$ of variable instances. The output nodes of
the given circuit which remain unlabeled by sign marks represent a
parametric polynomial application (in fact a morphism of algebraic
varieties) which maps for each admissible parameter instance the
corresponding locally closed set into a suitable affine space. We
shall interpret the system of polynomial equations and inequations
represented by the circuit as a {\em parametric family of systems}
in the variables of the circuit. The corresponding varieties
constitute a {\em parametric family of varieties}. The same point
of view is applied to the morphism determined by the unlabeled
output nodes of the circuit. We shall consider this morphism as a
{\em parametric family of morphisms}.\medskip

To a given essentially division--free arithmetic circuit we
may associate different complexity measures and models. In this
paper we shall be exclusively concerned with {\em sequential}
computing {\em time}, measured by the {\em size} of the circuit.
Our main complexity model is the non--scalar one, over the
parameter field $K$. Exceptionally we will also consider the
non--scalar complexity model over the ground field $k$. In the
non--scalar complexity model over $K$ we count only the {\em
essential} multiplications (i.e. multiplications between
intermediate results which actually involve variables and not
exclusively parameters). This means that $K$--linear operations
(i.e. additions and multiplications by arbitrary elements of $K$)
are {\em cost free}. Similarly, $k$--linear operations are not
counted in the non-scalar model over $k$.\medskip

Let $\theta_1\klk\theta_m$ be the elements of the parameter
field $K$ computed by the given circuit. Since this circuit is
essentially division--free we conclude that its outputs belong to
$k[\theta_1\klk\theta_m][Y_1\klk Y_t]$. Let $L$ be the non--scalar
size (over $K$) of the given circuit and suppose that  the circuit
contains $q$ output nodes. Then the circuit may be rearranged
(without affecting its non--scalar complexity nor its outputs) in
such a way that the condition
\begin{equation}\label{arith_circuits} m=L^2+(2t-1)L+q(L+t+1)
\end{equation}
is satisfied (see \cite[Chapter 9, Exercise 9.18]{BuClSh97}). In
the sequel we shall always assume that we have already performed
this rearrangement. Let \linebreak $Y:=(Y_1\klk Y_t)$,
$\theta:=(\theta_1\klk\theta_m)$ and let $f_1\klk f_q\in
k[\theta][Y]$ be the outputs of the given circuit. Let $Z_1\klk
Z_m$ be new indeterminates and write $Z:=(Z_1\klk Z_m)$. Then
there exist polynomials $F_1\klk F_q\in k[Z,Y]$ such that
$f_1=F_1(\theta,Y)\klk f_q=F_q(\theta,Y)$ holds. Let us write
$f:=(f_1\klk f_q)$ and $F:=(F_1\klk F_q)$. Consider the object
class $$\mathcal{O}:=\{F(\zeta,Y):\zeta\in\A^m\}$$ which we think
represented by the data structure $\mathcal{D}:=\A^m$ by means of
the obvious encoding which maps each code $\zeta\in\mathcal{D}$ to
the object $F(\zeta,Y)\in\overline{k}[Y]^q$.

For the moment, let us consider as input data structure the
Zariski open subset $\mathcal{U}$ where the rational map
$\theta=(\theta_1\klk\theta_m)$ is defined. Then the given
essentially division--free arithmetic circuit represents an
algorithm which computes for each input code $u\in\mathcal{U}$ an
output code $\theta(u)$ representing the output object
$f(u,Y)=F\big(\theta(u),Y\big)$. This algorithm is in the above
sense essentially division--free. From identity
(\ref{arith_circuits}) we deduce that the size $m$ of the data
structure $\mathcal{D}$ is closely related to the non--scalar size
$L$ of the given circuit. In particular we have the estimate
\begin{equation}\label{arith_circuits2} \sqrt{m}-(t+q)\le L.\end{equation}
Later we shall meet specific situations where we are able to
deduce from a previous (mathematical) knowledge of the
mathematical object \linebreak $f=(f_1\klk f_q)$ a lower bound for
the size of the output data structure of {\em any}
essentially division--free algorithm which computes for an
arbitrary input code $u\in\mathcal{U}$ the object $f(u,Y)$. Of
course, in such situations we obtain by means of
(\ref{arith_circuits2}) a lower bound for the non--scalar size
(over $K$) of {\em any} essentially division--free
arithmetic circuit which solves the same task. In particular we
obtain lower bounds for the total size and for the non--scalar
size over $k$ of all such arithmetic circuits.

%-----------------------------------------------------------------
%-----------------------------------------------------------------
%-----------------------------------------------------------------
%-----------------------------------------------------------------
%-----------------------------------------------------------------
%-----------------------------------------------------------------
%-----------------------------------------------------------------
%-----------------------------------------------------------------

\section{Holomorphic, continuous and robust encodings.}
\label{section3}

\subsection{Holomorphic and continuous encodings.}
\label{holocont}

Let $\mathcal{O}$ be an object class of polynomial functions
belonging to the polynomial ring $\overline{k}[Y_1\klk Y_t]$. We
shall say that $\mathcal{O}$ is {\em $k$--constructible} (or
$k$--{\em definable}) if the following conditions are satisfied:
\begin{itemize}
  \item[($i$)] The $\overline{k}$--vector space $W$ generated by
  the elements of $\mathcal{O}$ in $\overline{k}[Y_1\klk Y_t]$ is
  finite dimensional and there exists a $\overline{k}$--basis of
  $W$ consisting of polynomials which belong to $k[Y_1\klk Y_t]$
  (we call such a basis of $W$ {\em canonical}).
  \item[($ii$)] With respect to a given canonical basis of $W$, the
  object class $\mathcal{O}$ forms a $k$--constructible subset of
  $W$ (observe that this condition does not depend on the particular
  canonical basis we have chosen).
\end{itemize}

Suppose now that the object class $\mathcal{O}$ is
$k$--constructible and fix a canonical basis $P=(P_1\klk P_{N'})$
of $W$. Without loss of generality we may assume $P_1\klk
P_{N'}\in\mathcal{O}$. The evaluation map
$eval:W\times\A^t\to\A^1$ is defined by $eval(F,y):=F(y)$ for
$F\in W$ and $y\in\A^t$. With respect to the canonical basis $P$,
the evaluation map is $k$--definable and linear in its first
argument. Since $P_1\klk P_{N'}$ are polynomials of $k[Y_1\klk
Y_t]$ one sees easily that there exists a bound $\Delta\in\N$ with
$\deg F\le\Delta$ for any $F\in W$. Let
$N\ge\begin{pmatrix}\Delta+t\\t\end{pmatrix}$. Then we have $N'\le
N$ and there exist suitable (generic interpolation) points
$\eta_1\klk\eta_N\in k^t$ such that the map $\varphi:W\to\A^N$
defined for $F\in W$ by $\varphi(F):=\big(eval(F,\eta_1)\klk
eval(F,\eta_N)\big)=\big(F(\eta_1)\klk F(\eta_N)\big)$ induces a
$\overline{k}$--linear embedding of $W$ into the affine space
$\A^N$. Observe that $\varphi$ is $k$--definable with respect to
the canonical basis $P$ of $W$. In particular the image of
$\varphi$ is a $k$--definable linear subspace of $\A^N$ of
dimension $N'$. Under the embedding $\varphi$, the object class
$\mathcal{O}$ becomes a $k$--constructible subset of the {\em
ambient space} $\A^N$ and the evaluation map becomes a
$k$--definable morphism of algebraic varieties which is linear in
its first argument and whose domain of definition can be extended
(not uniquely) to the affine space $\A^N\times\A^t$.\medskip

This is the point of view we shall adopt in the sequel for
$k$--constructible object classes of polynomial functions.

In particular we consider $\mathcal{O}$ and $W$ as topological
spaces equipped with the Zariski (or, in case $k:=\Q$ and
$\overline{k}:=\C$, with the strong topology) induced from the
ambient space $\A^N$. Observe that the Zariski closure
$\overline{\mathcal{O}}$ of the object class $\mathcal{O}$ is a
$k$--definable closed subvariety of $\A^N$ whose degree does not
depend on the particular $\overline{k}$--linear embedding
$\varphi$ we have chosen. We denote this degree by $\deg
\overline{\mathcal{O}}$. Furthermore observe that any upper bound
for the degree of the polynomials of $\overline{k}[Y_1\klk Y_t]$
contained in $\mathcal{O}$ is also an upper bound for the degree
of the polynomials in $\overline{\mathcal{O}}$.

We say that $\mathcal{O}$ is a {\em cone} if for any
$\lambda\in\overline{k}$ the set $\lambda\mathcal{O}:=\{\lambda
f;f\in\mathcal{O} \}$ is contained in $\mathcal{O}$. Suppose that
$\mathcal{O}$ is a cone. One immediately verifies that the
$\overline{k}$--closure $\overline{\mathcal{O}}$ of $\mathcal{O}$
is a $k$--definable cone which is contained in $W$. Therefore the
evaluation map $eval:W\times\A^t\to\A^1$ induces a $k$--definable
morphism of algebraic varieties
$\overline{\mathcal{O}}\times\A^t\to\A^1$ which we denote also by
$eval$, which is homogeneous of degree one in its first argument
(i.e. for $f\in\mathcal{O}$, $y\in\A^t$ and
$\lambda\in\overline{k}$ we have $eval(\lambda f,y)=\lambda\,
eval(f,y)$).\medskip

Let $\mathcal{O}$ be an arbitrary (not necessarily
$k$--constructible) object class  of polynomial functions
belonging to the polynomial ring $\overline{k}[Y_1\klk Y_t]$. Let
$\gamma_1\klk \gamma_m\in\A^t$ and let
$\gamma:=(\gamma_1\klk\gamma_m)$. We say that $m$ is the length
of $\gamma$. We call $\gamma$ a {\em correct test sequence} for
the object class $\mathcal{O}$ if for any polynomial
$F\in\mathcal{O}$ the following implication holds:
$$F(\gamma_1)=\cdots=F(\gamma_m)=0\Rightarrow F=0.$$

We call $\gamma$ an {\em identification sequence} for
$\mathcal{O}$ if for any two polynomials $F_1,F_2\in\mathcal{O}$
the following implication holds $$F_1(\gamma_1)=F_2(\gamma_1)\klk
F_1(\gamma_m)=F_2(\gamma_m)\Rightarrow F_1=F_2.$$

\medskip

Now we suppose that there exists a bound $\Delta\in\N$ with $\deg
F\le\Delta$ for any $F\in\mathcal{O}$. Let
$N\ge\begin{pmatrix}\Delta+t\\t\end{pmatrix}$. We may interpret
$\mathcal{O}$ as a subset of $\A^N$. Suppose now that there is
given a $k$--definable data structure $\mathcal{D}\subset\A^L$
which encodes the object class $\mathcal{O}$ and contains a {\em
Zariski--dense set of $k$--rational points}. Let
$\omega:\mathcal{D}\to\mathcal{O}$ be this encoding and suppose
that there exists a $k$--definable polynomial map
$\rho:\A^L\times\A^t\to\A^1$ with $\rho(D,y)=\omega(D)(y)$ for any
$D\in\mathcal{D}$ and any $y\in\A^t$. In these circumstances
we say that $\rho$ allows
to answer the value question about the object class
$\mathcal{O}$
{\em holomorphically}.

\begin{remark} \label{holo}
Let assumptions and notations be as before. Then $\mathcal{O}$ is
$k$--constructible and $\omega: \mathcal{D}\to \mathcal{O}$ is the
restriction of a suitable $k$--definable polynomial map $\Omega:
\A^L \to \A^N$.
\end{remark}

\begin{prf}
Let $N'\le N$ be the dimension of the $\overline{k}$--vector
space $W$ generated in $\overline{k}[Y_1\klk Y_t]$ by the
elements of $\mathcal{O}$. Choose $N$ generic interpolation
points $\eta_1\klk\eta_{N}\in k^t$ for the polynomials of
$\overline{k}[Y_1\klk Y_t]$ of degree at most $\Delta$. Since
the $k$--rational points are Zariski--dense in $\mathcal{D}$ we
conclude that there exist $D_1\klk D_{N'}\in\mathcal{D}\cap k^L$
such that for any choice of indices $1\le k_1<\cdots<k_{N'}\le N$
the $N'\times N'$--matrix
$$\big(\omega(D_i)(\eta_{k_j})\big)_{1\le i,j\le
N'}=\big(\rho(D_i,\eta_{k_j})\big)_{1\le i,j\le N'}$$ is regular.
Since for any such index choice this matrix is $k$--rational we
deduce that $\omega(D_1)\klk\omega(D_{N'})$ are polynomials which
belong to $\mathcal{O}\cap k[Y_1\klk Y_t]$ and form a basis of the
$\overline{k}$--vector space $W$. In the same manner as before,
using the $k$--rational interpolation points $\eta_1\klk\eta_N$,
we may construct from $\rho$ a $k$--definable polynomial map
$\Omega:\A^L\to\A^N$ with $\Omega|_{\mathcal{D}}=\omega$ (here
$\Omega|_{\mathcal{D}}$ denotes the restriction of the map
$\Omega$ to the set $\mathcal{D}$). In particular we have
$\Omega(\mathcal{D})=\mathcal{O}$. Since the polynomial map
$\Omega$ is $k$--definable we conclude that $\mathcal{O}$ is a
$k$--constructible subset of $\A^N$.
\end{prf} \cqfd

The preceding considerations about object classes of
polynomial functions and their encodings lead us to the following
fundamental notions of this paper:
\begin{definition} \label{holodef}
Let be given a data structure $\mathcal{D}\subset \A^L$, an
object class $\mathcal{O}\subset \A^N$ and an encoding $\omega :
\mathcal{D}\to \mathcal{O}$.

We call $\omega$ a $k$--definable encoding if the graph of
$\omega$ is a $k$--constructible subset of the affine space
$\A^L\times \A^N$.

Similarly the data structure $\mathcal{D}$ and the object class
$\mathcal{O}$ are called $k$--constructible (or $k$--definable)
if they form $k$--constructible subsets of the affine spaces
$\A^L$ and $\A^N$ respectively.

Let $\omega : \mathcal{D}\to \mathcal{O}$ be $k$--definable. We
call $\omega$ a continuous encoding if $\omega$ is a continuous
map with respect to the Zariski topologies of $\mathcal{D}$ and
$\mathcal{O}$ (or, in case $k:=\mathbb{Q}$ and
$\overline{k}:=\mathbb{C}$, with respect to their strong
topologies).

We call $\omega$ a holomorphic encoding if there exists a
$k$--definable polynomial map $\Omega: \A^L\to \A^N$ with $\omega
=\Omega\mid_{\mathcal{D}}$.
\end{definition}

In case that the $k$--definable encoding $\omega : \mathcal{D}\to
\mathcal{O}$ is holomorphic, we observe that $\omega$ can be
extended uniquely to a morphism of algebraic varieties mapping
$\overline{\mathcal{D}}$ into $\overline{\mathcal{O}}$. We denote
this morphism also by $\omega$.

%-----------------------------------------------------------------
%-----------------------------------------------------------------
%-----------------------------------------------------------------
%-----------------------------------------------------------------

\subsection{Robust encodings.}
\label{subsection:robust_encodings} For data structures, object
classes and encodings which are defined over $\mathbb{Q}$ and
interpreted over $\mathbb{C}$, it may happen that an
unbounded sequence of codes produces a convergent sequence of
objects. This is for example a typical behaviour of circuit
encodings of polynomials (see Sections \ref{section4real} and
\ref{section4.2}).

A {\em continuous} encoding which does not admit this phenomenon
is called {\em robust}. Unfortunately, this notion of robustness
is only well defined in case $k:=\mathbb{Q}$ and
$\overline{k}:=\mathbb{C}$. In order to obtain a more operative
notion of robustness which is also applicable to ground fields of
arbitrary characteristic, we are going to analyze this notion of
robustness under the restriction that the given encoding is not
only continuous, but also {\em holomorphic}. This will lead us to
a new definition of robustness which is equivalent to the
previous one in case $k:=\mathbb{Q}$, $\overline{k}:=\mathbb{C}$
and in case that the given encoding is
holomorphic.\medskip

Let $\md\subset\A^L(\C)$ and $\mo\subset\A^N(\C)$ be
$\Q$--constructible sets and let $\omega:\omd\to\omo$ be a
$\Q$--definable map with $\omega(\md)=\mo$. Suppose that $\omega$
is continuous with respect to the strong topologies of
$\overline{\mathcal{D}}$ and $\overline{\mathcal{O}}$. Let us
consider $\md$ as a data structure, $\mo$ as an object class and
$\omega|_{\md}: \md\to\mo$ as a $\Q$--definable {\em continuous}
encoding of the object class $\mo$ by the data structure $\md$.
In order to simplify notations we denote the map $\omega|_{\md}$
just by $\omega:\md\to\mo$.

\begin{definition} \label{robustcont}
(Robustness for continuous encodings)

\noindent Let notations and assumptions be as before. We call the
continuous encoding $\omega:\md\to\mo$ robust if $\omega$
satisfies the following condition: \smallskip

\noindent let $(D_i)_{i\in\N}$ be an arbitrary sequence of
elements of $\mathcal{D}$ encoding a sequence $(O_i)_{i\in\N}$ of
objects of $\mathcal{O}$. Let $O\in\mo$ be an accumulation point
of $(O_i)_{i\in\N}$ (with respect to the strong topology of
$\A^N(\C)$). Then there exists in $\overline{\mathcal{D}}$ an
accumulation point $Q$ of the sequence $(D_i)_{i\in\N}$ with
$\omega(Q)=O$.
\end{definition}

\begin{remark} \label{remark:proper}
Let notations and assumptions be as before. Suppose
furthermore that the data structure $\mathcal{D}$ is a closed
subvariety of its affine ambient space. Then the robustness
of the encoding $\omega:\md\to\mo$ is equivalent to the condition
that $\omega$ is a surjective and $\Q$--definable proper
continuous map of topological spaces. If $\omega$ is robust, then
$\omega$ has finite, non-empty fibers. Moreover, if
$\mathcal{O}$ is a closed subvariety of its affine ambient
space, then $\omega$--preimages of compact subsets of
$\mathcal{O}$ are compact.
\end{remark}

\begin{prf}
One sees easily that properness of the $\Q$--definable, continuous
map $\omega$ implies its robustness.

Suppose now that $\omega$ is robust. Since $\mathcal{D}$ is
closed, we conclude that $\omega$ is a closed, continuous map with
(sequentially) compact fibers. Hence $\omega$ is
proper.

From the arguments used at the beginning of the proof of Lemma
\ref{lemma:robust} below, one deduces easily that $\omega$ has
finite fibers.

Suppose furthermore that $\mathcal{O}$ is a closed subvariety
of its affine ambient space. Then $\mathcal{D}$ and
$\mathcal{O}$ are locally compact topological spaces. Therefore,
since $\omega$ is a proper continuous map, we conclude that
$\omega$--preimages of compact subsets of ${\mathcal{O}}$ are
compact.
\end{prf}\cqfd

We are now going to discuss the notion of robustness in terms of
algebraic geometry in order to obtain a suitable and well
motivated definition of robustness for $k$--definable {\em
holomorphic} encodings over any ground field $k$. For the rest of
this
 subsection we assume that $\omega$ is a holomorphic encoding.
\medskip

We are going to use the following fact:

\begin{lemma}
\label{fact} Let $S$ be a locally closed subvariety of $\A^n(\C)$.
Suppose $dim\, S>0$. Then $S$ is unbounded in $\A^n(\C)$.
\end{lemma}

\begin{prf}
Let $r:=dim\, S$. From Noether's Normalization Lemma we deduce
that there exists a linear map $\varphi:\A^n(\C)\to\A^r(\C)$ with
$\varphi(\overline{S})=\A^r(\C)$ (see \cite[I.7]{Mumford88}).
Observe that the $\C$--Zariski closure $\overline{S}$ of $S$
coincides with the closure of $S$ in the strong topology of
$\A^n(\C)$ (see \cite[I.10, Corollary 1]{Mumford88}). Suppose that
$S$ is bounded. Then $\overline{S}$ is compact in the strong
topology and therefore also its image $\varphi(\overline{S})$.
However $\varphi(\overline{S})=\A^r(\C)$ is not  compact since $r$
is positive.
\end{prf}\cqfd

The following key result will lead us to the intended notion of
robustness for holomorphic encodings defined over ground fields
of arbitrary characteristic:

\begin{lemma}
\label{lemma:robust} Let notations and assumptions be as before
and suppose that $\omega:\mathcal{D}\to\mathcal{O}$ is a
holomorphic encoding. Suppose that $\omega$ is robust in the
sense of Definition \ref{robustcont}. Let $V$ be a closed
irreducible ($\C$--definable) subvariety of $\A^L(\C)$ and
suppose that there exists a nonempty Zariski open subset
$\mathcal{U}$ of $V$ such that $\mathcal{U}$ is contained in
$\md$. Let $W:=\overline{\omega(\mathcal{U})}$ and let $O$ be a
point of $\omega(\mathcal{U})$. Let $\mathfrak{m}$ be the maximal
ideal of $\C[W]$ which defines the point $O$. Then
$\C[V]_{\mathfrak{m}}$ is a finite $\C[W]_{\mathfrak{m}}$--module
(i.e. the $\C$--algebra extension
$\C[W]_{\mathfrak{m}}\to\C[V]_{\mathfrak{m}}$ induced by $\omega$
is integral).
\end{lemma}

\begin{prf}
Since $\mathcal{U}$ Zariski dense in $V$ and $\mathcal{U}$ is
contained in $\mathcal D$, we conclude that $V$ is
contained in $\overline{\mathcal{D}}$.

By assumption $\omega:\md\to\mo$ is a $\Q$--definable
morphism of algebraic varieties. Therefore there exists a
(unique) extension of $\omega|_{\mathcal{U}}:\mathcal{U}\to\mo$ to
a morphism of algebraic varieties which maps $V$ into $W$. We
denote this morphism by $\omega:V\to W$ and observe that it is a
dominant morphism of irreducible affine varieties. Thus
$\omega:V\to W$ induces an injective \linebreak $\C$--algebra
homomorphism $\C[W]\to\C[V]$ and consequently a field extension
$\C(W)\to\C(V)$. Observe that $\omega^{-1}(O)\cap \mathcal{U}$ is
a locally closed algebraic subvariety of $\A^L(\C)$. From
$O\in\omega(\mau)$ we deduce that $\omega^{-1}(O)\cap\mau$ is not
empty. Therefore $r:=dim\, \omega^{-1}(O)\cap\mau$ is nonnegative.
Suppose $r>0$. Then from Lemma \ref{fact} we deduce that
$\omega^{-1}(O)\cap\mau$ is unbounded. Thus there exists a
sequence $(D_i)_{i\in\N}$ of points of
$\omega^{-1}(\mathcal{O})\cap\mau \subset\overline{\md}$ which has
no accumulation point. On the other hand we have $O=\omega(D_i)$
for any $i\in\N$. Therefore $\big(\omega(D_i)\big)_{i\in\N}$ is a
sequence of elements of the object class $\mo$ which converges to
the point $O$. Since $\omega$ is by assumption a robust encoding
in the sense of Definition \ref{robustcont}, we conclude that
$(D_i)_{i\in\N}$ must contain an accumulation point in
$\overline{\mathcal{D}}$. This contradicts the choice of the
 sequence $(D_i)_{i\in\N}$.
Thus we conclude $r=0$.

From the Theorem of Fibers we deduce now that $dim\, W= dim\, V$
holds and that $\C(W)\hookrightarrow\C(V)$ is a finite field
extension. Following \cite[Chapter V]{Lang58} (see also
\cite[Chapter II, 5.2]{Shafarevich84}) we may choose a finite
morphism of irreducible affine varieties $\psi: \widetilde{W}\to
W$ such that the coordinate ring $\C[\widetilde{W}]$ is isomorphic
to the integral closure of $\C[W]$ in $\C[V]$. Observe  that there
exists a unique morphism of affine varieties
$\widetilde{\omega}:V\to\widetilde{W}$ such that the diagram
$$\xymatrix{ & &\widetilde{W}\ar[dd]^{\mbox{$\psi$}}\\
V\ar[urr]^{\mbox{$\widetilde{\omega}$}} \ar[drr]^{\mbox{$\omega$}}
& &
\\ & & W}$$
commutes. Since $\omega$ is dominant we conclude that
$\widetilde{\omega}$ is dominant too. Thus $\widetilde{\omega}$
induces an injective $\C$--algebra homomorphism
$\C[\widetilde{W}]\to\C[V]$ which maps $\C[\widetilde{W}]$ onto a
subring of $\C[V]$ which is integrally closed in $\C[V]$. In this
sense we shall say that {\em $\C[\widetilde{W}]$ is integrally
closed in $\C[V]$.}

Let now $P\in\widetilde{W}$ be an arbitrary  point with $\psi
(P)=O$ (observe that such a point exists since $\psi$ is
surjective). Since $\widetilde{\omega }$ is dominant, we conclude
that $\widetilde{\omega }(\mau)$ contains a nonempty Zariski open
subset of $\widetilde{W}$. Hence $\widetilde{\omega }(\mau)$ is
dense in the strong topology of $\widetilde{W}$. Therefore we may
choose a sequence $(D_i)_{i\in\N}$ of elements of $\mau\subset\md$
such that $\big(\widetilde{\omega }(D_i)\big)_{i\in\N}$ converges
in the strong topology of $\widetilde{W}$ to the point $P$. Thus
the sequence $\big(\omega(D_i)\big)_{i\in\N}$ is a sequence of
elements of the object class $\mo$ which converges to the object
$O\in\mo$. Since by assumption the encoding $\omega:\mathcal{D}\to
\mathcal{O}$ is robust in the sense of Definition
\ref{robustcont}, we conclude that there exists an accumulation
point $D\in\overline{\md}$ of the sequence $(D_i)_{i\in\N}$.
Without loss of generality we may assume that $(D_i)_{i\in\N}$
converges to $D$. This implies $D\in V$ and $\widetilde{\omega
}(D)=P$.

Thus the fiber $\widetilde{\omega}^{-1}(P)$ has nonnegative
dimension. Suppose now that the dimension of
$\widetilde{\omega}^{-1}(P)$ is positive. Then Lemma \ref{fact}
implies that $\widetilde{\omega }^{-1} (P)$ is unbounded.
Therefore we may choose a sequence $(Q_n)_{n\in\N}$ of points of
$\widetilde{\omega}^{-1}(P)$ which has no accumulation point.
Since $\mau$ is dense in the strong topology of $V$ there exists a
family $(D_i^{(n)})_{n,i\in\N}$ of elements of $\mau$ such that
for any $n\in\N$ the sequence $\iseq {D^{(n)}}$ converges to the
point $Q_n$. Without loss of generality we may suppose that this
convergence is uniform in the parameter $n$. Observe that we have
$\omega (Q_n)=\psi (P)=O$ for any index $n\in \N$. Therefore we
may assume without loss of generality that the
sequence $(\omega(D_n^{(n)}))_{n\in\N}$ converges to the object
$O$. From the robustness of $\omega $ we infer now that the
sequence $({D_n^{(n)}})_{n\in\N}$  has an accumulation point $Q$
in $\overline{\md}$. Since for any index $n\in\N$ the convergence
of the sequence $(D_i^{(n)})_{i\in\N}$ to $Q_n$ is uniform in $n$,
we conclude that $Q$ is an accumulation
 point of the sequence $(Q_n)_{n\in\N}$.
This contradicts the
choice of the sequence $\nseq Q$.

Therefore we have $dim\, \widetilde{\omega }^{-1}(P)=0$. Let
$\mathfrak{m}_P$ be the maximal ideal of $\C[\widetilde{W}]$ which
defines the point $P$ and consider $\C[V]$ as a
$\C[\widetilde{W}]$--module. Since $\C[\widetilde{W}]$ is
integrally closed in $\C[V]$, we deduce now from Zariski's Main
Theorem (see e.g. \cite[IV.2]{Iversen73}) that
\begin{equation}
\label{once}
\C[V]_{\mathfrak{m}_P}=\C[\widetilde{W}]_{\mathfrak{m}_P}\end{equation}
holds. Consider $\C[V]_\mathfrak{m}$ as a
$\C[\widetilde{W}]_\mathfrak{m}$--module and observe that the
maximal ideals of $\C[\widetilde{W}]_\mathfrak{m}$ correspond
bijectively to the maximal ideals of $\C[\widetilde{W}]$ of the
form $\mathfrak{m}_P$ with $P\in\widetilde{W}$ and $\psi(P)=O$.
From (\ref{once}) one deduces now
$\C[V]_\mathfrak{m}=\C[\widetilde{W}]_\mathfrak{m}$. Since by
definition of $\widetilde{W}$ the coordinate ring
$\C[\widetilde{W}]$ is a finite $\C[W]$--module, this implies that
$\C[V]_\mathfrak{m}$ is a finite $\C[W]_\mathfrak{m}$--module.
\end{prf}\cqfd

Let notations and assumptions be as before. From Lemma
\ref{lemma:robust} and its proof we infer that the encoding
$\omega $ satisfies the following conditions:
\begin{itemize}
\item[$(i)$] for any object $O\in\mo$ there are only finitely many
encodings $D\in\md$ with $\omega (D)=O$ (in this sense we shall
call the {\em ambiguity} of $\omega $ {\em finite}).
\item[$(ii)$] for any object $O\in\mo$ with maximal defining ideal
$\mathfrak{m}$ in $\C[\omo]$, the local ring
$\C[\omd]_\mathfrak{m}$ is a {\em finite}
$\C[\omo]_\mathfrak{m}$--module.
\end{itemize}

If the algebraic variety $\omd$ is irreducible then $\omo$ is
irreducible too and we may replace condition ($ii$) by the
following equivalent one:

\begin{itemize}
\item[($iii$)] Let $O$ be an arbitrary object of $\mo$.
Then any place $\varphi:\C(\omd)\to\C\cup\{\infty\}$ whose
valuation ring which contains the local ring of the variety $\omo$
at the point $O$, takes only finite values on $\C[\omd]$.
\end{itemize}

Although somewhat weaker and limited to the case
$\overline{\mathcal{D}}$ irreducible, condition ($iii$) is
the genuine algebraic--geometric counterpart of the notion of
robustness given by Definition \ref{robustcont}. This indicates
that the following definition of robustness for {\em holomorphic}
encodings (not simply continuous ones) captures the intuitive
meaning of the previous Definition \ref{robustcont} in case of a
ground field $k$ of arbitrary characteristic with arbitrary
algebraic closure $\overline{k}$.

\begin{definition}
\label{robust1} (Robustness of holomorphic encodings)

\noindent Let $\omega :\md\to\mo$ be a $k$--definable holomorphic encoding
of a $k$--constructible object class $\mo$ by a $k$--constructible
data structure $\md$. Then we call $\omega $ robust if
for any object $O\in\mo$ with maximal defining ideal
$\mathfrak{m}$ in $\overline{k}[\omo]$ the localization ring
$\overline{k}[\omd]_\mathfrak{m}$ is a finite
$\overline{k}[\omo]_\mathfrak{m}$--module.
\end{definition}

In case that a $k$--definable holomorphic encoding $\omega
:\md\to\mo$ induces a finite morphism of affine varieties which
maps $\omd$ onto $\omo$, we conclude that the encoding $\omega $
is robust in the sense of Definition \ref{robust1}.\medskip

We are now going to show that in case $k:=\Q$, $\ok :=
\C$ and $\omega$ holomorphic, Definition \ref{robustcont} and
Definition \ref{robust1} represent the same notion of robustness.

\begin{lemma} \label{equivrob}
Let $k:=\Q$, $\ok := \C$, and let $\mathcal{D}$ and
$\mathcal{O}$ be a $\Q$--constructible data structure and object
class respectively. Let $\omega : \mathcal{D} \to \mathcal{O}$ a
$\Q$--definable, holomorphic encoding. Then $\omega$ is
robust in the sense of Definition \ref{robustcont} (as a
continuous encoding with respect to the strong topologies of
$\mathcal{D}$ and $\mathcal{O}$) if and only if $\omega$ is
robust in the sense of Definition \ref{robust1} (as a holomorphic
encoding).
\end{lemma}

\begin{prf}
Suppose that $\omega$ is robust in the sense of Definition
\ref{robustcont}. Then, from the statement $(ii)$ above, we
deduce that $\omega$ is a robust encoding in the sense of
Definition \ref{robust1}.

Suppose now that $\omega$ is robust in the sense of Definition
\ref{robust1}. Let be given a sequence $(D_i)_{i\in \N}$ of
elements of the data structure $\mathcal{D}$ which encodes a
sequence $(O_i)_{i\in \N}$ of objects of $\mathcal O$. Let be
given an accumulation point $O\in \mathcal O$ of the sequence
$(O_i)_{i\in \N}$ with respect to the strong topology $\mathcal
O$. For the sake of simplicity we shall assume that $(O_i)_{i\in
\N}$ converges to $O$. Let $\mathfrak m$ be the maximal defining
ideal of $O$ in $\C[\overline{\mathcal{O}}]$.

Let us consider an arbitrary element $f$ of
$\C[\overline{\mathcal{O}}]$. In case that $f(O_i)=0$ holds for
infinitely many indices $i\in \N$, we conclude $f(O)=0$.
Therefore, if $f$ does not belong to the maximal ideal $\mathfrak
m$, then $f$ vanishes on all but finitely many entries of the
sequence $(O_i)_{i\in \N}$. Since $\omega$ is robust in the sense
of Definition \ref{robust1}, we may now conclude that there exists
an element $g$ of $\C[\overline{\mathcal{O}}]$ with the following
properties:

\begin{itemize}
\item[(a)] $g(O)\ne 0$ and $g(O_i)\ne 0$ for all but finitely
many indices $i\in \N$.
\item[(b)] $\C[\overline{\mathcal{D}}]_g$ is a finite
$\C[\overline{\mathcal{O}}]_g$--module.
\end{itemize}

For the sake of simplicity we shall suppose $g(O_i)\ne 0$ for any
$i\in\N$.

Consider now an arbitrary element $h$ of
$\C[\overline{\mathcal{D}}]$. Let $Y$ be an
indeterminate.

From properties (a) and (b) above we deduce that there exists a
monic polynomial $P\in \C[\overline{\mathcal{O}}]_g[Y]$ with
$P(h)=0$ and such that $P$ can be specialized for the object
$O$ and any index $i\in \N$ into well--defined
elements $P(O)$ and $P(O_i)$ of the polynomial ring $\C[Y]$.
Without loss of generality we may suppose that $P$ is the minimal
polynomial of $h$ over $\C(\overline{\mathcal{O}})$. Thus for any
$i\in \N$ we have $P(O_i)(h(D_i))=0$. Therefore the sequence
$(h(D_i))_{i\in \N}$ has an accumulation point which is a zero of
the polynomial $P(O)(Y)\in \C[Y]$. Since $P$ is the minimal
polynomial of $h$ over $\C(\overline{\mathcal{O}})$ we deduce from
property (b) above that there exists an element $Q\in
\overline{\mathcal{D}}$ with $g(Q)\ne 0$ and $\omega(Q)=O$ such
that $h(Q)$ is an accumulation point of the sequence
$(h(D_i))_{i\in \N}$.

Generalizing this argument to a finite set of generators of the
$\C[\overline{\mathcal{O}}]_g$--module
$\C[\overline{\mathcal{D}}]_g$ we conclude that the sequence
$(D_i)_{i\in \N}$ has an accumulation point in
$\overline{\mathcal{D}}$. Therefore the encoding $\omega$
is robust in the sense of Definition \ref{robustcont}.
\end{prf} \cqfd

\begin{remark} \label{remark:proper1}

Let $k:=\Q$ and $\ok =\C$. Then, in terms of algebraic geometry,
Lemma \ref{lemma:robust} and Remark \ref{remark:proper} imply the
following folkloric statement:
\smallskip

\noindent let $V$ and $W$ be closed, equidimensional subvarieties of
suitable complex affine spaces and let $\varphi:V\to W$ be a
morphism of affine varieties mapping $V$ onto $W$. Suppose that
$\varphi$ is a proper continuous map with respect to the strong
topologies of $V$ and $W$. Then $\varphi$ is a finite morphism of
affine varieties.
\end{remark}

%-----------------------------------------------------------------
%-----------------------------------------------------------------
%-----------------------------------------------------------------
%-----------------------------------------------------------------

\subsection{Correct test and identification sequences.}
\label{correct}
\subsubsection{Correct test and identification sequences for
holomorphic encodings.} \label{holocorrect}

We are now going to develop the fundamental technical tools we
shall need in Section \ref{valencode} for the
formulation and proof of the first main result of this paper,
namely Theorem \ref{theorem:encode1}.
\medskip

The following statement generalizes \cite[Theorem 4.4]{HeSc82}.
\begin{lemma}
\label{ts} Let $\mathcal{O}$ be a $k$--constructible object class
of polynomial functions belonging to $\overline{k}[Y_1\klk Y_t]$.
Let $\Delta\in\N$ be an upper bound for the degree of the
polynomials contained in $\mathcal{O}$. Suppose that there is
given a $k$--constructible data structure $\mathcal{D}\subset\A^L$
and a $k$--definable holomorphic encoding $\omega:\mathcal{D}\to
\mathcal{O}$. Suppose that there exists a quantifier--free
first--order formula which defines the data structure
$\mathcal{D}$ and whose equations involve only $K$ distinct
polynomials of degree at most $\Delta_1$ in $L$ indeterminates
over $k$. Moreover, assume that the encoding $\omega$ is definable
by polynomials of degree at most $\Delta_2\ge 1$ in $L$
indeterminates over $k$. Then the degrees of the algebraic
varieties $\overline{\mathcal{D}}$ and $\overline{\mathcal{O}}$
satisfy the estimates
$$\deg\overline{\mathcal{D}}\le(1+K\Delta_1)^L$$ and
$$\deg\overline{\mathcal{O}}\le(L+1)\Delta_2^L\deg\overline{\mathcal{D}}
\le(L+1)\big((1+K\Delta_1)\Delta_2\big)^L.$$ For
$\overline{\mathcal{O}}$ equidimensional this estimate may be
improved to
$$\deg\overline{\mathcal{O}}\le\Delta_2^L\deg\overline{\mathcal{D}}\le
\big((1+K\Delta_1)\Delta_2\big)^L.$$ Let $M$ be a finite subset of
$k$ having at least two elements. Suppose \linebreak $\#
M\ge\Delta^2(\deg\overline{\mathcal{O}})^{\frac{1}{L}}$ (observe
that this is the case if $$\#
M\ge\Delta^2(1+L)^{\frac{1}{L}}(1+K\Delta_1)\Delta_2$$ holds). Let
$m\ge 2L+2$. Then there exist points $\gamma_1\klk\gamma_m$ of
$M^t$ such that $\gamma:=(\gamma_1\klk\gamma_m)$ is a correct test
sequence for the object class $\overline{\mathcal{O}}$ (and hence
for $\mathcal{O}$).

Suppose that the points of the finite set $M^t$ are
equidistributed. Then the probability of finding in $M^{mt}$ by a
random choice such a correct test sequence is at least
$1-\frac{1}{\# M}\ge\frac{1}{2}$.
\end{lemma}

\begin{prf}
The proof is subdivided in three parts. Let us start with the first
one. With the terminology introduced before, suppose that the
object class $\mathcal{O}$ is given as a $k$--constructible subset
of some affine space $\A^N$. Let $Z_1\klk Z_L$ be the coordinate
functions of the affine space $\A^L$. By hypothesis there exists
a quantifier--free definition of $\mathcal{D}$ whose equations
involve only $K$ distinct polynomials $G_1\klk G_K\in k[Z_1\klk Z_L]$
of degree at most $\Delta_1$. Observe that for any irreducible component
$\mathcal{C}$ of $\overline{\mathcal{D}}$ there exists a subset
$\mathcal{G}$ of $\{G_1\klk G_K\}$ such that $\mathcal{C}$ is an
irreducible component of the closed subvariety $\{G=0;G\in
\mathcal{G}\}$ of $\A^L$. From \cite[Theorem 2]{JeSa00} (see also
\cite[Corollary 1]{Heintz83}) one deduces now easily the estimate
$$\deg\overline{\mathcal{D}} \le\sum_{h=0}^L
\begin{pmatrix}K\\h\end{pmatrix}\Delta_1^h
\le(1+K\Delta_1)^L.$$

By assumption there exist polynomials $\Omega_1\klk\Omega_N\in
k[Z_1\klk Z_L]$ of degree at most $\Delta_2$ such that
$\Omega:=(\Omega_1\klk \Omega_N)$ defines a polynomial map
$\Omega:\A^L\to\A^N$ with $\Omega|_{\mathcal{D}}=\omega$. Observe
that $\Omega$ induces a morphism of (possibly reducible) affine
varieties $\overline{\mathcal{D}}\to\overline{\mathcal{O}}$ which
we denote also by $\omega$. From $\omega(\mathcal{D})=\mathcal{O}$
we deduce that $\omega$ is dominant. This implies
$dim\,\overline{\mathcal{O}} \le dim\,\overline{\mathcal{D}}\le
L$.

Let $0\le h\le L$ and let $E_h$ be the union of the irreducible
components of $\overline{\mathcal{O}}$ of dimension $h$. Suppose
that $E_h$ is nonempty. Let $T_1\klk T_N$ be the coordinate
functions of $\A^N$. Since the morphism $\omega$ is dominant, we
may choose a nonempty, Zariski open subset $\mathcal{U}$ of $E_h$
which is contained in the image $\omega(\overline{\mathcal{D}})$
(see e.g. \cite[I.8, Theorem 3]{Mumford88}). On the other hand, we
may choose $N-h$ generic affine--linear equations $H_1\klk
H_{n-h}\in k[T_1\klk T_N]$ such that $E_h\cap\{H_1=0\klk
H_{N-h}=0\}$ consists of $\deg E_h$ points, all contained in
$\mathcal{U}$ and therefore in $\omega(\overline{\mathcal{D}})$
(see \cite{Heintz83}, Remark 2). Each of these points is the image
of a $\overline{k}$--irreducible component of the closed
subvariety $$\omega^{-1}(E_h)=\overline{\mathcal{D}}\cap
\{H_1(\Omega)=0\klk H_{N-h}(\Omega)=0\}$$ of $\A^L$. From the
B\'ezout Inequality (in the variant of \cite[Proposition
2.3]{HeSc82}) we conclude now $$\deg E_h\le\deg
\omega^{-1}(E_h)\le\deg\overline{\mathcal{D}}
\cdot\Delta_2^{dim\,\overline{\mathcal{D}}}\le\deg\overline{\mathcal{D}}
\cdot\Delta_2^L.$$ Thus, if $\overline{\mathcal{O}}$ is
equidimensional of dimension $h$, we have
$\overline{\mathcal{O}}=E_h$ and therefore
$$\deg\overline{\mathcal{O}}\le\deg\overline{\mathcal{D}}
\cdot\Delta_2^L\le \big((1+K\Delta_1)\Delta_2\big)^L.$$ In the
general case we obtain the following estimate:
$$\deg\overline{\mathcal{O}}=\sum_{h=0}^L\deg E_h\le
(L+1)\Delta_2^L\deg\overline{\mathcal{D}} \le
(L+1)\big((1+K\Delta_1)\Delta_2\big)^L.$$ This proves the first
statement of the Lemma.

In the second part of the proof we consider the closed subvariety
\[
\begin{array}{l}
V=\{(F,y^{(1)}\klk y^{(m)});F\in\overline{\mathcal{O}},
y^{(1)}\klk y^{(m)}\in\A^t,\\
\qquad\qquad\qquad\qquad\qquad\qquad\qquad\qquad\qquad\quad
F(y^{(1)})=\cdots=F(y^{(m)})=0\}\end{array}\] of the affine space
$\A^N\times\A^{mt}$ and the morphisms of algebraic varieties
\linebreak
$\pi_1:V\to \A^N$ and $\pi_2:V\to\A^{mt}$ induced by the canonical
projections of $\A^N\times\A^{mt}$ onto $\A^N$ and $A^{mt}$.

Since any polynomial of $\overline{\mathcal{O}}$ has degree at
most $\Delta$, we deduce from the B\'ezout Inequality the estimate
\begin{equation}
\label{uno} \deg V\le
\deg\overline{\mathcal{O}}\cdot\Delta^m.\end{equation}

Let $\mathcal{C}_1\klk\mathcal{C}_s$ be the irreducible components
of $V$ whose $\pi_1$--image contains at least one nonzero
polynomial of $\overline{\mathcal{O}}$. Let
$V^*:=\bigcup_{1\le j\le s} \mathcal{C}_j$. Thus
$\pi_2(V^*)=\bigcup_{1\le j\le s}\pi_2
(\mathcal{C}_j)$ is the set of all ``incorrect" test sequences of
length $m$ for the object class $\overline{\mathcal{O}}$. From
(\ref{uno}) we deduce the estimate
\begin{equation}
\label{dos} \deg V^*\le
\deg\overline{\mathcal{O}}\cdot\Delta^m.\end{equation}

Let $1\le j\le s$. There exists a polynomial
$F\in\overline{\mathcal{O}}$ with $F\not=0$ and
\linebreak
$F\in\pi_1(\mathcal{C}_j)$. Observe that the fiber $\pi_1^{-1}(F)$
is isomorphic to the equidimensional algebraic variety
$$\{(y^{(1)}\klk y^{(m)})\in\A^{mt}; y^{(1)}\klk
y^{(m)}\in\A^t,F(y^{(1)})=\cdots=F(y^{(m)})=0\}.$$ Thus $F\not=0$
implies $dim\,\pi_1^{-1}(F)=m(t-1)$. Applying the Theorem of
Fibers (see e.g. \cite[I.8, Corollary]{Mumford88}) to the morphism
of irreducible affine varieties
$$\pi_1|_{\mathcal{C}_j}:\mathcal{C}_j\to\overline{\pi_1(\mathcal{C}_j)}$$
we deduce
$$dim\,\mathcal{C}_j-dim\,\overline{\pi_1(\mathcal{C}_j)}\le
m(t-1).$$ Since $\overline{\pi_1(\mathcal{C}_j)}$ is contained in
the affine variety $\overline{\mathcal{O}}$ we conclude
$dim\,\overline{\pi_1(\mathcal{C}_j)} \le
dim\,\overline{\mathcal{O}}$ and therefore $dim\,\mathcal{C}_j-
dim\,\overline{\mathcal{O}}\le m(t-1)$. This implies
\begin{equation}
\label{tres}
dim\,\mathcal{C}_j\le m(t-1)+dim\,\overline{\mathcal{O}}.\end{equation}

By assumption the data structure $\mathcal{D}$ encodes the object
class $\mathcal{O}$ holomorphically by means of the encoding
$\omega$. This means that the encoding determines a morphism of
affine varieties $\overline{\mathcal{D}}\to\overline{\mathcal{O}}$
which contains $\mathcal{O}$ in its image. Therefore this morphism
is dominant and this implies $dim\,\overline{\mathcal{O}}\le
dim\,\overline{\mathcal{D}}.$ From (\ref{tres}) we conclude now
$$dim\,\mathcal{C}_j\le m(t-1)+dim\,\overline{\mathcal{D}}.$$
Since $1\le j\le s$ was arbitrary, we obtain the estimate
\begin{equation}
\label{cuatro}
dim\, V^*\le m(t-1)+dim\,\overline{\mathcal{D}}.
\end{equation}
This implies $dim\,\overline{\pi_2(V^*)}\le m(t-1)+dim\,
\overline{\mathcal{D}}\le m(t-1)+L.$
\smallskip

Before continuing with the proof, observe that by assumption $m\ge
2L+2>L$ and therefore $mt>m(t-1)+L$ holds. Hence
$\overline{\pi_2(V^*)}$ is a proper closed subset of $\A^{mt}$.
Thus  any element $\gamma:=(\gamma_1\klk\gamma_m)$ of the Zariski
open, dense subset
$\mathcal{U}:=\A^{mt}\setminus\overline{\pi_2(V^*)}$ of $\A^{mt}$
with $\gamma_1\klk\gamma_m\in\A^t$ is a correct test sequence for
the object class $\overline{\mathcal{O}}$.
\smallskip

Let us finally pass to the third and final part of the proof. For
\linebreak $1\le k\le m$ and $1\le\ell\le t$ let $Y_{ke}$ be a new
indeterminate and let $ H_{ke}:=\prod_{\mu\in M}
(Y_{k\ell}-\mu)$. Thus $H_{k\ell}$ is a univariate polynomial of
degree $\# M$ belonging to the polynomial ring $k[Y_{k\ell}]$. We
consider the indeterminates $Y_{k\ell}$ with $1\le k\le m$, $1\le
\ell\le t$ as coordinate functions of the affine space $\A^{mt}$.
Observe that $M^{mt}=\{H_{k\ell}=0;1\le k\le m, 1\le\ell\le t\}$
holds and that the set of ``incorrect" test sequences contained in
$M^{mt}$, namely $\pi_2(V^*)\cap M^{mt}=\pi_2\left(V^*\cap
\{H_{k\ell}=0;1\le k\le m, 1\le\ell\le t\}\right)$, is a finite
$k$--definable (and hence Zariski closed) subset of $\A^{mt}$.

From \cite[Proposition 2.3]{HeSc82} (i.e. from the B\'ezout
Inequality) and from (\ref{dos}), (\ref{cuatro}) we conclude now
\[
\begin{array}{rcl}
\#\left(\pi_2(V^*)\cap M^{mt}\right)&=&\# \pi_2\left(V^*\cap
\{H_{k\ell}=0;1\le k\le m, 1\le\ell\le t\}\right)\\ \\ &\le
&\deg(V^*\cap \{H_{k\ell}=0;1\le k\le m, 1\le\ell\le t\})\\ \\
&\le &\deg(V^*)\, (\# M)^{dim\, V^*}\\ \\ &\le &
\deg(\overline{\mathcal{O}})\, \Delta^m\,(\# M)^{m(t-1)+dim\,
\overline{\mathcal{D}}}\\ \\ &\le &
\deg(\overline{\mathcal{O}})\,\Delta^m \,(\# M)^{m(t-1)+L}.
\end{array}
\]

Suppose now that the points of the finite set $M^t$ are
equidistributed. By assumption we have $m\ge 2L+2$, $\# M\ge
\Delta^2(\deg \overline{\mathcal{O}})^{\frac{1}{L}}$ and
$\Delta\ge 1$. From the estimate $\#\left(\pi_2(V^*)\cap
M^{mt}\right)\le
\deg(\overline{\mathcal{O}})\,\Delta^m\,(\#M)^{m(t-1)+L}$ we
deduce that the probability of finding in $M^{mt}$ by a random
choice an ``incorrect" test sequence for the object class $\omo$
is at most
\[\begin{array}{rcl}
\displaystyle
\frac{\deg(\omo)\,\Delta^m}{(\# M)^{m-L}}&\le&
\displaystyle
\frac{\deg(\omo)\,\Delta^m}{\# M\, \left(
\Delta^2(\deg \mo)^{\frac{1}{L}}\right)^{m-L-1}} \\ \\
&=&\displaystyle
\frac{\deg(\omo)\,\Delta^{2(L+1)}}{\# M\,\Delta^m\,
(\deg\mo)^{\frac{1}{L}(m-L-1)}} \\ \\
&\le&\displaystyle
\frac{\deg(\omo)}{\# M\,\big(\deg(\omo)\big)^{1+\frac{1}{L}}}
\le \frac{1}{\# M}\le \frac{1}{2}
\end{array}\]
(recall that by assumption $M$ has at least two elements). Hence
the probability of finding in $M^{mt}$ by a random choice a
correct test sequence for the object class $\omo$ is at least
$$1-\displaystyle\frac{1}{\# M}\ge\frac{1}{2}.$$

Since this probability is positive, we conclude that $M^{mt}$
really contains a correct test sequence $\gamma=(\om{\gamma}{m})$
with $\om{\gamma}{m}\in\A^t$ for the object class $\omo$.
\end{prf}\cqfd

\begin{corollary}
\label{corollary:tcs2} Let notations and assumptions be as in
Lemma \ref{ts}. Let $M$ be a finite subset of $k$ of cardinality
at least $max\{\Delta^2(\deg\omo)^{\frac{1}{L}},2\}$ and let $m\ge
4L+2$. Then there exist points $\om\gamma m$ of $M^t$ such that
$\gamma:=(\om\gamma m)$ is an identification sequence for the
object class $\omo$ (and hence for $\mo$). Suppose that the points
of the finite set $M^t$ are equidistributed. Then the probability
of finding in $M^{mt}$ by a random choice such an identification
sequence is at least $1-\frac{1}{\# M}\ge\frac{1}{2}$.
\end{corollary}

\begin{prf}
We use the same notations and assumptions as in the proof of Lemma
\ref{ts}. Let $\omega:\md\to\mo$ be the given $k$--definable
holomorphic encoding of the object class $\mo$. Let
$\md_*:=\md\times\md$, $\mo_*:=\{F_1-F_2;F_1,F_2\in\mo\}$ and let
$\omega_*:\md_*\to\mo_*$ be the encoding of the object class
$\mo_*$ defined by $\omega_*(D_1,D_2)=\omega(D_1)-\omega(D_2)$ for
$(D_1,D_2)\in\md_*$.

One verifies immediately that the data structure $\md_*$ and the
object class $\mo_*$ are $k$--constructible subsets of $\A^{2L}$
and $\A^N$ respectively and that $\omega_*$ is a $k$--definable
holomorphic encoding of the object class $\mo_*$. In particular
$\mo_*$ turns out to be a $k$--constructible object class in the
sense introduced before. Furthermore $\Delta$ is an upper bound
for the degree of the $t$--variate polynomials over $\overline{k}$
contained in the object class $\mo_*$. From \cite{Heintz83},
Proposition 2 and Lemma 2 we deduce the estimate
$\deg\omo_*\le(\deg\omo)^2$. Hence $\#
M\ge\Delta^2(\deg\omo)^{\frac{1}{L}} $ implies $\# M\ge
\Delta^2(\deg \omo_*)^{\frac{1}{2L}}$.

Suppose now that the points of the finite set $M^t$ are
equidistributed. From Lemma \ref{ts} we deduce that the
probability of finding in $M^{mt}$ by a random choice a \cts\ for
the object class $\omo_*$ is at least $1-\frac{1}{\# M}\ge
\frac{1}{2}$.

Let $\gamma=(\om\gamma m)\in M^{mt}$ with $\om\gamma m\in M^t$
such a correct test sequence and let $F_1,F_2$ be given elements
of $\omo$ (thus $F_1$ and $F_2$ are $t$--variate polynomials over
$\overline{k}$). Suppose that $F_1(\gamma_1)=F_2(\gamma_1),\dots,
F_1(\gamma_m)=F_2(\gamma_m)$ holds. Hence, for $F:=F_1-F_2$, we
have $F(\gamma_1)=\cdots=F(\gamma_m)=0$. Since $F$ belongs to the
object class $\omo_*$ and $\gamma$ is a \cts\ for $\omo_*$ we
infer $F=0$. This implies $F_1=F_2$.

In conclusion, we see that $\gamma$ is an identification sequence for
the object class $\omo$. Since the probability of finding such identification
sequences in $M^{mt}$ is positive, we infer that $M^{mt}$ contains at
least one of them.
\end{prf}\cqfd

Let $\mo$ be $k$--definable object class of polynomial functions
and $\omega:\md\to\mo$ be a $k$--definable holomorphic encoding of
$\mo$ by a $k$--constructible data structure of size $L$. By means
of the data structure $\md$ we are able to answer the value
question about the object class $\mo$ holomorphically. In this
sense, an identification sequence $\gamma$ of length $m$ allows to
answer the {\em identity question} about the object class $\mo$
{\em holomorphically}. From Corollary \ref{corollary:tcs2} we
conclude that there exist always {\em short} identification
sequences (of length $m$ linear in $L$) and that they are easy to
find by means of a suitable random choice. This means that the
identity question about the object class $\mo$ can always be
answered ``efficiently".

%-----------------------------------------------------------------
%-----------------------------------------------------------------
%-----------------------------------------------------------------
%-----------------------------------------------------------------

\subsubsection{Correct test and identification sequences for
circuit encodings.} \label{circuitcorrect} In order to exemplify
the ideas behind Lemma \ref{ts} and Corollary
\ref{corollary:tcs2} of Section \ref{holocorrect} we
are now going to apply the concept of identification sequence to
circuit encoded object classes of polynomial functions.\medskip

Let $\varepsilon$ be a new indeterminate and let us consider
$\varepsilon$ as a parameter and $\om Y t$ as variables. Let
$F\in\overline{k}[\om Y t]$. We denote by $L(F)$ the minimal
nonscalar size over $\overline{k}$ of all totally division--free
arithmetic circuits with inputs $\om Y t$ and scalars in
$\overline{k}$ which evaluate the polynomial $F$. Moreover we
denote by $\overline{L}(F)$ the minimal nonscalar size over
$\overline{k}(\varepsilon)$ of all essentially division--free
arithmetic circuits which evaluate a rational function of the form
$F+\varepsilon Q$ with $Q$ belonging to
$\overline{k}[\varepsilon,\om Y t]_\varepsilon$. Obviously we have
$\overline{L}(F)\le L(F)$. We call $L(F)$ the {\em nonscalar
(sequential time) complexity of $F$ over $\overline{k}$} and
$\overline{L}(F)$ the corresponding {\em approximative}
complexity. Let $L\in\N$ and let $W_{L,t}:= \{F\in\overline{k}[\om
Y t];L(F)\le L\}$. From \cite[Chapter 9, Exercise 9.18]{BuClSh97}
(see also \cite[Theorem 3.2]{HeSc82}) we deduce that all
polynomials contained in $W_{L,t}$ have degree bounded by $2^L$
and that $W_{L,t}$ forms a $k$--constructible object class which
has a $k$--definable holomorphic encoding by the data structure
$\A^{(L+t+1)^2}$. Moreover any polynomial $F\in\okyot$ with
$\overline{L}(F)\le L$ has degree at most $2^L$.

Let $N\in\N$ with $N\ge 2^L$. Then $W_{L,t}$ can be considered as
a $k$--constructible subset of $\A^N$. From \cite{Alder84}, Lemma
2 and Satz 4 one deduces easily the following statement:
$$\overline{W}_{L,t}:=\{F\in\okyot; \overline{L}(F)\le L\}.$$ In
this sense the Zariski closure of the object class $W_{L,t}$ has a
natural interpretation as the set of polynomials of $\okyot$ which
have approximative nonscalar (sequential time) complexity over
$\ok$ at most $L$.

Finally observe that $W_{L,t}$ and $\overline{W}_{L,t}$ are cones
and contain the zero polynomial. In particular any identification
sequence of $W_{L,t}$ or $\overline{W}_{L,t}$ is a \cts.

\begin{corollary} {\rm (compare \cite[Theorem 4.4]{HeSc82} and
\cite[Lemma 3]{GiHe01})} \label{corollary:tcs4} Let notations be
as before and let $L$, $m$, $t$ be natural numbers with $m\ge
4(L+t+1)^2+2$. Let $M$ be a finite subset of $k$ of cardinality at
least $2^{4(L+1)}$. Then there exist points $\gammam$ of $M^t$
such that $\gamma:= (\gammam)$ is an identification sequence for
the object class $\overline{W}_{L,t}$ of all polynomials
$F\in\okyot$ which have approximative nonscalar (sequential time)
complexity over \ok\ at most $L$.

Suppose that the points of the finite set $M^t$ are
equidistributed. Then the probability of finding in $M^{mt}$ by a
random choice such an identification sequence is at least
$1-\frac{1}{\# M}\ge\frac{1}{2}$.
\end{corollary}

\begin{prf}
Let $N\ge 2^L$, $r:=(L+t+1)^2$ and let $\om Z r$ be new
indeterminates. From \cite[Chapter 9, Exercise 9.18]{BuClSh97}
(compare also \cite[Theorem 2.1]{Schnorr78}) we deduce that there
exist $N$ polynomials of $k[\om Z r]$ having degree at most
$L\,2^{L+1}+2$ which induce a $k$--definable holomorphic encoding
$\omega:\A^r\to W_{L,t}$ of the object class $W_{L,t}$ which we
consider as a $k$--constructible subset of $\A^N$.

Taking into account that
$\overline{W}_{L,t}=\overline{\omega(\A^r)}$ is irreducible, we
deduce from Lemma \ref{ts} the estimate $\deg
\overline{W}_{L,t}\le(L\, 2^{L+1}+2)^r$. This implies
\begin{equation}\label{nuevouno}
(\deg \overline{W}_{L,t})^{\frac{1}{r}}\le L\,2^{L+1}+2.
\end{equation}
Observe that by hypothesis $m\ge 4(L+t+1)^2+2=4r+2$ holds and that
any polynomial contained in $W_{L,t}$ has degree at most
$\Delta:=2^L$. From the assumption $\# M\ge 2^{4(L+1)}$ and
(\ref{nuevouno}) we deduce $\# M\ge
\Delta^2(\deg\overline{W}_{L,t})^{ \frac{1}{r}}$. The statement to
prove follows now immediately from Corollary \ref{corollary:tcs2}.
\end{prf}\cqfd

%-----------------------------------------------------------------
%-----------------------------------------------------------------
%-----------------------------------------------------------------
%-----------------------------------------------------------------

\subsection{Encodings of polynomial functions by values.}
\label{valencode}

In this subsection we
are going to prove the first main
result of
this paper.\medskip

Let $\mathcal{O}$ be a $k$--constructible object class of
polynomial functions, $\mathcal{D}$ a $k$--constructible data
structure and $\omega :\mathcal{D}\to \mathcal{O}$ a
$k$--definable holomorphic encoding. Our first main result
(Theorem \ref{theorem:encode1} below) may be stated succinctly as
follows:\smallskip

\noindent assume that $\mathcal{O}$ is a class of polynomial
functions and that its encoding by $\mathcal{D}$ is holomorphic.
Suppose furthermore that the ambient space of $\mathcal{O}$ is affine
and contains $\mathcal{O}$ as a cone (i.e. we assume that $\mathcal{O}$
is closed under multiplication by scalars). Then there exists
a $k$--definable data structure $\overline{\mathcal{D}}$ which encodes
the closure class $\overline{\mathcal{O}}$ of $\mathcal{O}$ {\em
continuously} (with respect to the Zariski topologies of
$\overline{\mathcal{D}}$ and $\overline{\mathcal{O}}$) and
{\em unambiguously}. In particular,
$\overline{\mathcal{O}}$ and $\overline{\mathcal{D}}$ are
homeomorphic topological spaces. Moreover the size of
$\overline{\mathcal{D}}$ (i.e. the dimension of its ambient
space) is {\em linear} in the size of $\mathcal{D}$.\smallskip

In other words, we may always replace {\em efficiently} the given
data structure $\mathcal{D}$ by an unambiguous one, say
$\overline{\mathcal{D}}$, if we are only interested in a {\em
topological} characterization of the object class
$\overline{\mathcal{O}}$ (or $\mathcal{O}$). By means of
$\overline{\mathcal{D}}$ we are able to answer efficiently the
{\em identity} question about $\overline{\mathcal{O}}$, but {\em
not necessarily} the {\em value} question. The assumption that
the object class forms a cone in case that $\mathcal{O}$ has
affine ambient space is not restrictive in the context of this
paper, since $\mathcal{O}$ will be typically a class of functions
closed under multiplication by scalars. On the other hand, this
assumption guarantees that the encoding of the object class
$\overline{\mathcal{O}}$ by the data structure
$\overline{\mathcal{D}}$ is not only continuous, but also a
closed map with respect to the Zariski topologies of
$\overline{\mathcal{D}}$ and $\overline{\mathcal{O}}$.

In the Appendix of this paper (Section \ref{unisequence}) we shall
formulate a slight generalization of Theorem
\ref{theorem:encode1} below. \medskip

First we synthesize the essence of the technical Lemma \ref{ts}
and its Corollary \ref{corollary:tcs2} of
Section \ref{holocorrect}
in terms of continuous encodings.\medskip

Let $\mo\subset\okyot$ be a $k$--constructible object class of
polynomial functions and let $\gamma=(\gammam)\in k^{mt}$ with
$\gammam\in k^t$ and $m\ge 1$ be an identification sequence for
$\omo$ (from Corollary \ref{corollary:tcs2} one deduces easily
that for $m\in\N$ sufficiently large such an identification
sequence always exists).

Suppose now that $\mo$ is a {\em cone} in $\okyot$. Then
$\omo$ is a cone too. Let $\sigma:\omo\to\A^m$ be the map defined
by $\sigma(F):= \big(F(\gamma_1)\klk F(\gamma_m)\big)$ for
$F\in\omo$. Observe that $\sigma$ is the restriction of a
$k$--definable linear map $\A^N\to\A^m$, where $\A^N$ with $N\ge
1$ is a suitable affine ambient space which contains $\mo$ and
$\omo$ as cones. Thus $\sigma$ is homogeneous of degree one and
represents an {\em injective}, $k$--definable morphism of affine
varieties. Therefore $\sigma(\omo)$ is a $k$--definable subset of
$\A^m$. Since $\sigma$ is homogeneous of degree one and $\omo$ is
a cone, the image $\sigma(\omo)$ is a cone too. Hence the Zariski
closure $\md^*$ of $\sigma(\omo)$ in $\A^m$ is a $k$--definable
cone of $\A^m$ and $\sigma$ induces a dominant morphism of affine
varieties which maps $\omo$ into $\md^*$ and is again
homogeneous of degree one. We denote this morphism by $\sigma:
\omo\to\md^*$.
\begin{lemma}
\label{lemma4:encode} Let notations and assumptions be as before.
Then $\sigma: \omo\to\md^*$ is a finite, bijective, $k$--definable
morphism of affine varieties. Let $\mathcal{C}$ be an arbitrary
$k$--definable irreducible component of $\omo$. Then
$\sigma|_{\mc}$ is a birational, $k$--definable (finite and
bijective) morphism of $\mc$ onto the Zariski closed set
$\sigma(\mc)$.
\end{lemma}
\begin{prf}
Let $\om ZN$ be the coordinate functions of $\A^N$. There exist
linear polynomials $\om S m\in k[\om Z N]$ such that $\sigma$ is
the restriction of the linear map $(\om Sm)$ to the closed
subvariety $\omo$ of $\A^N$. Since $\omo$ and $\md^*$ are
$k$--definable Zariski closed cones of the affine spaces $\A^N$
and $\A^m$ respectively, they are definable by homogeneous
polynomials over $k$. Moreover $\omo$ and $\md^*$ contain the
origins of the affine spaces $\A^N$ and $\A^m$ respectively. From
the injectivity of $\sigma:\omo\to\md^*$ we deduce therefore that
$\omo\cap\{S_1=0\klk S_m=0\}$ contains only the origin of $\A^N$.
This implies that the homogeneous map $\sigma$ induces a finite
morphism between the closed projective subvarieties of $\Pe^{N-1}$
and $\Pe^{m-1}$ associated to the cones $\omo$ and $\md^*$
respectively. In fact, the standard proof of this classical result
implies something more, namely that also the morphism $\sigma:
\omo\to\md^*$ is finite (see \cite{Shafarevich84}, I.5.3, Theorem
8 and proof of Theorem 7). In particular, $\sigma$ is a surjective
closed map. Since $\sigma$ is also injective we conclude that
$\sigma$ is bijective.

Let $\mc$ be an arbitrary $k$--definable irreducible component
of $\omo$. Since $\sigma$ is a closed map we conclude that
$\sigma(\mc)$ is a closed irreducible subvariety of $\md^*$.
Since $\sigma$ is injective we infer that $\sigma|_\mc:\mc\to
\sigma(\mc)$ is a bijective, $k$--definable morphism of affine
varieties. Since for any point $y\in\sigma(\mc)$ we have
$\#\big(\sigma^{-1}(y)\cap\mc\big)=1$ we deduce from
\cite[Proposition 3.17]{Mumford88} that $k\big(\sigma(\mc)\big)=
k(\mc)$ holds. Hence $\sigma|_\mc$ is a birational morphism.
\end{prf}\cqfd

From Lemma \ref{lemma4:encode} we deduce that with respect to the
Zariski topologies of $\omo$ and $\md^*$, the morphism
$\sigma:\omo\to\md^*$ is a homeomorphism and that
$\md^*=\sigma(\omo)$ holds. Consider now $\md^*\subset \A^m$ as a
data structure. Then $\omega^*:=\sigma^{-1}:\md^*\to\omo$ is an
{\em unambiguous} encoding of the object class $\omo$ which is
{\em continuous} with respect to the Zariski topologies of $\md^*$
and $\omo$. Suppose that $\omega^*$ allows to answer the value
question about the object class $\omo$ holomorphically. Then from
Remark \ref{holo} we deduce that $\omega^*:\md^*\to\omo$ is a
$k$--definable morphism of algebraic varieties and therefore
$\omega^*$ is an unambiguous, $k$--definable and
{\em(bi--)holomorphic} encoding of the object class $\omo$ by the
data structure $\md^*$. We shall see later that in general this
will not be the case (see Corollary \ref{remark2prima,3.1} and
Theorem \ref{theorem:section3.1}). Suppose for the moment $k:=\Q$
and $\ok:=\C$. Since $\sigma$, the inverse map of the unambiguous
encoding $\omega^*$, is a morphism of algebraic varieties, we
conclude that $\sigma $ is continuous with respect to the strong
topologies of $\omo$ and $\md^*$. If $\omega^*$ is continuous with
respect to the strong topology, this implies that $\omega^*$ is a
{\em robust} encoding in the sense of Definition \ref{robustcont}.

However, $\omega^*$ may be not continuous with respect to the
strong topologies of $\omo$ and $\md^*$. On the other hand,
$\omega^*$ induces a map $\Omega$ between the projective
subvarieties of $\Pe^{m-1}(\C)$ and $\Pe^{N-1}(\C)$ associated to
the cones $\md^*$ and $\omo$. The map $\Omega$ encodes the
projective variety associated to the cone $\omo$ by the projective
variety associated to the cone $\md^*$ and is continuous with
respect to the corresponding strong topologies.

We may summarize the main results of this section by the
following statement:
\begin{theorem}
\label{theorem:encode1} Let $\mo$ be a $k$--constructible object
class of polynomial functions belonging to $\okyot$. Let $\Delta$
be an upper bound for the degree of the polynomials contained in
$\mo$. Suppose that $\mo$ is a cone in $\okyot$. Assume that there
is given a $k$--constructible data structure $\md\subset\A^L$ and
a $k$--definable holomorphic encoding $\md\to\mo$. Let $m\ge 4L+2$
and let $M$ be a finite subset of $k$ of cardinality at least
$max\{\Delta^2(\deg\omo)^{ \frac{1}{L}},2\}$. Then there exist a
$k$--definable, Zariski closed cone $\md^*$ of $\A^m$ and a
continuous encoding $\omega^*:\md^*\to\omo$ of the object class
$\omo$ by the data structure $\md^*$ which satisfies the following
conditions:
\begin{itemize}
\item[$(i)$] $\omega^*$ is a homeomorphism between the data
structure $\md^*$ and the object class $\omo$,
\item[$(ii)$] there  exist a point $\gamma:=(\gamma_1\klk\gamma_m)\in M^{mt}$
with $\gamma_1\klk\gamma_m\in M^t$ such that for any $F\in\omo$
the identity $(\omega^*)^{-1}(F)=\big(F(\gamma_1)\klk
F(\gamma_m)\big)$ holds,
\item[$(iii)$] $(\omega^*)^{-1}:\omo\to\md^*$ is a $k$--definable, bijective
and finite morphism of affine varieties. The morphism
$(\omega^*)^{-1}$ is homogeneous of degree one,
\item[$(iv)$] for any $k$--definable irreducible component $\mc$ of
$\omo$ the restriction map
$(\omega^*)^{-1}|_\mc:\mc\to(\omega^*)^{-1}(\mc)$ is a birational
$k$--definable (finite and surjective) morphism of $\mc$ onto the
irreducible Zariski closed set $(\omega^*)^{-1}(\mc)$.
\end{itemize}
In particular $\omega^*$ is an unambiguous continuous encoding of
the object class $\omo$ by the data structure $\md^*$. The
encoding $\omega^*$ is holomorphic if and only if $\omega^*$
allows to answer holomorphically the value question about the
object class $\omo$.

In case  $k:=\Q$, $\ok:=\C$ and $\omega^*$
continuous with respect to the strong topology, the encoding
$\omega^*$ is robust (in the sense of Definition
\ref{robustcont}).

Suppose that the elements of the finite set $M^t$ are
equidistributed. Then the probability of finding by a random
choice a point $\gamma:=(\gammam)\in M^{mt}$ with $\gammam\in M^t$
such that the map $\sigma_\gamma:\omo\to\A^m$ defined by
$\sigma_\gamma(F):=\big(F(\gamma_1)\klk F(\gamma_m)\big)$ for
$F\in\mo$ induces a $k$--definable, bijective morphism of $\omo$
onto a Zariski closed cone $\md_\gamma^*$ of $\A^m$ is at least
$1-\frac{1}{\# M}\ge \frac{1}{2}$. Any such morphism
$\sigma_\gamma:\omo\to\md_\gamma^*$ defines by
$\omega^*_\gamma:=\sigma_\gamma^{-1}$ a continuous unambiguous
encoding of the object class $\omo$ by the data structure
$\md_\gamma^*$. This encoding satisfies conditions $(i)$--$(iv)$.
\end{theorem}

The proof of Theorem \ref{theorem:encode1} is  an immediate
consequence of Corollary \ref{corollary:tcs2}, Lemma
\ref{lemma4:encode} and the subsequent considerations.

The following statement represents a version of Theorem
\ref{theorem:encode1} for object classes of
arithmetic--circuit--represented polynomials.

\begin{corollary} {\rm \cite[Lemma 4]{GiHe01}}
\label{encode3} Let notions and notations be as in Corollary
\ref{corollary:tcs4}. Let $L$, $m$, $t$ be natural numbers with
$m\ge 4(L+t+1)^2+2$. Let $M$ be a finite subset of cardinality at
least $2^{4(L+1)}$. Let $\overline{W}_{L,t}$ be the object class
of all polynomials $F\in\okyot$ which have approximative nonscalar
(sequential) complexity over $\ok$ at most $L$. Then
$\overline{W}_{L,t}$ is a cone and there exists a $k$--definable,
Zariski closed cone $\md^*_{L,t}$ of $\A^m$ and a continuous
encoding $\omega^*:\md^*_{L,t}\to\overline{W}_{L,t}$ of the object
class $\overline{W}_{L,t}$ by the data structure $\md^*_{L,t}$
which satisfies the conditions (i)--(iv) of Theorem
\ref{theorem:encode1}. The encoding $\omega^*$ is holomorphic if
and only if $\omega^*$ allows to answer holomorphically the value
question about the object class $\overline{W}_{L,t}$.

In case $k:=\Q$, $\ok:=\C$ and $\omega^*$ continuous
with respect to the strong topology, the encoding $\omega^*$ is
robust (in the sense of Definition \ref{robustcont}).

Suppose
that the elements of the finite set $M^t$ are equidistributed.
Then we may find by a random choice with probability of success
at least $1-\frac{1}{\# M}$ a point $\gamma=(\gammam)\in M^{mt}$
with $\gammam\in M^t$ such that the map
$\sigma_\gamma:\overline{W}_{L,t}\to\A^m$ defined by
$\sigma_\gamma(F):=\big(F(\gamma_1)\klk F(\gamma_m)\big)$ for
$F\in \overline{W}_{L,t}$ produces as in Theorem
\ref{theorem:encode1} a $k$--definable, Zariski closed cone
$\md^*_{L,t,\gamma}$ of $\A^m$ and a continuous encoding
$\omega^*_\gamma:\md^*_{L,t,\gamma}\to\overline{W}_{L,t}$.
\end{corollary}
\begin{prf}
Since in the nonscalar complexity model $\ok$--linear operations
are free, we conclude that ${W}_{L,t}=\{F\in\okyot;L(F)\le L\}$ is
a cone of $\okyot$. Therefore its closure $\overline{W}_{L,t}$ is
a cone too. The statement of Corollary \ref{encode3} follows now
immediately from Corollary \ref{corollary:tcs4} and Lemma
\ref{lemma4:encode}.
\end{prf}\cqfd

We call a continuous encoding of an object class of polynomial
functions as in Theorem \ref{theorem:encode1} and Corollary
\ref{encode3} of this section and Corollary \ref{coro:encode2} of
Section \ref{unisequence} an {\em encoding by an identification
sequence} or simply an {\em encoding by values}. An encoding by
an identification sequence allows us to answer the identity
question about the object class $\mo$. However, the corresponding
value question requires a holomorphic encoding. In the next
section we shall exhibit an example of a $\Q$--constructible
object class $\mo$ of univariate polynomials which has a
$\Q$--definable, holomorphic, {\em robust} but {\em ambiguous}
encoding by a data structure of small size. However we shall show
that any {\em holomorphic} encoding of $\mo$ by an identification
sequence requires a data structure of (exponentially) big size.

A given object class of polynomial functions has many, mostly
artificial encodings. However, encodings by values seem
particularly natural. This becomes evident in the situation of
Corollary \ref{encode3}. Encodings of object classes of polynomial
functions by arithmetic circuits are typically ambiguous. In
Corollary \ref{encode3} a given encoding of an object class of
polynomial functions by arithmetic circuits is replaced by an {\em
unambiguous} continuous and robust encoding by means of an
identification sequence (observe that evaluation is particularly
well--adapted to circuit encoding).

%-----------------------------------------------------------------
%-----------------------------------------------------------------
%-----------------------------------------------------------------
%-----------------------------------------------------------------

\subsection{Unirational encodings.}
\label{uniencode}
Let $\mo$ be a $k$--constructible object class
and let $\md$ be a $k$--constructible data structure of size $L$,
contained in the ambient space $\A^L$ or $\Pe^L$. Let
$\omega:\md\to\mo$ be a $k$--definable holomorphic encoding of
the object class $\mo$ by the data structure $\md$. We call
$\omega$ {\em unirational} if $\md$ contains a nonempty, Zariski
open set of its ambient space. Suppose that $\omega$ is
unirational. Then $\omd$ equals its ambient space $\A^L$ or
$\Pe^L$ and $\omo$ is an irreducible $k$--Zariski closed set in
some suitable affine or projective space. We call the encoding
$\omega$ {\em rational} if it defines a birational map between
the ambient space and $\omo$.

Let $L$ and $t$ be natural numbers. Then the generic computation
scheme of length $L$ in the nonscalar sequential complexity model
(see \cite{BuClSh97}, Chapter 9, Theorem 9.9 and Exercise 9.18,
or \cite{Heintz89}) defines a unirational encoding of the object
classes $W_{L,t}$ and $\overline{W}_{L,t}$ of polynomial
functions of $\okyot$ having exact or approximative (sequential)
nonscalar complexity over $\ok$ at most $L$. Analogously, the
standard representation of polynomials of $\okyot$ of degree at
most $d$ by their coefficients is a rational encoding of size
$\begin{pmatrix}d+t\\t\end{pmatrix}$. Similarly the $L$--sparse
polynomials of $\okyot$ containing only a previously fixed set of
$L$ monomials are rationally encoded by the data structure $\A^L$.

One may ask why we do not limit our attention exclusively to
unirational encodings of object classes. A technical reason for
this is that  such encodings represent only a limited range of
object classes. In order to exemplify this, let us observe that a
limitation to unirational data structures would automatically
exclude from our considerations important object classes as e.g.
the set of all $k$--definable equidimensional projective varieties
of dimension $r$ and degree $d$ contained in a projective space
$\Pe^n$ with $n\ge r$. The traditional data structures for these
object classes are the Chow varieties which encode (unambiguously)
a given object by its Chow coordinates.

Similarly the Hilbert varieties are data structures which encode
unambiguously the (reduced) projective subvarieties of a given
projective space with previously fixed Hilbert polynomial. The
natural topology of Chow and Hilbert varieties induces a topology
on the object classes they represent and hence a notion of limit
object. Typical Chow varieties, encoding complete intersection
varieties, are unirational and it is not clear whether they could
be also rational. In general, Hilbert varieties cannot be expected
to be unirational.

%-----------------------------------------------------------------
%-----------------------------------------------------------------
%-----------------------------------------------------------------
%-----------------------------------------------------------------
%-----------------------------------------------------------------
%-----------------------------------------------------------------
%-----------------------------------------------------------------
%-----------------------------------------------------------------

\section{Two paradigmatic object classes.}
\label{section4real} In this section we are going to exhibit two
paradigmatic object classes of polynomial functions and to discuss
different holomorphic encodings of them. We shall always assume
$k:=\Q$ and $\overline{k}:=\C$.

%-----------------------------------------------------------------
%-----------------------------------------------------------------
%-----------------------------------------------------------------
%-----------------------------------------------------------------

\subsection{First paradigm.}
\label{section3.1} Let $d$ be a natural number, let $U$ and $Y$ be
indeterminates over $\Q$ and let
$F_d:=\sum_{j=0}^d(U^d-1)U^jY^j\in\Q[U,Y]$. We are going to
interpret $U$ as parameter and $Y$ as variable. Let us consider
the object class of univariate polynomials
$\mo_d:=\{F_d(u,Y);u\in\A^1\}$ and the encoding
$\omega_d:\A^1\to\mo_d$ defined for $u\in\C$ by
$\omega_d(u):=F_d(u,Y)$. Representing the polynomials belonging to
$\mo_d$ by their coefficients, we identify the object class
$\mo_d$ with the corresponding subset of $\A^{d+1}$. With this
interpretation $\omega_d$ becomes a polynomial map which is
defined for $u\in\A^1$ by
$$\omega_d(u):=\big(u^d-1,(u^d-1)u,\dots,(u^d-1)u^d\big).$$
Therefore $\omega_d$ is a finite morphism of algebraic varieties
which maps the affine space $\A^1$ onto its image, namely $\mo_d$.
Hence $\mo_d$ is a closed, rational (and hence irreducible),
$\Q$--definable curve contained in the affine ambient space
$\A^{d+1}$. The coordinate ring of the curve $\mo_d$ is
canonically isomorphic to the $\Q$--algebra
$\Q[U^d-1,(U^d-1)U,\dots,(U^d-1)U^d]$. Therefore, the encoding
$\omega_d:\A^1\to\mo_d$ of the object class $\mo_d$ is
$\Q$--definable, holomorphic and {\em robust}.

Let $M:=\{e^{\frac{2\pi i}{d}k};0\le k<d\}$ and denote by
$0:=(0\klk 0)$ the origin of the affine space $\A^{d+1}$. Observe
that the point $0\in\A^{d+1}$ belongs to the curve $\mo_d$,
because $\omega_d$ maps any point of $M$ onto the origin of
$\A^{d+1}$. Thus $\omega_d$ represents an {\em ambiguous} robust
encoding of the object class $\mo_d$. One verifies easily that the
point $0$ is the only (ordinary) singularity of the rational curve
$\mo_d$ and that this singularity can be resolved by a single
blowing up. Moreover $\omega_d$ induces an isomorphism between the
affine curves $\A^1\setminus M$ and $\mo_d\setminus\{0\}$.

Suppose now that there is given a $\Q$--definable data structure
$\md_d$ and a $\Q$--definable holomorphic encoding
$\sigma_d:\md_d\to\mo_d$. Let us denote the size of $\md_d$ by
$m_d$. Suppose furthermore that there is given a \Q--definable
polynomial map $\theta_d:\A^1\to\A^{m_d}$ with
$\theta_d(\A^1)\subset\md_d$ and $\sigma_d\circ\theta_d=\omega_d$.
We interpret the polynomial map $\theta_d$ as a branching--free
algorithm which transforms the encoding $\omega_d$ into the
unambiguous encoding $\sigma_d$ (see Section \ref{section3.3} for
a motivation of this notion of algorithm).

Although the object class $\mo_d$ admits an (ambiguous) robust
encoding by a data structure of size one, namely $\omega_d$, any
{\em unambiguous} holomorphic encoding $\sigma_d$ of $\mo_d$,
obtained by an algorithmic transformation of $\omega_d$, requires
a data structure of large size (of approximately the dimension of
the ambient space of the object class $\mo_d$). This is the
content of the following result:

\begin{proposition}
\label{remark1,3.1} Let notations and assumptions be as before.
Suppose that \linebreak $\sigma_d:\md_d\to\mo_d$ is an unambiguous
holomorphic encoding
of the object class $\mo_d$. Then the size $m_d$ of the data
structure $ \md_d$ satisfies the estimate $$m_d\ge d.$$
\end{proposition}

\begin{prf}
Let $0\le k_1<k_2<d$. Since $\sigma_d:\md_d\to\mo_d$ is injective
we deduce from $\omega_d(e^{\frac{2\pi
i}{d}k_1})=\omega_d(e^{\frac{2\pi i}{d}k_2})=0$ and from
$\sigma_d\circ\theta_d=\omega_d$ that $\theta_d(e^{\frac{2\pi
i}{d}k_1})=\theta_d(e^{\frac{2\pi i}{d}k_2})$ holds. Therefore
there exists a code $\alpha\in\md_d$ satisfying the condition
$\alpha=\theta_d(e^{\frac{2\pi i}{d}k})$ for any $0\le k<d$. The
encoding $\sigma_d:\md_d\to\mo_d$ is induced by a polynomial map
$\A^{m_d}\to\A^{d+1}$ which we also denote by $\sigma_d$.

Let $0\le k<d$. Denote by $(D\sigma_d)_\alpha$ the derivative of
this polynomial map in the point $\alpha\in\A^{m_d}$ and by
$\theta_d'(e^{\frac{2\pi i}{d}k})\in\A^{m_d}$ and
$\omega_d'(e^{\frac{2\pi i}{d}k})\in\A^{d+1}$ the derivatives of
the polynomial maps $\theta_d$ and $\omega_d$ in the point
$e^{\frac{2\pi i}{d}k}\in\A^1$. Observe that the encoding
$\omega_d$ is represented by the $(d+1)$--tuple of univariate
polynomials $\big(U^d-1,(U^d-1)U\klk(U^d-1)U^d\big)$. Deriving
this representation with respect to the parameter $U$ and
evaluating the result in the point $e^{\frac{2\pi i}{d}k}\in\A^1$,
we conclude that $$(de^{\frac{2\pi i}{d}kj};-1\le
j<d)=\omega_d'(e^{\frac{2\pi
i}{d}k})=(D\sigma_d)_\alpha\big(\theta_d'(e^{\frac{2\pi
i}{d}k})\big)$$ holds.

One sees easily that the matrix $A:=\big(d(e^{\frac{2\pi
i}{d}kj})\big)_{0\le k<d,-1\le j<d }$ has maximal rank $d$.
Indeed, the $d\times d$ submatrix of the matrix $A$ consisting of
the last $d$ columns of the matrix $A$ is nonsingular, because it
is the product of a $d\times d$ nonsingular diagonal matrix by a
$d\times d$ nonsingular Vandermonde matrix. Therefore the $d$
tangent vectors $\omega_d'(e^{\frac{2\pi i}{d}k})$, $0\le k<d$, of
the curve $\mo_d$ at the point $0$ are \C--linearly independent.
Since $(D\sigma_d)_\alpha:\A^{m_d}\to\A^{d+1}$ is a \C--linear
map, we conclude that the $d$ points $\theta_d'(e^{\frac{2\pi
i}{d}k})$, $0\le k<d$ of the \C--linear space $\A^{m_d}$ are
linearly independent too. This implies $m_d\ge d$.
\end{prf}\cqfd

We observe that the proof of Proposition \ref{remark1,3.1}
implies that the local embedding dimension of the curve
$\mo_d$ at the point $0$ (and hence the global embedding
dimension of $\mo_d$) is at least $d$. We are now going to apply
the conclusion of Proposition \ref{remark1,3.1} to the
arithmetic circuit complexity model.
\medskip

Let $\ma:=\Q[U]$ and
$\mathcal{B}_d:=\Q[U^d-1,(U^d-1)U\klk(U^d-1)U^d]$. Assume $d\ge
3$. Observe that $F_d=\sum_{0\le j\le d}(U^d-1)U^jY^j$ belongs to
the polynomial rings $\ma[Y]$ and $\mathcal{B}_d[Y]$ and that
$\mathcal{B}_d$ is isomorphic to the coordinate ring of the curve
$\mo_d$.

For $R\in\{\ma,\mb_d\}$ denote by $L_R(F_d)$ the minimal
non--scalar size of the totally division--free arithmetic circuits
with single input $Y$ which evaluate the polynomial $F_d$ using
only scalars from $R$. In case $d=2^{r+1}-1$ for some $r\in\N$,
one infers from the representation $$F_d=(U^d-1)\prod_{0\le k\le
r}\big(1+(UY)^{2^k}\big)$$ the estimate
$$L_\ma(F_d)=r+1=\log(d+1)$$ (here by $\log$ we denote the
logarithm to the base 2).

In a similar way one sees easily that $L_\ma(F_d)=O(\log d)$ holds
for arbitrary $d\in\N$. From the trivial lower bound
$L_\ma(F_d)\ge\log d$ (see \cite[Chapter 8, 8.1]{BuClSh97}) one
deduces finally that the functions $L_\ma(F_d)$ and $\log d$ have
the same asymptotic growth (in symbols: $L_\ma(F_d)=\Theta(\log
d)$).

Let us now analyze $L_{\mb_d}(F_d)$. Since $F_d$ is a polynomial
of degree $d$ in the variable $Y$ we deduce from \cite[Chapter 9,
Proposition 9.1]{BuClSh97} the estimate $L_{\mb_d}(F_d)\le
2\sqrt{d}$. Let $L_d:=L_{\mb_d}(F_d)$. Then there exists a totally
division--free circuit $\beta_d$ of non--scalar size $L_d$ with
single input $Y$ which evaluates the polynomial $F_d$ using only
scalars from $\mb_d$. From \cite[Chapter 9, Theorem 9.9]{BuClSh97}
we deduce that without loss of generality the circuit $\beta_d$
may be supposed to use only $m_d:=L_d^2+2L_d+2$ scalars
$\theta_1^{(d)}\klk\theta_{m_d}^{(d)}$ from
$\mathcal{B}_d:=\Q[U^d-1,(U^d-1)U\klk(U^d-1)U^d]$. Let
$\theta_d:\A^1\to\A^{m_d}$ be the polynomial map defined by
$\theta_d:=(\theta_1^{(d)}\klk\theta_{m_d}^{(d)})$ and let $\md_d$
be the image of $\theta_d$. Observe that $\md_d$ is a
$\Q$--constructible subset of $\A^{m_d}$. Again from \cite[Chapter
9, Theorem 9.9]{BuClSh97} we infer that there exists a polynomial
map $\sigma_d:\A^{m_d}\to\A^{d+1}$ which satisfies the condition
$$\sigma_d(\theta_1^{(d)}\klk\theta_{m_d}^{(d)})=\big(U^d-1,
(U^d-1)U\klk(U^d-1)U^d\big).$$ Thus we have
$\sigma_d\circ\theta_d=\omega_d$ and $\sigma_d(\md_d)=\mo_d$. Let
us denote the restriction of the polynomial map
$\sigma_d:\A^{m_d}\to\A^{d+1}$ to $\md_d$ by
$\sigma_d:\md_d\to\mo_d$. Since
$\theta_1^{(d)}\klk\theta_{m_d}^{(d)}$ are polynomials in the
coefficients $U^d-1,(U^d-1)U\klk(U^d-1)U^d$ of $F_d\in\mb_d[Y]$ we
conclude that $\sigma_d:\md_d\to\mo_d$ is an unambiguous
holomorphic encoding of the object class $\mo_d$. From
Proposition \ref{remark1,3.1} we deduce now
$L_d^2+2L_d+2=m_d\ge d$. This implies the lower bound
$L_{\mb_d}(F_d)=L_d\ge\sqrt{d-2}$. In summary, we obtain the
following complexity result:
\begin{corollary}
\label{remark2,3.1} Let notations be as before. Then we have
$L_\ma(F_d)=\Theta(\log d)$ and $L_{\mb_d}(F_d)=\Theta(\sqrt{d})$.
\end{corollary}

In terms of \cite{Heintz89} this result means that the
sequence of polynomials $\mathcal{F}:=(F_d)_{d\in\N}$ is easy to
evaluate in $\ma[Y]$, whereas $\mathcal{F}$ becomes difficult to
evaluate if we require that for any $d\in\N$ the univariate
polynomial $F_d\in\mb_d[Y]$ has to be computed by a totally
division--free arithmetic circuit whose scalars belong to
$\mb_d$. In conclusion, the evaluation complexity of a polynomial
depends strongly on the ring of scalars admitted.
\medskip

We are now going to describe another application of
Proposition \ref{remark1,3.1}. Let $m_d$ be a natural number
and let
$\gamma_d:=(\gamma_1^{(d)}\klk\gamma_{m_d}^{(d)})\in\Q^{m_d}$ be
an identification sequence of length $m_d$ for the Zariski closure
$C_d$ of the cone generated by the object class $\mo_d$ in the
$(d+1)$--dimensional \C--linear subspace of polynomials of $\C[X]$
having degree at most $d$. Observe that $C_d$ is a \Q--definable,
closed, irreducible subvariety of $\A^{d+1}$. Let
$\md_d:=\{\big(G(\gamma_1^{(d)})\klk
G(\gamma_{m_d}^{(d)})\big);G\in\mo_d\}$ and let
$\tau_d:\mo_d\to\md_d$ be the bijective map defined for
$G\in\mo_d$ by $\tau_d(G):=\big(G(\gamma_1^{(d)})\klk
G(\gamma_{m_d}^{(d)})\big)$. One sees easily that
$\tau_d:\mo_d\to\md_d$ is induced by a \Q--definable linear map
from $\A^{d+1}$ to $\A^{m_d}$. On the other hand, this linear map
induces an injective, homogeneous morphism from the cone $C_d$
into the affine space $\A^{d+1}$. From Lemma \ref{lemma4:encode}
we deduce now that this morphism is closed with respect to the
Zariski topologies of $C_d$ and $\A^{d+1}$. Therefore
$\md_d=\tau_d(\mo_d)$ is a closed, \Q--definable, irreducible
curve contained in $\A^{m_d}$ and $\tau_d:\mo_d\to\md_d$ is a
bijective, birational morphism of \Q--definable, irreducible
curves. Let $\sigma_d:\md_d\to\mo_d$ be the inverse map of
$\tau_d$. We consider $\md_d$ as a \Q--constructible data
structure of size $m_d$ and $\sigma_d$ as an encoding by
values of the object class $\mathcal{O}_d$ in the sense of
Section \ref{valencode}. In particular $\sigma_d$ is a
$\Q$--definable, continuous encoding. With these notations
we are able to state the following result:
\begin{corollary}
\label{remark2prima,3.1} Suppose that the encoding by values
$\sigma_d:\md_d\to\mo_d$ is holomorphic.
Then the size $m_d$
of the data structure $\md_d$ satisfies the
estimate
$$m_d\ge d.$$
\end{corollary}

\begin{prf}
Since $\sigma_d$ and $\tau_d$ are inverse morphisms of
\Q--definable, irreducible curves, there exists a polynomial map
$\theta_d:\A^1\to\A^{m_d}$ with
$\theta_d(u)=\tau_d\big(\omega_d(u)\big)$ for any $u\in\A^1$. This
implies $\theta_d(\A^1)\subset\md_d$. Moreover we have
$\sigma_d\circ\theta_d=\sigma_d\circ(\tau_d\circ\omega_d)
=\omega_d$. Since $\sigma_d:\md_d\to\mo_d$ is an unambiguous
\Q--definable holomorphic encoding of the object class $\mo_d$, we
deduce from Proposition \ref{remark1,3.1} that $m_d\ge d$ holds.
\end{prf}\cqfd

Corollary \ref{remark2prima,3.1} says that there exists a family
of object classes, namely $(\mo_d)_{d\in\N}$, encoded by a single
data structure of size one, namely $\A^1$, such that any {\em
holomorphic} encoding of these object classes by values becomes
necessarily large, namely of size at least $d$ for any object
class $\mo_d$. Nevertheless in view of Theorem
\ref{theorem:encode1}, the object class $\mathcal{O}_d$ admits a
{\em continuous} robust encoding of {\em constant} length (in fact
of length 2).
\medskip

From Corollary \ref{remark2prima,3.1} we infer the following
general result:
\begin{theorem}
\label{theorem:section3.1} Let $L,m$ be natural numbers and let
$\gamma=(\om\gamma m)\in\Q^m$ be an identification sequence for
the object class $\overline{W}_{L,1}$ of all univariate
polynomials over \C\ which have approximative non--scalar
sequential time complexity at most $L$. Let
$\md^*:=\{\big(F(\gamma_1)\klk
F(\gamma_m)\big);F\in\overline{W}_{L,1}\}$ and let
$\tau:\overline{W}_{L,1}\to\md^*$ be the bijective map defined by
$\tau(F):=\big(F(\gamma_1)\klk F(\gamma_m)\big)$. Then $\md^*$ is
a \Q--definable closed cone of $\A^m$ and
$\tau:\overline{W}_{L,1}\to\md^*$ is a \Q--definable, bijective,
finite morphism of algebraic varieties. Let
$\sigma:\md^*\to\overline{W}_{L,1}$ be the inverse map of $\tau$.
Consider $\md^*$ as a \Q--definable data structure of size $m$ and
suppose that $\sigma:\md^*\to\overline{W}_{L,1}$ is a
\Q--definable, holomorphic encoding by values of the object class
$\overline{W}_{L,1}$. Then the size $m$ of the data structure
$\md^*$ satisfies the estimate $m\ge 2^{cL}$ for a suitable
universal constant $c>0$.
\end{theorem}
\begin{prf} There exists a constant $c'>0$ such that $L(G)\le
c'\log d$ holds for any $d\in\N$ and any univariate polynomial $G$
belonging to the object class $\mo_d$ (recall that $L(G)$ denotes
the non--scalar time complexity of the polynomial $G$).

Let $d:=\lfloor 2^{\frac{L}{c'}}\rfloor:=max\{z\in\Z;z\le
2^{\frac{L}{c'}}\}$. Then we have $G\in\overline{W}_{L,1}$ for any
$G\in\mo_d$. Therefore $\gamma$ is an identification sequence for
the object class $\mo_d$. From Corollary
\ref{remark2prima,3.1} we deduce now $m\ge d\ge
2^{\frac{L}{c'}}-1$. Choose now any constant $c>0$ with
$2^{\frac{1}{c'}}-1\ge 2^c$. Then we have $m\ge 2^{cL}$.
\end{prf}\cqfd

One proves easily a similar complexity result for multivariate
polynomials. This question will be reconsidered in a  forthcoming
paper.

%-----------------------------------------------------------------
%-----------------------------------------------------------------
%-----------------------------------------------------------------
%-----------------------------------------------------------------

\subsection{Second paradigm.}
\label{section3.2} Let $n$ be a fixed natural number and let $T$,
$U_1,\ldots,U_n$ and $Y$ be indeterminates over $\Q$. Let
$U:=(U_1,\ldots,U_n)$. We are going to consider $T$,
$U_1,\ldots,U_n$ as parameters and $Y$ as variable. In the sequel
we shall use the following notation: for arbitrary natural numbers
$i$ and $j$ we shall denote by $[j]_i$ the $i$th digit of the
binary representation of $j$.  Let $P_n$ be the following
polynomial of $\Q[T,U,Y]$:
\begin{equation}
  \label{equation:aux}
  P_n(T,U,Y):=\prod_{j=0}^{2^n-1}\big(Y-
  (j+T\prod_{i=1}^nU_i^{[j]_i})\big).
\end{equation}

We observe that the dense representation of $P_n$ with respect to
the variable $Y$ takes the form
$$P_n(T,U,Y)=Y^{2^n}+B_1^{(n)}Y^{2^n-1}+\cdots +B_{2^n}^{(n)},$$
where $B_1^{(n)},\ldots,B_{2^n}^{(n)}$ are suitable polynomials of
$\Q[T,U]$.

Let $1\le k\le 2^n$. In order to determine the polynomial
$B_k^{(n)}$, we observe, by expanding the right hand side of
(\ref{equation:aux}), that $B_k^{(n)}$ collects the contribution
of all terms of the form
$$\prod_{h=1}^k\big(-(j_h+T\prod_{i=1}^nU_i^{[j_h]_i})\big)$$ with
$0\le j_1<\cdots<j_k\le 2^n-1$. Therefore the polynomial
$B_k^{(n)}$ can be expressed as follows:
\[\begin{array}{rcl}
B_k^{(n)}&=&\displaystyle\sum_{\quad \ 0\le j_1<\cdots<j_k< 2^n\ }
  \prod_{h=1}^k\big(-(j_h+T\prod_{i=1}^nU_i^{[j_h]_i})\big)\\
  &=&\displaystyle\sum_{\quad \ 0\le j_1<\cdots<j_k< 2^n\ }
  (-1)^k\prod_{h=1}^k\big(j_h+T\prod_{i=1}^nU_i^{[j_h]_i}\big).
\end{array}\]
Observe that for $0\le j_1<\cdots<j_k< 2^n$ the expression $$
\prod_{h=1}^k\big(j_h+T\prod_{i=1}^nU_i^{[j_h]_i}\big)$$ can be
rewritten as: $$j_1\cdots
j_k+T\big(\sum_{h=1}^kj_1\cdots\widehat{j_h}\cdots j_k
\prod_{i=1}^nU_i^{[j_h]_i}\big)+\mbox{terms of higher degree in
}T.$$ Therefore, we conclude that $B_k^{(n)}$ has the form:
\begin{equation}
\label{equation:aux-1}
\begin{array}{rcl}
  B_k^{(n)}&=&\displaystyle\sum_{\quad 0\le j_1<\cdots<j_k< 2^n\ }j_1\cdots
  j_k\\ \qquad \\
  &&\quad\displaystyle+\ T\left(\hspace*{-9pt}\sum_{\quad 0\le j_1<\cdots<j_k< 2^n\ }
    \sum_{h=1}^k j_1\cdots\widehat{j_h}\cdots
    j_k \prod_{i=1}^nU_i^{[j_h]_i}\right)  \\ \qquad \\
  &&\quad+\ \displaystyle\mbox{terms of higher degree in }T.
\end{array}
\end{equation}

Let us denote by $L_k^{(n)}$ the coefficient of $T$ in the
representation (\ref{equation:aux-1}), namely:
$$L_k^{(n)}:=\sum_{\quad 0\le j_1<\cdots<j_k< 2^n\ } \sum_{h=1}^k
j_1\cdots\widehat{j_h}\cdots j_k \prod_{i=1}^nU_i^{[j_h]_i}.$$ We
shall need the following technical result of \cite{GiHe01}. In
order to maintain this paper self--contained we are going to
reproduce its proof here.

\begin{lemma}
\label{lemma:independence} The polynomials
$L_1^{(n)},\ldots,L_{2^n}^{(n)}$ are $\Q$-linearly independent in
$\Q[U]$.
\end{lemma}

\begin{prf}
Let us abbreviate $N:=2^n-1$ and $L_1:=L_1^{(n)}\klk
L_{N+1}:=L_{2^n}^{(n)}$. We observe that for $1\le k\le {N+1}$ and
$0\le j\le N$ the coefficient $\ell_{k,j}$ of the monomial
$\prod_{i=1}^nU_i^{[j]_i}$ occuring in the polynomial $L_k$ can be
represented as $$\ell_{k,j}= \sum_{\stackrel{\scriptstyle\ 0\le
j_1<\cdots<j_{k-1}\le N\ }{\scriptstyle j_r\not= j\ {\rm for}\
r=1,\ldots,k-1}}j_1\cdots j_{k-1}.$$

\noindent \textit{\textbf{Claim:} For fixed $N$ and $k$, the
coefficient $\ell_{k,j}$ can be written as a polynomial expression
of degree exactly $k-1$ in the index $j$. Moreover, this
polynomial expression for $\ell_{k,j}$ has integer coefficients}.
\medskip

\noindent {\em Proof of the Claim.} We proceed by induction on the
index parameter $k$.

For $k=1$ we have $\ell_{1,j}=1$ for any $0\le j\le N$ and
therefore $\ell_{1,j}$ is a polynomial of degree ${k-1}=0$ in the
index $j$.

Let $1\le k\le {N+1}$. Assume inductively that $\ell_{k,j}$ is a
polynomial of degree exactly $k-1$ in the index $j$ and that the
coefficients of this polynomial are integers. We are now going to
show that $\ell_{k+1,j}$ is  a polynomial of degree exactly $k$ in
$j$ and that the coefficients of this polynomial are integers too.
Observe that
\[\begin{array}{rcl}
    \ell_{k+1,j}&=&\displaystyle
    \sum_{\stackrel{\ \scriptstyle 0\le j_1<\cdots<j_{k}\le N\
}{\scriptstyle
         j_r\not= j\ {\rm for}\ r=1,\ldots,k}}j_1\cdots
j_{k}\\
         \qquad \\
    &=&\displaystyle\sum_{\ 0\le j_1<\cdots<j_{k}\le N\ }j_1\cdots j_{k}-
    j\Bigg(\sum_{\stackrel{\ \scriptstyle 0\le j_1<\cdots<j_{k-1}\le N\
}{\scriptstyle
         j_r\not= j\ {\rm for}\ r=1,\ldots,k-1}}j_1\cdots
j_{k-1}\Bigg).
 \end{array}\]
holds. Since the term $$ \sum_{0\le  j_1<\cdots<j_{k}\le
N}j_1\cdots j_{k}$$ does not depend on $j$ and since by induction
hypothesis $$\ell_{k,j}=\sum_{\stackrel{\ \scriptstyle 0\le
j_1<\cdots<j_{k-1}\le N\ }{\scriptstyle j_r\not= j\ {\rm for}\
r=1,\ldots,k-1}}j_1\cdots j_{k-1}$$ is a polynomial of degree
exactly $k-1$ in $j$, we conclude that $\ell_{k+1,j}$ is a
polynomial of degree exactly $k$ in $j$. Moreover, the
coefficients of this polynomial are integers. This proves our
claim. \medskip

It is now easy to finish the proof of Lemma
\ref{lemma:independence}. By our claim there exist for arbitrary
$1\le k\le {N+1}$ integers $c_0^{(k)},\cdots,c_{k-1}^{(k)}$ with
$c_{k-1}^{(k)} \not = 0$ such that for any $0\le j\le N$ the
identity $\ell_{k,j}=c_0^{(k)}+\cdots+c_{k-1}^{(k)}j^{k-1}$ holds.
Hence for arbitrary $0\le k\le N$ there exist rational numbers
$\lambda_1^{(k)},\ldots,\lambda_{k+1}^{(k)}$ (not depending on
$j$) such for any $0\le j\le N$ the condition
$$j^k=\lambda_1^{(k)}\ell_{1,j}+\cdots+\lambda_{k+1}^{(k)}\ell_{k+1,j}$$
is satisfied (here we use the convention $0^{0}:=1$). This implies
for any index $0\le k\le N$ the polynomial identity
$$\lambda_1^{(k)}L_{1}+\cdots+\lambda_{k+1}^{(k)}L_{k+1}=
\sum_{0\le j\le N}j^k\prod_{i=1}^nU_i^{[j]_i}.$$

Hence for any $0\le k\le N$ the polynomial $Q_k:=\displaystyle
\sum_{0\le j\le N}j^k\prod_{i=1}^nU_i^{[j]_i}$ belongs to the
$\Q$--vector space generated by $L_1,\ldots,L_{N+1}$. On the other
hand, we deduce from the nonsingularity of the Vandermonde matrix
$\left( j^k \right)_{0 \leq k,j \leq N}$ that the polynomials
$Q_0,\ldots,Q_N$ are $\Q$-linearly independent. Therefore the
$\Q$--vector space generated by $L_1,\ldots,L_{N+1}$ in $\Q[U]$
has dimension ${N+1}=2^n$. This implies that $L_1,\ldots,L_{N+1}$
are $\Q$-linearly independent.
\end{prf}\cqfd

Let us now consider the object class of univariate polynomials
$\mo^{(n)}:=\{P_n(t,u,Y);t\in\A^1,u\in\A^n\}$ and the encoding
$\omega^{(n)}:\A^1\times\A^n\to\mo^{(n)}$ defined for $t\in\A^1$,
$u\in\A^n$ by $\omega^{(n)}(t,u):=P_n(t,u,Y)$. We are going to
analyze the object class $\mo^{(n)}$ and its encoding
$\omega^{(n)}$ in the same way as in Section \ref{section3.1}.

Representing the univariate polynomials belonging to $\mo^{(n)}$
by their coefficients, we identify the object class $\mo^{(n)}$
with the corresponding subset of the ambient space $\A^{2^n}$.
With this interpretation $\omega^{(n)}$ becomes a polynomial map
over \Q\ which is defined for $t\in\A^1$, $u\in\A^n$ by
$\omega^{(n)}(t,u):=\big(B_1^{(n)}(t,u)\klk
B_{2^n}^{(n)}(t,u)\big)$. Thus
$\omega^{(n)}:\A^1\times\A^n\to\mo^{(n)}$ is a \Q--definable
holomorphic encoding of the object class $\mo^{(n)}$. Let
$\beta^{(n)}:=(\beta_1^{(n)}\klk\beta_{2^n}^{(n)})$ with
$\beta_k^{(n)}:=\sum_{0\le j_1<\cdots<j_k< 2^n\ }j_1\cdots j_k$
for $1\le k\le 2^n$. From (\ref{equation:aux-1}) one deduces
immediately that $\beta^{(n)}$ belongs to $\mo^{(n)}$ and that
$P_n(0,u,Y)=Y^{2^n}+\beta_1^{(n)}Y^{2^n-1}\plp\beta_{2^n}^{(n)}$
holds for any $u\in\A^n$. Hence the fiber
$(\omega^{(n)})^{-1}(\beta^{(n)})$ contains the hyperplane
$\{0\}\times \A^n$ of the affine space $\A^1\times\A^n$. This
implies that the encoding $\omega^{(n)}$ is ambiguous and {\em not
robust}.

Suppose now that there is given a \Q--definable, holomorphic
encoding $\sigma^{(n)}:\md^{(n)}\to\mo^{(n)}$. Let us denote the
size of $\md^{(n)}$ by $m^{(n)}$. Suppose furthermore that there
is given a \Q--definable polynomial map
$\theta^{(n)}:\A^1\times\A^n\to\A^{m^{(n)}}$ with
$\theta^{(n)}(\A^n)\subset\md^{(n)}$ and
$\sigma^{(n)}\circ\theta^{(n)}=\omega^{(n)}$. As before, we
interpret the polynomial map $\theta^{(n)}$ as a branching--free
algorithm which transforms the encoding $\omega^{(n)}$ into the
encoding $\sigma^{(n)}$. Although the object class $\mo^{(n)}$
admits a (non--robust) encoding of small  size (i.e. small in
comparison with the embedding dimension of the object class
$\mo^{(n)}$), the requirement of {\em robustness} for the encoding
$\sigma^{(n)}$ entails that the size of the data structure
$\md^{(n)}$ must be necessarily large. This is the content of the
following result.
\begin{proposition}
\label{remark3,3.2} Let notations and assumptions be as before.
Suppose that \linebreak $\sigma^{(n)}:\md^{(n)}\to\mo^{(n)}$ is a
robust, holomorphic encoding of the object class
$\mo^{(n)}$. Then the size of the data structure $\md^{(n)}$
satisfies the estimate $$m^{(n)}\ge 2^n.$$
\end{proposition}

\begin{prf}
Since $\sigma^{(n)}$ is a robust encoding we conclude that
$(\sigma^{(n)})^{-1}(\beta^{(n)})$ is a nonempty finite subset of
$\md^{(n)}$. From
$\{0\}\times\A^n\subset(\omega^{(n)})^{-1}(\beta^{(n)})$ and
\linebreak
$\omega^{(n)}=\sigma^{(n)}\circ\theta^{(n)}$ we infer
$\theta^{(n)}(\{0\}\times\A^n)\subset(\sigma^{(n)})^{-1}(\beta^{(n)})$.
Since $(\sigma^{(n)})^{-1}(\beta^{(n)})$ is finite and
$\{0\}\times\A^n$ irreducible there exists a point
$\alpha\in(\sigma^{(n)})^{-1}(\beta^{(n)})$ with
$\theta^{(n)}(\{0\}\times\A^n)=\{\alpha\}$. Let $u$ be arbitrary
point of $\A^n$ and let $\gamma_u:\A^1\to\A^{m^{(n)}}$ and
$\delta_u:\A^1\to\A^{2^n}$ be the polynomial maps defined for
$t\in\A^1$ by $\gamma_u(t):=\theta^{(n)}(t,u)$ and
$\delta_u(t):=\omega^{(n)}(t,u)$. Then we have
$\gamma_u(0)=\alpha$, $\delta_u(0)=\beta^{(n)}$ and
$\sigma^{(n)}\circ\gamma_u=\delta_u$. From (\ref{equation:aux-1})
we deduce now $$\big(L_1^{(n)}(u)\klk
L_{2^n}^{(n)}(u)\big)=\frac{\partial}{\partial
t}\omega^{(n)}(0,u)=\delta_u'(0)=(D\sigma^{(n)})_\alpha
\big(\gamma_u'(0)\big).$$ Lemma \ref{lemma:independence} implies
that there exist points $u_1\klk u_{2^n}\in\A^n$ such that the
$(2^n\times 2^n)$--matrix $\big(L_i^{(n)}(u_j)\big)_{1\le i,j\le
2^n}$ is nonsingular. Therefore
$\delta_{u_1}'(0)\klk\delta_{u_{2^n}}'(0)$ are linearly
independent elements of the \C--vector space $\A^{2^n}$. Since
$(D\,\sigma^{(n)})_\alpha:\A^{m^{(n)}}\to\A^{2^n}$ is a \C--linear
map, we conclude that $\gamma_{u_1}'(0)\klk\gamma_{u_{2^n}}'(0)$
are linearly independent elements of the \C--linear space
$\A^{m^{(n)}}$. This implies $m^{(n)}\ge 2^n$.
\end{prf}\cqfd

We observe that the proof of Proposition \ref{remark3,3.2}
implies that the local embedding dimension of the closed
algebraic variety $\omo$ at the point $\beta^{(n)}$ (and hence
the global embedding dimension of $\mo^{(n)}$) is exactly $2^n$.
\medskip

Let $\mathcal{B}^{(n)}:=\Q[B_1^{(n)}\klk B_{2^n}^{(n)}]$ and let
us denote by $L_{\mathcal{B}^{(n)}}(P_n)$ the minimal non--scalar
size of the totally division--free arithmetic circuit with single
input $Y$ which evaluates the polynomial $P_n$ using only scalars
belonging to the \Q--algebra $\mathcal{B}^{(n)}$. In the same way
as in Section \ref{section3.1} we may deduce from Proposition
\ref{remark3,3.2} the following result:
\begin{corollary}
\label{remark4,3.2} With the notations introduced before we have
$$L_{\mathcal{B}^{(n)}}(P_n)=\Theta(2^{\frac{n}{2}}).$$
\end{corollary}

Corollary \ref{remark4,3.2} says that the sequence of
polynomials $(P_n)_{n\in\N}$ becomes hard to evaluate, if we
require that for any $n\in\N$ the univariate polynomial
$P_n\in\mathcal{B}^{(n)}[Y]$ has to be evaluated by a totally
division--free arithmetic circuit whose scalars belong only to the
\Q--algebra $\mathcal{B}^{(n)}$.

%-----------------------------------------------------------------
%-----------------------------------------------------------------
%-----------------------------------------------------------------
%-----------------------------------------------------------------

\subsection{Rationality considerations.}
\label{section3.3} In this section we motivate the algorithmic
model used in Sections \ref{section3.1} and \ref{section3.2} for
the algorithmic transformation of encodings of a given object
class. For this purpose we are going to discuss the effect of
certain rationality conditions on the encoding of an object class.
Our first rationality condition requires to fix not only the
ground field, namely \Q, but also its algebraic closure, namely
\C.

Let be given a data structure $\md$ and an object class $\mo$ and
suppose that $\md$ and $\mo$ are \Q--constructible subsets of the
ambient spaces $\A^L$ and $\A^N$ respectively. Let be given an
encoding $\omega:\md\to\mo$ and suppose that $\omega$ is
\Q--definable and holomorphic. Let us denote by $m$ the maximal
local embedding dimension of the \Q--Zariski closure $\omo$ of the
object class $\mo$ at any point of $\omo$ (i.e. $m$ is the maximal
\C--vector space dimension of the Zariski tangent space of the
algebraic variety $\omo$ at any point). \bigskip

The first rationality condition we are going to consider is the
following: \smallskip

\noindent {\em for any object $\beta=(\beta_1\klk \beta_N)\in\mo$
there exists a code
$\alpha=(\alpha_1\klk\alpha_L)\in\omega^{-1}(\beta)$ with
$\alpha_1\klk\alpha_L\in\Q[\beta]$.} \smallskip

Suppose now that $\omega$ satisfies this rationality condition,
that $\md$ and $\mo$ are \Q--definable closed subvarieties of
$\A^L$ and $\A^N$ and that $\mo$ is \Q--irreducible. Since the
transcendence degree of \C\ over \Q\ is infinite, there exists a
generic element $b=(b_1\klk b_N)$ of $\mo$ such that the canonical
specialization of the coordinate ring $\Q[\mo]$ of the irreducible
algebraic variety $\mo$ onto $\Q[b]$ is injective. Therefore we
have $\Q[\mo]\cong\Q[b]$. By hypothesis there exists a code
$a=(\om a L)\in\omega^{-1}(b)$ with $\om a L\in\Q[b]$. Denote by
$Y_1\klk Y_N$ the coordinate functions of the affine space $\A^N$.
Then there exist polynomials $\om\psi L\in\Q[Y_1\klk Y_N]$ with
$a_k=\psi_k(b_1\klk b_N)$ for $1\le k\le L$.

Since $\mo$ and $\md$ are closed subvarieties of $\A^N$ and $\A^L$
respectively, this implies that $\psi:=(\om\psi L)$ induces a
\Q--definable morphism of algebraic varieties which maps $\mo$
into $\md$ and which we denote by $\psi:\mo\to\md$. Moreover we
have $\omega\circ\psi=id_\mo$. Therefore for any $\beta\in\mo$ the
\C--linear map $T_{\psi(\beta)}(\omega):T_{\psi(\beta)}(\md)\to
T_\beta(\mo)$ is surjective (here $T_{\psi(\beta)}(\md)$ and
$T_\beta(\mo)$ denote the Zariski tangent spaces of the algebraic
varieties $\md$ and $\mo$ at the points $\psi(\beta)$ and $\beta$
respectively and $T_{\psi(\beta)}(\omega)$ denotes the tangent map
induced by $\omega$ at the point $\psi(\beta)$).

Since the object class $\mo$ is irreducible and Zariski closed in
$\A^N$, we may choose a point $\beta_0\in\mo$ with $dim_\C
T_{\beta_0}(\mo)=m$ (here $dim_\C T_{\beta_0}(\mo)$ denotes the
\C--vector space dimension of the Zariski tangent space
$T_{\beta_0}(\mo)$). Since
$T_{\psi(\beta_0)}(\omega):T_{\psi(\beta_0)}(\md)\to
T_{\beta_0}(\mo)$ is a surjective \C--linear map, we conclude
$L\ge dim_\C T_{\psi(\beta_0)}(\md)\ge m$. Therefore any encoding
of the data structure $\mo_d$ of Section \ref{section3.1} or of
the data structure $\mo^{(n)}$ of Section \ref{section3.2} which
satisfies the rationality condition above has at least size $d$ or
size $2^n$ respectively.
\medskip

Now we are going to discuss a second rationality condition which
comes much closer to the usual requirements in the design of
practical algorithms. For the sake of succinctness of exposition
we shall omit proofs (they are based on Hilbert's Irreducibility
Theorem and L\"uroth's Theorem and will be published in
forthcoming paper).
\medskip

Informally, we may state our second rationality requirement as
follows:

\smallskip \noindent {\em suppose that the object class $\mo$
contains ``many" integer objects (i.e. points which belong to
$\Z^N$). Then there exist ``sufficiently many" integer objects of
$\mo$ such that for each such object $O\in\mo\cap\Z^N$ there
exists an integer code $D\in\md\cap\Z^L$ with $\omega(D)=O$. In
order to guarantee the existence of sufficiently many integer
objects in $\mo$ we require that there is given an encoding
$\omega^*:\A^{L^*}\to\mo$ of the object class $\mo$ such that
$\omega^*$ is definable by polynomials with integer coefficients.
Thus $\omega^*$ maps integer codes of $\A^{L^*}$ onto integer
objects of $\mo$.}
\smallskip

For technical reasons we shall need the following additional
assumptions:

\smallskip \noindent we suppose that $\md$ is a \Q--definable,
\Q--irreducible subvariety of $\A^L$ and that the given encoding
$\omega:\md\to\mo$ is definable by polynomials with integer
coefficients. This implies that the closed subvariety $\omo$ of
$\A^N$ is \Q--definable and \Q--irreducible too. Moreover we
suppose $dim\,\md=dim\,\omo$. Therefore there exists a Zariski open
subset of $\omo$ which is contained in $\mo=\omega(\md)$, such
that each point of this subset has a nonempty, finite
$\omega$--fiber. With these notations and assumptions we are able
to state the following result:

\begin{proposition}
\label{proposition:section4.3} Suppose that there exists a
nonempty \Q--definable, Zariski open subset $\mau$ of $\A^{L^*}$
with the following property: for any code $D^*\in\mau\cap\Z^{L^*}$
there exists a code $D\in\md\cap\Z^L$ with
$\omega^*(D^*)=\omega(D)$.

Then there exists a \Q--definable morphism $\theta:\A^{L^*}\to\md$
and a subset $\mau_0$ of $\mau\cap\Z^{L^*}$ such that the
following conditions are satisfied:
\begin{itemize}
  \item[($i$)] $\mau_0$ is Zariski--dense in $\A^{L^*}$.
  \item[($ii$)] $\omega\circ\theta=\omega^*$.
  \item[($iii$)] $\theta(D^*)\in\md\cap\Z^L$ for any
  $D^*\in\mau_0$.
\end{itemize}
Suppose additionally $L^*:=1$ and that the leading coefficients of
all nonconstant polynomials occurring in the definition of
$\omega^*$ have greatest common divisor one. Suppose furthermore
that for any integer object $O\in\mo\cap\Z^N$ there exists an
integer code $D\in\md\cap\Z^L$ with $O=\omega(D)$. Then $\omega$
is a robust encoding.
\end{proposition}

We may paraphrase the first part of Proposition
\ref{proposition:section4.3} as follows:

\smallskip \noindent if the encoding $\omega:\md\to\mo$ admits for
any object of $\mo$, which allows an {\em integer} encoding by
$\omega^*$, an {\em integer} encoding by $\omega$, then the
encoding $\omega^*$ may be transformed into the encoding $\omega$
by means of an algorithm in the sense of Sections \ref{section3.1}
and \ref{section3.2}. This motivates the notion of
(branching--free) algorithm which we introduced before and which
we shall continue to use in the remaining part of this paper.

The second part of Proposition \ref{proposition:section4.3} says
roughly the following: if a given \Q--definable holomorphic
encoding of an infinite object class by a data structure of size
one satisfies our second rationality condition, then this encoding
is necessarily robust. Therefore the paradigm of Section
\ref{section3.1} is representative for this type of encodings.

%-----------------------------------------------------------------
%-----------------------------------------------------------------
%-----------------------------------------------------------------
%-----------------------------------------------------------------
%-----------------------------------------------------------------
%-----------------------------------------------------------------
%-----------------------------------------------------------------
%-----------------------------------------------------------------

\section{The complexity of elimination algorithms.}
\label{section4}

%-----------------------------------------------------------------
%-----------------------------------------------------------------
%-----------------------------------------------------------------
%-----------------------------------------------------------------

\subsection{Flat families of zero--dimensional elimination
problems.} \label{section4.1} Let, as before, $k$ be an infinite
and perfect field with algebraic closure $\bar{k}$ and let $\om U
r, \xon, Y$ be indeterminates over $k$. In the sequel we shall
consider $\xon$ and $Y$ as variables and $\om U r$ as parameters.
Let $U:=(\om U r)$ and $X:=(\om X n)$ and let $\om G n$ and $F$ be
polynomials belonging to the $k$-algebra $k[U,X]:=k[\om U r,
\xon]$. Suppose that the polynomials $\om G n$ form a regular
sequence in  $k[U,X]$ defining thus an equidimensional subvariety
$V:=\{ G_1=0\klk G_n=0\}$ of the $(r+n)$--dimensional affine space
$\A^{r} \times \A^{n}$. The algebraic variety $V$ has dimension
$r$. Let $\delta$ be the (geometric) degree of $V$ (observe that
this degree does not take into account multiplicities or
components at infinity). Suppose furthermore that the morphism of
affine varieties $\pi:V\rightarrow \A^{r}$, induced by the
canonical projection of $ \A^{r} \times \A^{n}$ onto $ \A^{r}$, is
finite and generically unramified (this implies that $\pi$ is flat
and that the ideal generated by $\om G n$ in $k[U,X]$ is radical).
Let $\tilde{\pi}:V\rightarrow \A^{r+1}$ be the morphism defined by
$\tilde{\pi}(z) := (\pi(z),F(z))$ for any point $z$ of the variety
$V$. The image of $\tilde{\pi}$ is a hypersurface of $\A^{r+1}$
whose minimal equation is a polynomial of $k[U,Y]:=k[\om U r, Y]$
which we denote by $P$. Let us write $\deg P$ for the total degree
of the polynomial $P$ and $\deg_Y P$ for its partial degree in the
variable $Y$. Observe that $P$ is monic in $Y$ and that $\deg
P\leq \delta \deg F$ holds. Furthermore, for a Zariski dense set
of points $u$ of $\A^{r}$, we have that $\deg _Y P$ is the
cardinality of the image of the restriction of $F$ to the finite
set $\pi^{-1}(u)$. The polynomial $P(U,F)$ vanishes on the variety
$V$.\medskip

Let us consider an arbitrary point $u:=(u_1,\ldots ,u_r)$ of $\A
^r$. For arbitrary polynomials $A\in k[U,X]$ and $B\in k[U,Y]$ we
denote by $A^{(u)}$ and $B^{(u)}$ the polynomials $A(u_1\klk
u_r,\xon)$ and $B(u_1\klk u_r,Y)$ which belong to
$k(u)[X]:=k(u_1\klk u_r)[\xon]$ and $k(u)[Y]:=k(u_1\klk u_r)[Y]$
respectively. Similarly we denote for an arbitrary polynomial
$C\in k[U]$ by $C^{(u)}$ the value $C(u_1,\ldots ,u_r)$ which
belongs to the field $k(u):=k(u_1,\ldots ,u_r)$. The polynomials
$G_1^{(u)}\klk G_n^{(u)}$ define a zero dimensional subvariety
$V^{(u)}:= \{G_1^{(u)}=0\klk G_n^{(u)}=0\}= \pi^{-1}(u)$ of the
affine space $\A^{n}$. The degree (i.e. the cardinality) of
$V^{(u)}$ is bounded by $\delta$. Denote by $\tilde{\pi}^{(u)}:\
V^{(u)} \rightarrow \A^{1}$ the morphism induced by the polynomial
$F^{(u)}$ on the variety $V^{(u)}$. Observe that the polynomial
$P^{(u)}$ vanishes on the (finite) image of the morphism
$\tilde{\pi}^{(u)}$. Observe also that the polynomial $P^{(u)}$ is
not necessarily the minimal equation of the image of
$\tilde{\pi}^{(u)}$.\medskip

We call the equation system $G_1 =0\klk G_n =0$ and the
polynomial $F$ {\em a flat family of zero--dimensional
elimination problems depending on the parameters $\om U r$} and
we call $P$ the associated {\em elimination polynomial}. An
element $u\in\A^r$ is considered as a {\em parameter point} which
determines a {\em particular problem instance}. The equation
system $G_1 =0\klk G_n =0$ together with the polynomial $F$ is
called the {\em general instance} of the given flat family of
elimination problems and the elimination polynomial $P$ is called
the {\em general solution} of this flat family. A branching--free
algorithm which in terms of suitable data structures computes
from a given representation of the general problem instance
$G_1=0\klk G_n=0,F$ a representation of its general solution $P$
is called a {\em Kronecker--like elimination procedure} (see
Section \ref{section1.2}).\medskip

The {\em particular problem instance} determined by the parameter
point $u\in \A^r$ is given by the equations $G_1^{(u)}=0\klk
G_n^{(u)}=0$ and the polynomial $F^{(u)}$. The polynomial
$P^{(u)}$ is called {\em a solution} of this particular problem
instance. We call two parameter points $u,u' \in \A^r$ {\em
equivalent} (in symbols: $u \sim u'$) if $G_1^{(u)}=G_1^{(u')}
\klk G_n^{(u)}=G_n^{(u')}$ and $F^{(u)} = F^{(u')}$ holds.
Observe that $u \sim u'$ implies $P^{(u)} = P^{(u')}$. We call
polynomials $A\in k[U,X]$, $B\in k[U,Y]$ and $C\in k[U]$ {\em
invariant} (with respect to $\sim$) if for any two parameter
points $u,u'$ of $\A^r$ with $u\sim u'$ the respective identities
$A^{(u)} = A^{(u')}$, $B^{(u)} = B^{(u')}$
and $C^{(u)}=C^{(u')}$ hold.\medskip

Let us consider the set of parameter points of $\A^r$ as data
structure which encodes the object class $\mo:=\{(G_1^{(u)}\klk
G_n^{(u)},F^{(u)});u\in\A^r\}$. The corresponding encoding
$\omega:\A^r\to\mo$ is defined for $u\in\A^r$ by
$\omega(u):=(G_1^{(u)}\klk G_n^{(u)},F^{(u)})$. Observe that
$\omega$ is \Q--definable and holomorphic. Let $\md^*$ be a
$k$--constructible data structure of size $L^*$ which encodes the
object class $\mo^*:=\{P^{(u)};u\in\A^r\}$ by means of a given
$k$--definable, holomorphic encoding $\omega^*:\md^*\to\mo^*$. Let
us consider $\mo$ as input and $\mo^*$ as output object class of a
given branching--free Kronecker--like elimination procedure.
Suppose that this elimination procedure is determined by
polynomials $\theta_1\klk\theta_{L^*}\in k[U]$ such that
$\theta:=(\om\theta {L^*})$ induces a $k$--definable map from
$\A^r$ into $\md^*$ which we denote by $\theta:\A^r\to\md^*$. This
means that for any parameter point $u\in\A^r$, the given
elimination procedure, which we denote also by $\theta$, satisfies
the condition $\omega^*\big(\theta(u)\big)=P^{(u)}$. Suppose
furthermore that the elimination procedure $\theta$ is totally
{\em division--free} (this means that the general solution $P$ of
the given elimination problem belongs to $k[\theta][Y]$; see
Section \ref{models}). We call $\theta$ {\em invariant} (with
respect to the equivalence relation $\sim$) if $\om\theta{L^*}$
are invariant polynomials. The invariance of the elimination
procedure $\theta$ means that for any input code $u\in\A^r$ the
code $\theta(u)\in\md^*$ of the corresponding output object
$P^{(u)}$ depends only on the input {\em object}, namely
$(G_1^{(u)}\klk G_n^{(u)},F^{(u)})\in\big(k(u)[X]\big)^{n+1}$ and
not on its particular representation $u$. Said otherwise, an
invariant elimination procedure produces the solution of a
particular problem instance in a way which is independent of the
possibly different representations of the given problem
instance.\medskip

Since all known Kronecker--like elimination procedures produce for
flat families of zero--dimensional elimination problems a
branching and totally division--free representation of the output
polynomial, and since they are based on the manipulation of the
input objects (and not on their particular representations) by
means of linear algebra or comprehensive Gr\"obner basis
techniques, we conclude that these algorithms are in fact
invariant elimination procedures.

Typical examples of such procedures are furnished by black--box
algorithms. With the notations introduced before, we call the
elimination procedure $\theta$ a {\em black--box algorithm} if
for any input code $u\in\A^r$, the procedure $\theta$ calls only
for evaluations of the input object $(G_1^{(u)}\klk
G_n^{(u)},F^{(u)})$ on specializations of the variables $\xon$ to
assignment values which belong to suitable commutative
$k[u]$--algebras.

A (branching--parsimonious) computer program for elimination
tasks which calls its input polynomials only by their
specification as evaluation procedures, represents necessarily a
black--box algorithm.\medskip

We are now going to introduce a slight generalization of the
notion of invariance of the elimination procedure $\theta$.

Let $\md_\theta:=\left\{\big(\omega(u),\theta(u)\big);u\in\A^r\right\}$ and
let $\omega_\theta:\md_\theta\to\mo$ be the canonical first
projection of $\md_\theta$ onto the object class $\mo$. One
verifies immediately that $\md_\theta$ is a \Q--definable data
structure and that $\omega_\theta$ is a \Q--definable holomorphic
encoding of the object class $\mo$. Observe that
$\overline{\md}_\theta$ and $\omo$ are irreducible closed
subvarieties of their corresponding affine ambient spaces.

\begin{definition} \label{robustel}
Let notations and assumptions be as before. We call the
elimination procedure $\theta$ robust if the following condition
is satisfied:
\smallskip

\noindent let $u\in\A^r$ be a given parameter point determining
the input object \linebreak $(G_1^{(u)}\klk
G_n^{(u)},F^{(u)})\in\mo$ and let $\mathfrak{m}_u$ be the maximal
defining ideal of this input object in $\C[\omo]$. Then the local
ring $\C[\overline{\md}_\theta]_{\mathfrak{m}_u}$ is a finite
$\C[\omo]_{\mathfrak{m}_u}$--module.
\end{definition}

In other words,  the elimination procedure $\theta$ is robust
if and only if $\omega_{\theta}$ is a robust holomorphic encoding.

Observe that an invariant elimination procedure is
robust.

If $\theta$ is a robust elimination procedure, then one sees
easily that the following condition is satisfied:

\begin{itemize}
  \item[$(i)$]{\em for any parameter point $u\in\A^r$, the set
  $$\{
  \theta(v);v\in\A^r,G_1^{(v)}=G_1^{(u)},\dots,
  G_n^{(v)}=G_n^{(u)},F^{(v)}=F^{(u)}\}$$ is finite.}
\end{itemize}

In case $k:=\Q$ and $\ok:=\C$, one deduces easily from Lemma
\ref{equivrob} that Definition \ref{robustel} is equivalent to
the following condition:

\begin{itemize}
\item[$(ii)$]{\em let $\iseq u$ be a sequence of parameter points of $\A^r$
encoding a sequence of input objects $\big((G_1^{(u_i)}\klk
G_n^{(u_i)},F^{(u_i)})\big)_{i\in\N}$. Suppose that there
exists a parameter point $u\in\A^r$ such that $(G_1^{(u)}\klk
G_n^{(u)},F^{(u)})\in\mo$ is an limit point of the sequence of
input objects $\big((G_1^{(u_i)}\klk
G_n^{(u_i)},F^{(u_i)})\big)_{i\in\N}$ (with respect to the
strong topology). Then the sequence $(\theta(u_i))_{i\in\N}$ has
an accumulation point.}
\end{itemize}

Observe that condition $(ii)$ gives an intuitive meaning to
the technical Definition \ref{robustel}.

%-----------------------------------------------------------------
%-----------------------------------------------------------------
%-----------------------------------------------------------------
%-----------------------------------------------------------------

\subsection{Parametric greatest common divisors and their computation.}
\label{section5.2} Let us now introduce the notion of parametric
greatest common divisor of a given algebraic family of polynomials
and let us consider the corresponding algorithmic problem. We are
going to use the same notations as in Sections \ref{models}
and \ref{section4.1}.

Suppose that there is given a
positive number $s$ of nonzero polynomials, say
${B_1,\ldots,B_s}\in k[{U_1,\ldots,U_r},Y]$. Let $V:= \{B_1
=0,\ldots,B_s =0 \}$. Suppose that $V$ is nonempty. We consider
now the morphism of affine varieties $\pi:V\longrightarrow
\A^{r}$, induced by the canonical projection of $ \A^{r} \times
\A^{1}$ onto $\A^{r}$. Let $S$ be the Zariski closure of $\pi(V)$
and suppose that $S$ is an {\em irreducible\/} closed subvariety
of $ \A^{r}$. Let us denote by $k[S]$ the coordinate ring of $S$.
Since $S$ is irreducible we conclude that $k[S]$ is a domain with
a well defined function field which we denote by $k(S)$.

Let ${b_1,\ldots,b_s}\in k[S][Y]$ be the polynomials in the
variable $Y$ with coefficients in $k[S]$, induced by
${B_1,\ldots,B_s}$. Suppose that there exists an index $1\leq
k\leq s$ with $b_k \neq 0$. Without loss of generality we may
suppose that for some index $1\leq q\leq s$ the polynomials
${b_1,\ldots,b_q}$ are exactly the non--zero elements of
${b_1,\ldots,b_s}$. Observe that each polynomial
${b_1,\ldots,b_q}$ has positive degree (in the variable $Y$).

We consider ${b_1,\ldots,b_q}$ as an {\em algebraic family of
polynomials\/} (in the variable $Y$) and ${B_1,\ldots,B_q}$ as
their representatives. The polynomials ${b_1,\ldots,b_q}$ have in
$k(S)[Y]$ a well defined {\em normalized\/} (i.e. monic) greatest
common divisor, which we denote by $h$. Let $D$ be the degree of
$h$ (with respect to the variable $Y$).\medskip

We are now going to describe certain geometric requirements which
will allow us to consider $h$ as a parametric greatest common
divisor of the algebraic family of polynomials ${b_1,\ldots,b_q}$.

Our first requirement is $D\geq 1$. Moreover we require that for
any point $u\in S$ and any place
$\varphi:\overline{k}(S)\to\overline{k}\cup\{\infty\}$, whose
valuation ring contains the local ring of the variety $S$ at the
point $u$, the values of the coefficients of the polynomial $h\in
k(S)[Y]$ under $\varphi$ are {\em finite} and {\em uniquely
determined} by the point $u$. In this way the place $\varphi$ maps
the polynomial $h$ to a monic polynomial of degree $D$ in $Y$ with
coefficients in $\overline{k}$. This polynomial depends only on
the point $u\in S$ and we denote it therefore by $h(u)(Y)$. In
analogy with this notation we write $b_k(u)(Y):= B_k(u)(Y)$ for
$1\leq k\leq q$. Since $h$ is monic one concludes easily that
$h(u)(Y)$ divides the polynomials $b_1(u)(Y)\klk b_q(u)(Y)$ (and
hence their greatest common divisor if not all of them are
zero).\medskip

We say that a polynomial $H$ of $k({U_1,\ldots,U_r})[Y]$ with
${deg_Y}H =D$ {\em represents\/} the greatest common divisor
$h\in {k(S)[Y]}$ if the coefficients of $H$ with respect to the
variable $Y$ induce well--defined rational functions of the
variety $S$ and if these rational functions are exactly the
coefficients of $h$ (with respect to the variable $Y$).\medskip

Suppose now that the polynomials ${B_1,\ldots,B_s}\in
k[{U_1,\ldots,U_r},Y]$ satisfy all our requirements for any point
$u\in S$. Then we say that for the algebraic family of polynomials
$b_1\klk b_q\in k(S)[Y]$ a parametric common divisor exists and we
call $h\in k(S)[Y]$ the {\em parametric greatest common divisor}
of $b_1\klk b_q$. Any polynomial $H\in k({U_1,\ldots,U_r})[Y]$
which represents $h$ is said to {\em represent the parametric
greatest common divisor associated to the polynomials
$B_1,\ldots,B_s$\/}.\medskip

A monic {\em squarefree\/} polynomial $\widehat{h}\in {k(S)[Y]}$
with the same zeroes as $h$ in an algebraic closure of $k(S)$, is
called the {\em generically squarefree\/} parametric greatest
common divisor of the algebraic family ${b_1,\ldots,b_q}\in
k[S][Y]$ if $\widehat{h}$ satisfies the requirements imposed
above on $h$. In this case we say that for the algebraic family
of polynomials $b_1\klk b_q\in k[S][Y]$ a generically squarefree
parametric greatest common divisor exists. The notion of a
representative of $\widehat{h}$ is defined in the
same way as for $h$.\medskip

Let us consider $\A^r$ as input data structure of size $r$ with
$S$ the set of admissible input instances and let us consider the
problem of computing the parametric greatest common divisor $h$
by means of an essentially division--free algorithm for any
admissible input instance $u\in S$.

Such an algorithm, with output data structure of size $m$, is
represented by $m$ rational functions $\theta_1\klk\theta_m\in
k(S)$ such that the parametric greatest common divisor $h$ belongs
to the $k$--algebra $k[\theta_1\klk\theta_m][Y]$. Recall that our
assumptions on $h$ imply that for any input instance $u\in S$ the
polynomial $h(u)(Y)\in\overline{k}[Y]$ is well defined.
Consequently we shall require that for any input instance $u\in S$
and any place $\varphi:k(S)\to\overline{k}\cup\{\infty\}$, whose
valuation ring contains the local ring of the variety $S$ at the
point $u$, the values $\varphi(\theta_1)\klk\varphi(\theta_m)$ are
finite and uniquely determined by the input instance $u$. If this
requirement is satisfied we shall say that our algorithm {\em
computes} the parametric greatest common divisor $h$ of the
algebraic family of polynomials $b_1\klk b_q$ {\em for any
admissible input instance} $u\in S$. Observe that in this case the
rational functions $\theta_1\klk\theta_m$ belong to the integral
closure of $k[S]$ in $k(S)$.\medskip

In concrete situations it is reasonable, however not required by
the mathematical arguments we will apply in this paper, to
include the following items in the notion of an algorithm which
computes the parametric greatest common divisor $h$ of the
algebraic family of polynomials $b_1\klk b_q$:
\begin{itemize}
\item an explicit representation of the rational functions $\om\theta m$ by
numerator and denominator polynomials belonging to $k[U_1\klk
U_r]$,
\item an explicit definition of the closed subvariety $S$ of $\A^r$ by
polynomials belonging to $k[\om U r]$.
\end{itemize}
The numerator and denominator polynomials representing the
rational functions $\om\theta m$ and the polynomials of $k[U_1\klk
U_r]$ defining the closed variety $S$ should then be
holomorphically encoded by a suitable data structure.

If there exists for the algebraic family of polynomials $b_1\klk
b_q$ a generically squarefree parametric greatest common divisor
$\widehat{h}$, we shall apply the same terminology to any
essentially division--free algorithm which computes $\widehat{h}$.

%-----------------------------------------------------------------
%-----------------------------------------------------------------
%-----------------------------------------------------------------
%-----------------------------------------------------------------

\subsection{A particular flat elimination problem.}
\label{section4.2}

Changing slightly the notations of Section \ref{section4.1} put
now $r:=n+1$, $T:=U_{n+1}$, $U:=(\om Un)$. Let us consider the
following polynomials of $\Q[T,U,X]$: $$G_1:=X_1^2-X_1\klk
G_n:=X_n^2-X_n,$$
\begin{equation}
\label{equation:uno4.2}
F_n:=\sum_{i=1}^n2^{i-1}X_i+T\prod_{i=1}^n\big(1+(U_i-1)X_i\big).
\end{equation}

It is clear from their definition that the polynomials $\om Gn$
and $F$ can be evaluated by a totally division--free arithmetic
circuit $\beta$ of size $O(n)$ in $\Q [T,U,X]$. Observe that the
polynomials $G_1,\ldots ,G_n$ do not depend on the parameters
$T,U_1,\ldots ,U_n$ and that their degree is two. The polynomial
$F_n$ is of degree $2n+1$. More precisely, we have $\deg_XF_n=n$,
$\deg_UF_n=n$, and $\deg_{T}F_n=1$. Although the polynomial $F_n$
may be evaluated by a totally division--free circuit of size
$O(n)$, the sparse representation of $F_n$, as a polynomial over
$\Q$ in the variables $T,\om Un,\om Xn$, contains $3^n$ nonzero
monomial terms and, as a polynomial over $\Q[T,\om Un]$ in the
variables $\om Xn$, it contains $2^n$ nonzero terms.\medskip

Let us now verify that the polynomials $\om G{n}$ and $F_n$ form
a flat family of elimination problems depending on the parameters
$T,U_1,\ldots ,U_n $.

The variety $V:= \{ G_1=0\klk G_n = 0\}$ is nothing but the union
of $2^n$ affine linear subspaces of $\A^{n+1} \times \A^{n}$, each
of them of the form $\A^{n+1} \times \{ \xi\}$, where $\xi$ is a
point of the hypercube $\{ 0,1 \} ^n$. The canonical projection
$\A^{n+1} \times \A^{n} \rightarrow \A^{n+1}$ induces a morphism
$\pi:\ V \rightarrow \A^{n+1}$ which glues together the canonical
projections $\A^{n+1} \times \{ \xi\} \rightarrow \A^{n+1}$ for
any $\xi$ in $\{ 0,1 \} ^n$. Obviously the morphism $\pi$ is
finite and  unramified. In particular $\pi$ has constant fibres
which are all canonically isomorphic to the
hypercube $\{ 0,1 \}^n$.\medskip

Let $(j_1\klk j_n)$ be an arbitrary point of $\{0,1\}^n$ and let
$j:=\sum_{1\leq i \leq n} j_i 2^{i-1}$ be the integer $0\leq j <
2^n$ whose bit representation is $j_n j_{n-1}\ldots j_1$. One
verifies immediately the identity $$F_n(T,U_1 ,\ldots , U_n , j_1
,\ldots , j_n) = j+T\prod_{i=1}^n U_i^{j_i}.$$ Therefore for any
point $(t,u_1,\ldots,u_n,j_1 ,\ldots,j_n)\in V$ with
$j:=\sum_{i=1}^{n} j_i 2^{i-1}$ we have $$F_n(t,u_1 ,\ldots,u_n,
j_1,\ldots,j_n)=j+t\prod_{i=1}^n u_i^{j_i}.$$

From this observation we deduce easily that the elimination
polynomial associated with the flat family of zero--dimensional
elimination problems determined by the polynomials $G_1 ,\ldots
,G_n$ and $F$ is in fact the polynomial $$P_n =
\prod_{j=0}^{2^{n}-1} (Y-(j+T\prod_{i=1}^n {U_i}^{[j]_i}))$$ of
Section \ref{section3.2}. With the notations of Section
\ref{section3.2}, this polynomial has the form
\begin{equation}
\label{equation:dos4.2} \begin{array}{rcl}P_n&=& Y^{2^n} + \sum_{1
\leq k \leq 2^n } B_k^{(n)} Y^{2^n -k} \\ &\equiv& Y^{2^n} +
\sum_{1 \leq k \leq 2^n } (\beta_k^{(n)} + TL_k^{(n)} ) Y^{2^n -k}
\ \textrm{modulo} \ T^2,
\end{array}
\end{equation}
with $\beta_k^{(n)}:= \sum_{1 \leq j_1 < \cdots < j_k
  < 2^n} j_1 \cdots j_k\ $ for $1 \leq k \leq 2^n$.\medskip

Let us consider $$\begin{array}{c}
\mo_n:=\Big\{F_n^{(t,u)};t\in\A^1,u:=(\om
un)\in\A^n,\qquad\qquad\qquad\qquad\qquad
\\ \qquad\qquad\qquad\qquad\qquad
F_n^{(t,u)}:=\sum_{i=1}^n2^{i-1}X_i+t\prod_{i=1}^n\big(1+(u_i-1)X_i\big)\Big\}
\end{array}$$ as input object class of our flat family of
zero--dimensional elimination problems, the affine space
$\A^1\times\A^n$ as input data structure and the map
$\omega_n:\A^1\times\A^n\to\mo_n$ defined for
$(t,u)\in\A^1\times\A^n$ by $\omega_n(t,u):=F_n^{(t,u)}$ as a
\Q--definable holomorphic encoding of the input object class
$\mo_n$.\medskip

Let us consider the set of univariate polynomials
$$\begin{array}{c}\mo_n^*:=\Big\{P_n^{(t,u)};t\in\A^1,u:=(\om
un)\in\A^n,\qquad\qquad\qquad\qquad\qquad
\\ \qquad\qquad\qquad\qquad\qquad\qquad P_n^{(t,u)}:=
\prod_{j=0}^{2^{n}-1} \big(Y-(j+t
\prod_{i=1}^n {u_i}^{[j]_i})\big)\Big\}
\end{array}$$ as output
object class and let be given a \Q--constructible output data
structure $\md_n^*$ of size $m_n^*$ and a \Q--definable,
holomorphic encoding $\omega_n^*:\md_n^*\to\mo_n^*$. Finally let
be given a totally division--free elimination procedure
$\theta_n:\A^1\times\A^n\to\md_n^*$ (in the sense of Section
\ref{section4.1}) which solves the zero--dimensional elimination
problem determined by the polynomials $G_1\klk G_n$ and $F_n$.
Observe that the size of our input data structure is $n+1$. With
this notations we have the following result:

\begin{theorem}
\label{lemma:section4.2} Assume that the elimination procedure
$\theta_n$ is robust in the sense of Definition \ref{robustel}.
Then the size $m_n^*$ of the output data structure $\md_n^*$
satisfies the estimate $$m_n^*\ge 2^n.$$
\end{theorem}
\begin{prf}
Since the arguments of this proof are similar to those used in
Section \ref{section3.2}, we shall be concise in our presentation.

Representing the univariate polynomials belonging to $\mo_n^*$ by
their coefficients we may identify the output object class $\mo^*$
with the corresponding subset of the ambient space $\A^{2^n}$.
With this interpretation the encoding $\omega_n^*:\md^*\to\mo^*$
becomes induced by a polynomial map from the affine space
$\A^{m_n^*}$ to the affine space $\A^{2^n}$. Therefore
$\omega_n^*\circ\theta_n:\A^{n+1}\to \A^{2^n}$ is a polynomial map
too.

Let $f_n:=\sum_{i=1}^n2^{i-1}X_i$, let
$$\displaystyle\beta_n:=(\beta_1^{(n)}\klk\beta_{2^n}^{(n)})=
\left(\sum_{1\le j_1\le\cdots\le j_k< 2^n}j_1\cdots j_k
\right)_{1\le k\le 2^n}$$ and let $u$ be an arbitrary point of
$\A^n$. From (\ref{equation:uno4.2}) and (\ref{equation:dos4.2})
we deduce $F_n^{(0,u)}=f_n$ and
$P_n^{(0,u)}=Y^{2^n}+\beta_1^{(n)}Y^{2^n-1}\plp\beta_{2^n}^{(n)}$.
In particular we have $\beta_n=\omega_n^*\big(\theta_n(0,u)\big)$.
This implies that the fiber
$(\omega_n^*\circ\theta_n)^{-1}(\beta_n)$ contains the
hypersurface $\{0\}\times\A^n$ of the affine space
$\A^1\times\A^n$. Moreover, since the elimination algorithm
$\theta_n$ is robust, we deduce from condition $(i)$ of Section
\ref{section4.1} and from $F_n^{(0,u)}=f_n$ that there are only
finitely many possible values for $\theta_n(0,u)$. More precisely,
the set $\{\theta_n(0,v);v\in\A^n\}$ is finite. Since
$\{0\}\times\A^n$ is irreducible we conclude now that there exists
an output code $\alpha\in\md^*$ with
$\theta_n(\{0\}\times\A^n)=\{\alpha\}$. From $\{0\}\times\A^n
\subset(\omega_n^*\circ\theta_n)^{-1}(\beta_n)$ we deduce
$\omega_n^*(\alpha)=\beta_n$.

Let $\gamma_u:\A^1\to\A^{m_n^*}$ and $\delta_u:\A^1\to\A^{2^n}$ be
the polynomial maps defined for $t\in\A^1$ by
$\gamma_u(t):=\theta_n(t,u)$ and
$\delta_u(t):=\omega_n^*\big(\theta_n(t,u)\big)$. We have
$\gamma_u(0):=\alpha$, $\delta_u(0)=\beta_n$ and
$\omega_n^*\circ\gamma_u=\delta_u$.

The following argumentation is exactly the same as in the proof of
Proposition \ref{remark3,3.2} of Section
\ref{section3.2}. First we deduce from (\ref{equation:dos4.2})
that
$$\big(L_1^{(n)}(u)\klk
L_{2^n}^{(n)}(u)\big)=\delta_u'(0)=(D\omega_n^*)_\alpha\big(
\gamma_u'(0)\big)$$ holds. Then we infer from Lemma
\ref{lemma:independence} that there exist points $u_1\klk
u_{2^n}\in\A^n$ such that the $(2^n\times 2^n)$--matrix
$\big(L_i^{(n)}(u_j)\big)_{1\le i,j\le 2^n}$ is nonsingular. This
implies that $\delta_{u_1}'(0)\klk\delta_{u_{2^n}}'(0)$ are
linearly independent elements of the \C--vector space $\A^{2^n}$.
Since $(D\,\omega_n^*)_\alpha:\A^{m_n^*}\to\A^{2^n}$ is a
\C--linear map, we conclude that
$\gamma_{u_1}'(0)\klk\gamma_{u_{2^n}}'(0)$ are linearly
independent elements of the \C--linear space $\A^{m_n^*}$ and
finally that $m_n^*\ge 2^n$ holds.
\end{prf}\cqfd

Suppose that there is given a procedure $\mathcal{P}$ which finds
for suitable encodings of input and output objects the solution
for each instance of any flat family of zero--dimensional
elimination problems. Suppose furthermore that the procedure
$\mathcal{P}$, applied to any flat family of zero--dimensional
elimination problems produces a robust (e.g. black box) algorithm
in the sense of Section \ref{section4.1} and that $\mathcal{P}$
can be applied to the encoding $\omega:\A^1\times\A^n\to\mo_n$ of
the input object class of the flat family of zero--dimensional
elimination problems (\ref{equation:uno4.2}). Then Theorem
\ref{lemma:section4.2} implies that $\mathcal{P}$ requires {\em
exponential} sequential time on infinitely many inputs. On the
other hand one sees easily that there do exist single exponential
time procedures of this kind (see \cite[Section 3.4]{GiHe01} and
the references cited there). Therefore the sequential time
complexity of zero--dimensional (parametric) elimination performed
by this kind of procedures is {\em intrinsically exponential}.
Observe in particular that this conclusion is valid for suitable
circuit encodings of input and output objects (see \cite[Theorem
1]{HeMaPaWa98} and \cite[Theorem 2]{GiHe01}).
\medskip

The sparse encoding of the object class $\mo_n$, defined by the
polynomial
$F_n=\sum_{i=1}^n2^{i-1}X_i+T\prod_{i=1}^n\big(1+(U_i-1)X_i\big)$,
is of size $3^n$. Therefore, from the point of view of
``classical'' parametric (i.e. branching--free) elimination
procedures (based on the sparse or dense encoding of polynomials
by their coefficients), it is not surprising that the sequential
 time becomes exponential in $n$ for the computation of the
solution of the general problem instance (\ref{equation:uno4.2}),
even if we change the data structure representing the output
objects (see e.g. \cite{GiHe93} and \cite{KrPa94}, \cite{KrPa96}
for this type of change of data structures).

Let us therefore look at the following flat family of
zero--dimensional elimination problems $\widetilde{G}_1\klk
\widetilde{G}_{3n-1},\widetilde{F}\in\Q[T,\om U n,\om X {3n-1}]$
in the parameters $T,\om U n$ and the variables $\om X {3n-1}$.
This family contains only sparse polynomials of at most four
monomial terms: $$\begin{array}{rrl}
\widetilde{G}_1&:=&X_1^2-X_1\klk\widetilde{G}_n:=X_n^2-X_n,\\ \\
\widetilde{G}_{n+1}&:=&X_{n+1}-2^1X_2-X_1,\\ \\
\widetilde{G}_{n+2}&:=& X_{n+2}-X_{n+1}-2^2X_3, \\ &\vdots & \\
\widetilde{G}_{2n-1}&:=&X_{2n-1}-X_{2n-2}-2^{n-1}X_n\\
\\ \widetilde{G}_{2n}&:=&X_{2n}-U_1X_1+X_1-1
=X_{2n}-\big(1+(U_1-1)X_1\big),\\
\\
\widetilde{G}_{2n+1}&:=&X_{2n+1}-U_2X_{2n}X_2+X_{2n}X_2-X_{2n}\\
&=&X_{2n+1}- X_{2n}\big(1+(U_2-1)X_2\big),\\ & \vdots \\
\widetilde{G}_{3n-1}&:=&X_{3n-1}-U_nX_{3n-2}X_n+X_{3n-2}X_n-X_{3n-2}\\
& = & X_{3n-1}-X_{3n-2}\big(1+(U_n-1)X_n\big)\\ \\
\widetilde{F}_n&:=&X_{2n-1}+TX_{3n-1}.\end{array}$$

One sees easily that the solution of the general problem instance
\linebreak
$\widetilde{G}_1=0\klk \widetilde{G}_{3n-1}=0,\widetilde{F}_n$ is
again the polynomial $P_n\in\Q[T,U,Y]$ of Section
\ref{section3.2}. The polynomials $\widetilde{G}_1\klk
\widetilde{G}_{3n-1},\widetilde{F}_n$ determine, with the
notations of Section \ref{section4.1}, the input object class
$$\widetilde{\mo}_n:=\left\{\big(\widetilde{G}_1^{(t,u)}\klk
\widetilde{G}_{3n-1}^{(t,u)},\widetilde{F}_n^{(t,u)}\big);
(t,u)\in\A^1\times\A^n\right\}$$ and the encoding
$\widetilde{\omega}_n:\A^1\times\A^n\to\widetilde{\mo}_n$ which
for $t\in\A^1$, $u\in\A^n$ is defined by
$\widetilde{\omega}_n(t,u):=\big(\widetilde{G}_1^{(t,u)}\klk
\widetilde{G}_{3n-1}^{(t,u)},\widetilde{F}_n^{(t,u)}\big)$. Since
these polynomials contain altogether exactly $9n-3$
monomials in the variables $X_1\klk X_{3n-1}$, we may consider
the input object class $\widetilde{\mo}_n$ as a \Q--constructible
subset of the affine space $\A^{9n-3}$. With this interpretation,
the encoding
$\widetilde{\omega}_n:\A^1\times\A^n\to\widetilde{\mo}_n$ becomes
\Q--definable and holomorphic. Applying to this situation the same
argumentation as in the proof of Theorem
\ref{lemma:section4.2} we conclude again that any branching--
and totally division--free, robust elimination procedure, which
finds from any input code $(t,u)\in\A^1\times\A^n$ the code of the
output object $P_n^{(t,u)}$ in a given data structure, requires
an output data structure of size at least $2^n$.
\medskip

Let us turn back to the polynomial
$F_n=\sum_{i=1}^n2^{i-1}X_i+T\prod_{i=1}^n\big(1+(U_i-1)X_i\big)$
of (\ref{equation:uno4.2}), to the object class
$\mo_n:=\big\{F_n^{(t,u)};t\in\A^1,u:=(\om un)\in\A^n\big\}$
defined by $F_n$ and to its encoding
${\omega}_n:=\A^1\times\A^n\to{\mo}_n$. Since the polynomial $F_n$
contains in the variables $\xon$ exactly $2^n$ nonzero monomial
terms we may consider $\mo_n$ as a \Q--constructible subset of
the affine space $\A^{2^n}$ and
${\omega}_n:=\A^1\times\A^n\to{\mo}_n$ as \Q--definable,
holomorphic encoding of the object class $\mo_n$. One sees easily
that $\mo_n$ is a closed, irreducible and \Q--definable subvariety
of $\A^{2^n}$ and that $\omega_n$ induces a robust encoding of
the object class $\mo_n\setminus\{\sum_{i=1}^n2^{i-1}X_i\}$ by
the data structure
$\A^1\times\A^n\setminus\big(\{0\}\times\A^n\big)$. On the other
hand $\{0\}\times\A^n$ is an exceptional fiber of the morphism of
algebraic varieties ${\omega}_n:\A^1\times\A^n\to{\mo}_n$.
Therefore the encoding $\omega_n$ of the object class $\mo_n$ is
not robust. On the other hand, by similar arguments as in Section
\ref{section3.2}, we may show that any {\em robust} encoding of
$\mo_n$ has size at least $2^n$. Therefore Theorem
\ref{lemma:section4.2} says only that any branching-- and
totally division--free {\em robust} elimination procedure
necessarily transfers a certain obstruction hidden in the given
encoding of the input object $F_n$ to the encoding of the output
object $P_n$.

However, Theorem \ref{lemma:section4.2} does {\em not}
say that the process of elimination creates a {\em genuine}
complexity problem for the encoding of the output object. In
particular we are not able to deduce from Theorem
\ref{lemma:section4.2} that the sequence of polynomials $\nseq
P$ is hard to evaluate. In fact, for $n\in\N$ the polynomial $P_n$
admits a short, \Q--definable, holomorphic encoding by the data
structure $\A^1\times\A^n$ and the sequence of polynomials $\nseq
P$ may in principle be easy to evaluate. However, in the latter
case, no branching-- and totally division--free, robust
elimination procedure will be able to discover this fact.

%-----------------------------------------------------------------
%-----------------------------------------------------------------
%-----------------------------------------------------------------
%-----------------------------------------------------------------

\subsection{The hardness of universal elimination.}
\label{section5.3} In this section we are going to show the second
main result of this paper, namely Theorem \ref{theorem:section5.3}
below, which says that there exists no universal polynomial
sequential time elimination algorithm $\mathcal{P}$ satisfying the
following condition:
\smallskip

\noindent $\mathcal{P}$ is able to compute equations for the
Zariski closure of any given constructible set and the generically
square--free parametric greatest common divisor of any given
algebraic family of univariate polynomials (see Sections
\ref{models} and \ref{section4.1} for the computational
model).
\medskip

The following considerations are devoted to the precise statement
and the proof of Theorem \ref{theorem:section5.3} below.

Let us suppose again $k:=\Q$ and $\ok:=\C$. Let $n$ be a fixed
natural number, let $m(n):=4n+10$ and let $T,\om Un,\om Xn$ and
$\om S{m(n)},Y$ be indeterminates over \Q. Let $U:=(\om Un)$,
$X:=(\xon)$, $S:=(\om S{m(n)})$ and let
$$R_n:=Z\left(\sum_{i=1}^n2^{i-1}X_i+T\prod_{i=1}^n
\big(1+(U_i-1)X_i\big)\right)\in\Q[Z,T,U,X].$$ One sees easily
that the polynomial $R_n$ may be evaluated by a totally
division--free arithmetic circuit of size $O(n)$.

Let $\widehat{\mo}_n$ be the Zariski closure of the set
$$\begin{array}{c}\displaystyle\Bigg\{R_n^{(z,t,u)};(z,t)\in\A^2,u=(\om
un)\in\A^n,\qquad\qquad\qquad\qquad\qquad\\ \\
\displaystyle\qquad\qquad\qquad\qquad
R_n^{(z,t,u)}:=z\left(\sum_{i=1}^n2^{i-1}X_i+t\prod_{i=1}^n
\big(1+(u_i-1)X_i\big)\right)\Bigg\}\end{array}$$ in a suitable
finite dimensional \C--linear subspace of $\C[X]$ and let
$\gamma_n:= (\gamma_1^{(n)}\klk
\gamma_{m(n)}^{(n)})\in\Z^{m(n)\times n}$ be an identification
sequence for $\widehat{\mo}_n$. From Corollary
\ref{corollary:tcs2} we deduce that such an identification
sequence exists and that we may assume without loss of generality
that the absolute values of the entries of the $\big(m(n)\times n
\big)$--matrix $\gamma_n$ are bounded by $3n^3$. Observe that
$\widehat{\mo}_n$ is a \Q--definable, irreducible, closed cone of
dimension at most $n+2$. We shall consider $\widehat{\mo}_n$ as
object class of $n$--variate polynomial functions.\medskip

Let us now consider the following prenex existential
formula $\Phi_n(S,Y)$ in the free variables $\om S{m(n)},Y$ and
the bounded variables $\xon,Z,$ $T,\om Un$:
$$\begin{array}{c}\displaystyle(\exists X_1) \cdots(\exists
X_n)(\exists Z)(\exists T)(\exists U_1) \cdots(\exists
U_n)\qquad\qquad\qquad\qquad\qquad\qquad\qquad\\ \displaystyle
\bigg(\bigwedge_{i=1}^nX_i^2-X_i=0\ \wedge \bigwedge_{k=1}^{m(n)}
S_k=R_n(Z,T,\om Un,\gamma_k^{(n)})\ \wedge \\
\displaystyle\qquad\qquad\qquad\qquad\qquad\qquad\wedge\
 Y=R_n(Z,T,\om Un,\xon)
\bigg).\end{array}$$

Using the previously mentioned arithmetic circuit encoding of the
polynomial $R_n$ and the bit encoding for integers, we see that
the length $|\Phi_n|$ of the formula $\Phi_n(S,Y)$ is $O(n^2)$.

Observe that the quantifier free formula $$Y=R_n(Z,T,\om Un,\om
Xn)$$ is equivalent to the following formula $$\Pi_n(Z,T,\om
Un,\om Xn,Y)$$ in the free variables $Z,T,\om Un,\xon,Y$ and the
bounded variables $X_{n+1}\klk X_{3n-1}$, i.e. both formulas
define the same subset of $\A^{2n+2}$ (compare Section
\ref{section4.2}):

$$\begin{array}{c}\displaystyle(\exists X_{n+1}) \cdots(\exists
X_{3n-1})\bigg(X_{n+1}-2X_2-X_1=0\
\wedge\qquad\qquad\qquad\qquad\qquad
\\\displaystyle\wedge\ \bigwedge_{j=n+2}^{2n-1}
X_j-X_{j-1}-2^{j-n}X_{j-n+1}=0\ \wedge
 X_{2n}-U_1X_1+X_1-1=0\ \wedge\ \\ \displaystyle\wedge\
\bigwedge_{k=2n+1}^{3n-1}X_k-U_{k-2n+1}X_{k-1}X_{k-2n+1}+X_{k-1}
X_{k-2n+1}-X_{k-1}=0\ \wedge \\
\qquad\qquad\qquad\qquad\qquad\qquad\qquad\qquad\displaystyle\wedge
\ Y=ZX_{2n-1} +ZTX_{3n-1} \bigg).\end{array} $$

Replacing now in the formula $\Phi_n(S,Y)$ for $1\le k\le
m(n)$ the occurrencies  of the subformulas $S_k=R_n(Z,T,\om
Un,\gamma_k^{(n)})$ by $\Pi_n(Z,T,U,\gamma_k^{(n)},$ $S_k)$ and the
occurrency of $Y=R_n(Z,T,\om Un,\xon)$ by \linebreak
$\Pi_n(Z,T,U,X,Y)$, we
obtain another prenex existential formula
$\widetilde{\Phi}_n(S,Y)$ in the free variables $\om S{m(n)},Y$
and $8n^2+20n-9$ bounded variables. The formula
$\widetilde{\Phi}_n(S,Y)$ has length
$|\widetilde{\Phi}_n|=O(n^2)$ for the sparse encoding of
polynomials and the bit representation of integers. Observe that
the formulas ${\Phi}_n(S,Y)$ and $\widetilde{\Phi}_n(S,Y)$ are
equivalent and asymptotically of the same length $O(n^2)$.
Thus the formulas $\Phi_n$ and $\widetilde{\Phi}_n$ are logical
expressions of asymptotically the same length which describe the
same constructible  subset of the affine space
$\mathbb{A}^{m(n)}\times \A^1$. The polynomials occurring in
$\Phi_n(X,Y)$ are given in arithmetic circuit encoding, whereas
the polynomials occurring in $\widetilde{\Phi}_n(X,Y)$ are given
in sparse encoding. The formulas ${\Phi}_n(X,Y)$ and
$\widetilde{\Phi}_n(X,Y)$ will be the inputs for an elimination
problem which we are now going to describe in detail.\medskip

In the sequel we shall restrict our attention to the
formula ${\Phi}_n(S,Y)$. Our considerations will be
identically valid for the formula $\widetilde{\Phi}_n(S,Y)$.\medskip

Let $\widehat{\sigma}_n:\widehat{\mo}_n\to\A^{m(n)}$ be the map
defined for $R\in\widehat{\mo}_n$ by
$\widehat{\sigma}_n(R):=\big(R(\gamma_1^{(n)})\klk
R(\gamma_{m(n)}^{(n)})\big)$ and let $\widehat{\md}_n$ be the
image of $\widehat{\sigma}_n$. From Lemma \ref{lemma4:encode} we
deduce that $\widehat{\md}_n$ is a \Q--definable, irreducible,
closed cone of $\A^{m(n)}$ and that $\widehat{\sigma}_n$ induces a
finite, bijective morphism of algebraic varieties
$\widehat{\sigma}_n:\widehat{\mo}_n\to\widehat{\md}_n$ which is
therefore a homeomorphism with respect to the Zariski topologies
of $\widehat{\mo}_n$ and $\widehat{\md}_n$. In particular
$\widehat{\sigma}_n:\widehat{\mo}_n\to\widehat{\md}_n$ is a
homogeneous, birational map. We consider $\hmd_n$ as a
\Q--definable data structure and
$(\widehat{\sigma}_n)^{-1}:\hmd_n\to\widehat{\mo}_n$ as a
\Q--definable, continuous encoding of the object class
$\widehat{\mo}_n$. From a similar argument as in the proof of
Proposition \ref{remark3,3.2} we deduce that
$\widehat{\sigma}_n:\widehat{\mo}_n\to\widehat{\md}_n$ is not an
isomorphism of affine varieties. Thus
$\widehat{\sigma}_n^{-1}:\widehat{\md}_n\to\widehat{\mo}_n$ is not
a {\em holomorphic} encoding of the object class $\widehat{\mo}_n$
by the data structure $\widehat{\md}_n$, but only a {\em
continuous} one (compare Theorem \ref{theorem:encode1}). This
circumstance contributes to a certain technical intricateness of
the argumentation which now follows.\medskip

Observe first that the prenex existential formula $(\exists
Y)\Phi_n(S,Y)$ describes a \Q--constructible subset of $\A^{m(n)}$
whose Zariski closure is $\widehat{\md}_n$. Observe then that
$dim\, \widehat{\md}_n\le dim\, \widehat{\mo}_n\le
n+2<4n+10=m(n)$ holds. Therefore $\widehat{\md}_n$ is strictly
contained in the affine space $\A^{m(n)}$. Thus the formula
$\Phi_n(S,Y)$ introduces an implicit semantical dependence
between the indeterminates $\om S {m(n)}$. In the sequel we shall
consider the indeterminates $S_1\klk$ $S_{m(n)}$ as parameters
and $Y$ as variable.
\medskip

Let us now consider an arbitrary point $s=(\om s
{m(n)})\in\A^{m(n)}$ which satisfies the formula $(\exists
Y)\Phi_n(S,Y)$. Then there exist points $(z,t)\in\A^2$ and $u=(\om
u {n})\in\A^{n}$ such that the $n$--variate polynomial
$$R_n^{(z,t,u)}=z\left(\sum_{i=1}^n2^{i-1}X_i+
t\prod_{i=1}^n\left(1+(u_i-1)X_i\right)\right)$$ satisfies the
condition $$s_1=R_n^{(z,t,u)}(\gamma_1^{(n)})\klk s_{m(n)} =
R_n^{(z,t,u)}(\gamma_{m(n)}^{(n)}).$$

Since $\gamma_n=(\gamma_1^{(n)}\klk\gamma_{m(n)}^{(n)}) \in
\Z^{m(n)\times n}$ is an identification sequence for the object
class $\mo_n$, we conclude that the polynomial $R_n^{(z,t,u)}\in
\C[\xon]$ depends only on the point $s\in\A^{m(n)}$ and not on its
particular encoding $(z,t,u)$ belonging to the data structure
$\A^2\times\A^n$. We write therefore $R_n^{(s)}:=R_n^{(z,t,u)}$.
Let $$\displaystyle\widehat{P}_n^{(s)}:=\prod_{(\om\varepsilon n)
\in\{0,1\}^n}\left(Y-R_n^{(s)}(\om\varepsilon n)\right)$$ and let
us write $\Phi_n(s,Y)$ for the formula of the elementary language
of algebraically closed fields of characteristic zero with
constants in \C\ which is obtained by specializing in the formula
$\Phi_n(S,Y)$ the variables $S_1\klk S_{m(n)}$ into the values
$s_1\klk s_{m(n)}\in\C$. Observe that the polynomial
$\widehat{P}_n^{(s)}$ is monic of degree $2^n$ and $\Phi_n(s,Y)$
contains a single free variable, namely $Y$. One sees easily that
the formula $\Phi_n(s,Y)$ is equivalent to the quantifier--free
formula $\widehat{P}_n^{(s)}(Y)=0$.

Observe that for a suitable point $s=(\om s {m(n)})$ of
$\A^{m(n)}$ satisfying the formula $(\exists Y) \Phi_n(S,Y)$
(e.g. choosing $s$ such that for $f_n:=\sum_{i=1}^n 2^{i-1}X_i$
the condition $s_1=f_n(\gamma_1^{(n)})\klk s_{m(n)}=
f_n(\gamma_{m(n)}^{(n)})$ is satisfied) we obtain a univariate
{\em separable} polynomial $\widehat{P}_n^{(s)}$ of degree $2^n$.
This implies that there exists a nonempty Zariski open subset
$\mau$ of the closed, \Q--definable, irreducible subvariety
$\widehat{\md}_n$ of the affine space $\A^{m(n)}$ such that
$\mau$ is contained in the \Q--constructible subset of
$\A^{m(n)}$ defined by the formula $(\exists Y)\Phi_n(S,Y)$ and
such that for any point $s\in\mau$ the polynomial
$\widehat{P}_n^{(s)}\in\C[Y]$ is monic and separable of degree
$2^n$. Since the \Q--definable morphism
$\widehat{\sigma}_n:\widehat{\mo}_n\to\widehat{\md}_n$ is finite,
bijective and birational, there exists a polynomial
$\widehat{R}_n\in\Q(\widehat{\md}_n)[X]$ satisfying the following
two conditions:

\begin{itemize}
  \item for any point $s\in\widehat{\md}_n$, any coefficient
  $\rho$ of $\widehat{R}_n$ and any place \linebreak
$\varphi:\C(\widehat{\md}_n)
   \to\C\cup\{\infty\}$ whose valuation ring
contains the local
   ring of $\hmd_n$ at the point $s$, the value
$\varphi(\rho)$ is finite
   and uniquely determined by $s$.
  \item if additionally the point $s$ satisfies the formula
  $(\exists Y)\Phi_n(S,Y)$, then the polynomial $\varphi(\widehat{R}_n)
  \in\C[X]$, obtained by specializing the coefficients of
$\widehat{R}_n$ by means of the place $\varphi$, satisfies the
equation
 $\varphi(\widehat{R}_n)=R_n^{(s)}$.
\end{itemize}

With these notations, we shall write
$R_n^{(s)}:=\varphi(\widehat{R}_n)$
 also if $s$ does not satisfy the
formula $(\exists Y)\Phi_n(S,Y)$.
 Let $$\widehat{P}_n:=
\prod_{(\om\varepsilon n)
\in\{0,1\}^n}\left(Y-\widehat{R}_n(\om\varepsilon
n)\right)\in\Q(\widehat{\md}_n)[Y].$$ One sees easily that
 $\widehat{P}_n$
satisfies mutatis mutandis the above two
 conditions (note that in the second
condition $\varphi(\widehat{R}_n)=
 R_n^{(s)}$ has to be replaced by
$\varphi(\widehat{P}_n)=P_n^{(s)}$).

In particular the coefficients of $\widehat{P}_n$ belong to the
integral closure of the domain $\Q[\widehat{\md}_n]$ in its
function field $\Q(\widehat{\md}_n)$.

Since for any $s\in\mau$ the polynomial $\widehat{P}_n^{(s)}$ is
separable of degree $2^n$, we conclude that $\widehat{P}_n$ is a
monic, separable polynomial of degree $2^n$ in the variable $Y$.\medskip

Let us denote by $V_n$ the Zariski closure of the
\Q--constructible subset of $\A^{m(n)}\times\A^1$ defined by the
formula $\Phi_n (S,Y)$ and by
$\pi_n:\A^{m(n)}\times\A^1\to\A^{m(n)}$ the canonical projection
which maps each point of $\A^{m(n)}\times\A^1$ on its first $m(n)$
coordinates. Observe that $V_n$ is nonempty and that the
$\Q$--definable, irreducible variety $\widehat{\md}_n$ is the
Zariski closure of $\pi_n(V_n)$ in $\A^{m(n)}$. Let $C$ be any
irreducible component of $V_n$ satisfying the condition
$\overline{\pi_n(C)}=\widehat{\md}_n$ (observe that such an
irreducible component exists). Let us now fix a point $s\in\hmd_n$
which we think chosen generically between the points of $\hmd_n$.
From this choice we infer immediately that the set
$\pi_n^{-1}(s)\cap C$ is not empty and that its elements satisfy
the formula $\Phi_n(S,Y)$. One now sees easily that
$\pi_n^{-1}(s)\cap C$ is a nonempty and finite set. This implies
$dim\, C=dim\,\widehat{\md}_n$.\medskip

Observe that for any point $s\in \A^{m(n)}$ satisfying $(\exists
Y)\Phi_n(S,Y)$, the formula ${\Phi}_n(s,Y)$ is equivalent to the
quantifier free formula $\widehat{P}_n^{(s)}(Y)=0$.

Therefore we shall consider from now on $\widehat{P}_n$ as the
canonical output object associated to the elimination problem
given by the formula $\Phi_n(S,Y)$ in the parameters $S_1,\ldots,
S_{m(n)}$ and the single variable $Y$. More precisely, our
elimination task will consists in the computation of the
polynomial $\widehat{P}_n^{(s)}\in \C[Y]$ for any input instance
$s\in \widehat{\mathcal{D}}_n$. In this sense, we are looking for
output data structures which solve problem $(ii)$ of Section
\ref{section1.2} for the object class defined by the polynomial
$\widehat{P}_n$ in the parameter instances defined by the formula
$(\exists Y)\Phi_n(S,Y)$.\medskip

Suppose now that there is given an elimination procedure
$\mathcal{P}$ which is universal and branching--parsimonious
in the sense of Section \ref{section1.2}. Suppose furthermore
that $\mathcal{P}$ accepts as inputs prenex existential input
formulas of the elementary theory of algebraically closed fields
of characteristic zero, whose terms are polynomials in arithmetic
circuit representation (or alternatively polynomials in sparse
representation). Assume that $\mathcal{P}$ is associated with a
suitable output data structure which allows the holomorphic
encoding of polynomials and with a monotone sequential time
measure $\mathcal{T}$, and suppose that $\mathcal{P}$ and
$\mathcal{T}$ satisfy the following conditions:
\begin{itemize}
  \item[(1)] Let $\Phi$ be a given prenex existential formula
  of the elementary language $\mathcal{L}$ of algebraically closed
  fields of characteristic zero with constants 0, 1. Suppose
  that the polynomial terms occurring in the formula $\Phi$
  are encoded by the {\em input data structure} associated with the
  elimination procedure $\mathcal{P}$. Then the elimination procedure $\mathcal{P}$
  produces a quantifier--free formula $\Psi$ whose polynomial terms
  are (holomorphically) encoded by the output data structure
  associated with $\mathcal{P}$, such that $\Phi$ and $\Psi$ are equivalent
  formulas. The length $|\Psi|$ of the output
  formula $\Psi$ satisfies the estimate $|\Psi|\le
  \mathcal{T}(|\Phi|)$.
  \item[(2)] Let $\Xi\in\mathcal{L}$ be a quantifier--free formula
  whose polynomial terms are encoded by the {\em output data
  structure} associated with $\mathcal{P}$. Then the procedure
  $\mathcal{P}$ produces from the input $\Xi$ a system of polynomial
  equations $\mathcal{F}$, encoded by the output data structure associated with
  $\mathcal{P}$, such that $\mathcal{F}$ defines the Zariski closure of
  the \Q--constructible set defined by $\Xi$. The size
  $|\mathcal{F}|$ of the system of polynomial equations $\mathcal{F}$
  satisfies the estimate $|\mathcal{F}|\le\mathcal{T}(|\Xi|)$.
  \item[(3)] Let $\mathcal{B}$ be a system of polynomials
  encoded by the {\em output data structure} associated with the
  elimination procedure $\mathcal{P}$. Suppose that $\mathcal{B}$
  represents an algebraic family of univariate polynomials, for which
  a generically square--free parametric greatest common divisor
  $\widehat{h}$ in the sense of Section \ref{section2} exists.
  Then the procedure $\mathcal{P}$ produces from the input $\mathcal{B}$
  an algorithm in the sense of Section \ref{section5.2} which computes
  for any admissible input instance of $\mathcal{B}$ the generically
  square--free greatest common divisor $\widehat{h}$ of the algebraic
  family of univariate polynomials represented by $\mathcal{B}$. Here
  we assume implicitly  that $\widehat{h}$ is represented by a
  polynomial $H$ which is encoded by the output data structure
  associated with the procedure $\mathcal{P}$. With respect to this data structure
  the size $|H|$ of the polynomial $H$ satisfies the estimate
  $|H|\le\mathcal{T}(|\mathcal{B}|)$.
\end{itemize}

Let us remark that condition (1) above characterizes $\mathcal{P}$
as a universal elimination procedure in the usual sense, whereas
conditions (2) and (3) state that $\mathcal{P}$ solves suitable
elimination problems of type $(ii)$ of Section \ref{section1.2}.
In principle, input and output formulas mentioned in condition (1)
may be represented by algorithms which admit branchings. If for
example $\mathcal{P}$ uses as input and output data structures for
the encoding of polynomials arithmetic circuits, quantifier--free
(sub--)formulas in condition (1) may be represented by arithmetic
networks (arithmetic--boolean circuits, see \cite{Gathen86},
\cite{Gathen93}). Nevertheless we require that the outputs
mentioned in conditions (2) and (3) represent branching--free
evaluation procedures. All known universal elimination procedures
satisfy with respect to a suitably defined sequential time
complexity measure conditions (1), (2), (3) above.\medskip

Let us now apply the given elimination procedure
$\mathcal{P}$ to the input formula $\Phi_n(S,Y)$. Since
$\mathcal{P}$ satisfies condition (1), the output is a
quantifier--free formula $\Psi_n(S,Y)$ of the elementary language
$\mathcal{L}$, such that $\Psi_n(S,Y)$ is equivalent to
$\Phi_n(S,Y)$. Moreover, the polynomial terms occurring in the
formula $\Psi(S,Y)$ are represented by the output data structure
associated with $\mathcal{P}$.\medskip

We apply now the procedure $\mathcal{P}$ to the
quantifier--free formula $\Psi_n(S,Y)$. Since $\mathcal{P}$
satisfies condition (2), the output is a finite set $\mb_n$ of
polynomials of $\Q[S,Y]$ which define the algebraic variety $V_n$.
Again, the polynomials contained in $\mathcal{B}_n$ are
represented by the output data structure associated with
$\mathcal{P}$.
\medskip

Recall that $V_n$ is the Zariski closure of the
\Q--constructible subset of $\A^{m(n)}\times\A^1$ defined by the
formula $\Phi_n(S,Y)$ (and hence by the formula $\Psi_n(S,Y)$),
that $V_n$ is nonempty, that $\widehat{\md}_n$ is the Zariski
closure of the image of $V_n$ under the canonical projection
$\pi_n:\A^{m(n)} \times\A^1\to\A^{m(n)}$ and that any irreducible
component $C$ of $V_n$ with $\overline{\pi_n(C)}=\widehat{\md}_n$
satisfies the condition $dim\,C=dim\,\widehat{\md}_n$. Let
$B_1^{(n)}\klk B_{q_n}^{(n)}$ be the elements of $\mathcal{B}_n$
which do not vanish identically on the algebraic variety
$\widehat{\md}_n\times\A^1$. Let $b_1^{(n)}\klk b_{q_n}^{(n)}$ be
the univariate polynomials of $\Q[\widehat{\md}_n][Y]$ induced by
$B_1^{(n)}\klk B_{q_n}^{(n)}$ on $\hmd_n\times\A^1$. Observe that
$b_1^{(n)}\not=0\klk b_{q_n}^{(n)}\not=0$ holds. Since any
irreducible component $C$ of $V_n$ with
$\overline{\pi_n(C)}=\widehat{\md}_n$ satisfies the condition
$dim\,C=dim\,\hmd_n$ and since such an irreducible component
exists, we conclude $q_n\ge 1$. Therefore $b_1^{(n)}\klk
b_{q_n}^{(n)}$ is an algebraic family of univariate polynomials in
the sense of Section \ref{section5.2}.

Let $h\in\Q(\hmd_n)[Y]$ be the greatest common divisor of the
polynomials $b_1^{(n)}\klk b_{q_n}^{(n)}$ in $\Q(\hmd_n)[Y]$.
Since for any point $s\in\mau$ the formula $\Phi_n(s,Y)$ (and
hence the formula $\Psi_n(s,Y)$) is equivalent to the formula
$\widehat{P}_n^{(s)}(Y)=0$ and since $\mau$ is a nonempty Zariski
open subset of $\hmd_n$, we conclude that the monic polynomials
$h$ and $\widehat{P}_n$ of $\Q(\hmd_n)[Y]$ have the same roots in
any algebraic closure of the field $\Q(\hmd_n)$.

Therefore we have $\deg h\ge\deg\widehat{P}_n=2^n$. Thus the
degree of $h$ in the variable $Y$ is positive.

Since the univariate polynomial $\widehat{P}_n$ is separable, we
conclude, from our previous considerations concerning the
definition of $P_n^{(s)}$ for arbitrary $s\in\hmd_n$, that there
exists a generically square--free greatest common divisor for the
algebraic family of univariate polynomials $b_1^{(n)}\klk
b_{q_n}^{(n)}$ and that this greatest common divisor is
$\widehat{P}_n$.\medskip

Finally we apply the procedure $\mathcal{P}$ to the finite
set of polynomials $\mb_n$. Since $\mathcal{P}$ satisfies
condition (3), the output are rational functions
$\widehat{\theta}_1^{(n)}\klk\widehat{\theta}_{\widehat{m}_n}^{(n)}$
of $\Q(\hmd_n)$, such that $\widehat{\theta}_n:=
(\widehat{\theta}_1^{(n)} \klk
\widehat{\theta}_{\widehat{m}_n}^{(n)})$ represents an algorithm
in the sense of Section \ref{models} which computes the
generically square--free parametric greatest common divisor
$\widehat{P}_n$ of the algebraic family of univariate polynomials
$b_1^{(n)}\klk b_{q_n}^{(n)}$ for any admissible input instance
(which necessarily belongs to $\hmd_n$). Let us make this
statement more precise:
\smallskip

\noindent let $\md_n^*$ be the $\Q$--constructible output data
structure associated with the procedure $\mathcal{P}$, when
$\mathcal{P}$ is applied to the input $\mb_n$. Then the size of
$\md_n^*$ is $\widehat{m}_n$ and we may suppose without loss of
generality that $\md_n^*$ is a closed subvariety of the affine
space $\A^{\widehat{m}_n}$. Then the closure of the image of
$\widehat{\theta}_n$ is contained in $\md_n^*$ and therefore we
may interpret $\widehat{\theta}_n$ as a dominant rational map from
$\widehat{\md}_n$ to $\md_n^*$.
\smallskip

The output data structure $\md_n^*$ encodes a suitable output
object class $\mo_n^*$ of univariate polynomials by means of a
\Q--definable, holomorphic encoding $\omega_n^*:\md_n^*\to
\mo_n^*$. The object class $\mo_n^*$ contains the set
$\{\widehat{P}_n^{(s)};s\in\widehat{\md}_n\}$. Since
$\widehat{\theta}_n$ is a dominant rational map, the composition
$\omega_n^*\circ\widehat{\theta}_n$ is well defined and
$\omega_n^*\circ\widehat{\theta}_n$ is a rational map from
$\hmd_n$ to the Zariski closure of the object class $\mo_n^*$ in a
suitable affine ambient space. By assumption the algorithm
$\widehat{\theta}_n$ computes the generically square--free
parametric greatest common divisor $\widehat{P}_n\in\Q(\md_n)[Y]$
of the algebraic family of univariate polynomials $b_1^{(n)}\klk
b_{q_n}^{(n)}$ for any admissible input instance. Thus for any
input instance $s\in\hmd_n$ and any place $\varphi:\C(\hmd_n)
\to\C\cup\{\infty\}$ whose valuation ring contains the local ring
of $\hmd_n$ at $s$, the values $\varphi(\widehat{\theta}_1^{(n)})
\klk \varphi(\widehat{\theta}_{\widehat{m}_n}^{(n)})$ are finite
and uniquely determined by the input instance $s$. With these
notations we may therefore consistently write
$\widehat{\theta}_1^{(n)}(s):=\varphi(\widehat{\theta}_1^{(n)})\klk
\widehat{\theta}_{\widehat{m}_n}^{(n)}(s):=
\varphi(\widehat{\theta}_{\widehat{m}_n}^{(n)})$ and
$\widehat{\theta}_n(s):=(\widehat{\theta}_1^{(n)}(s)\klk
\widehat{\theta}_{\widehat{m}_n}^{(n)}(s))$. Since $\md_n^*$ is a
closed subvariety of $\A^{\widehat{m}_n}$ we have
$\widehat{\theta}_n(s)\in\md_n^*$ for any $s\in\hmd_n$. We may
therefore interpret $\widehat{\theta}_n$ as a {\em total} map from
$\hmd_n$ to $\md_n^*$ whose value is defined for {\em any}
argument from $\hmd_n$. With this interpretation $\omega_n^*\circ
\widehat{\theta}_n$ is a total map from $\hmd_n$ to $\mo_n^*$
satisfying the condition $\omega_n^*\circ\thetah_n(s)=\omega_n^*
\big(\thetah_n(s)\big)=\widehat{P}_n^{(s)}$ for any
$s\in\hmd_n$.\medskip

The above considerations imply that the rational functions
$\thetah_1^{(n)}\klk \thetah_{\widehat{m}_n}^{(n)}$ belong to the
integral closure of the domain $\Q[\hmd_n]$ in its fraction field
$\Q(\hmd_n)$. Moreover they imply that the rational function
$\thetah_n$, which we may suppose well--defined for the Zariski
open subset $\mau$ of $\hmd_n$, represents an essentially
division--free algorithm in the sense Section \ref{section5.2}
which computes for each input instance $s\in\mau$ a code
$\thetah_n(s)$ for the output object $\widehat{P}_n^{(s)}$ and
which can be uniquely extended to the
limit data structure $\hmd_n$ of $\mau$.\medskip

Observe now, that specializing in the polynomial
$$R_n=Z\left(\sum_{i=1}^n2^{i-1}X_i+T\prod_{i=1}^n
\big(1+(U_i-1)X_i\big)\right)\in\Q[Z,T,U,X]$$ the variable $Z$
into the value one, we obtain the polynomial
$$F_n=\sum_{i=1}^n2^{i-1}X_i+ T\prod_{i=1}^n
\big(1+(U_i-1)X_i\big)\in\Q[T,U,X]$$ introduced in Section
\ref{section4.2}. Therefore the object class $\mo_n:=\{
F_n^{(t,u)};t\in\A^1,u\in\A^n\}$ is contained in the object class
$\widehat{\mo}_n$. Since $\widehat{\mo}_n$ is Zariski closed in
its ambient space we have $\omo_n\subset\widehat{\mo}_n$. Let
$\md_n:=\widehat{\sigma}_n(\overline{\mo}_n)$. Since
$\overline{\mathcal{\mo}}_n$ is a \Q--definable, closed,
irreducible subvariety of $\widehat{\mo}_n$ and
$\widehat{\sigma}_n:\widehat{\mo}_n\to\widehat{\md}_n$ is a
\Q--definable, finite, bijective morphism of algebraic varieties,
we conclude that $\md_n$ is a (nonempty) \Q--definable, closed,
irreducible subvariety of $\hmd_n$. For any point $s\in\md_n$ and
any place $\varphi:\C(\hmd_n)\to\C\cup\{\infty\}$ whose valuation
ring contains the local ring of $\hmd_n$ at $s$, and any
coefficient $\beta$ of $\widehat{P}_n\in\Q(\hmd_n)[Y]$, the values
of $\varphi(\beta)$ and of $\varphi(\widehat{\theta}_1^{(n)})\klk
\varphi(\widehat{\theta}_{\widehat{m}_n}^{(n)})$ are finite and
{\em uniquely} determined by $s$. Therefore there exists a monic
polynomial $\check{P}_n\in\Q(\md_n)[Y]$ of degree $2^n$, rational
functions
$\check{\theta}_1^{(n)}\klk\check{\theta}_{\widehat{m}_n}^{(n)}$
and a nonempty Zariski open subset $\mau_0$ of $\md_n$ such that
$\check{P}_n$ and $\check{\theta}_n:=(\check{\theta}_1^{(n)}\klk
\check{\theta}_{\widehat{m}_n}^{(n)})$ are well defined in any
point $s$ of $\mau_0$ and such that the conditions
$\check{P}_n^{(s)}=\widehat{P}_n^{(s)}$ and
$\check{\theta}_n^{(s)}=\widehat{\theta}_n^{(s)}$ are satisfied.

Since $\Q[\md_n]$ is a holomorphic image of $\Q[\hmd_n]$ and since
the coefficients of $\widehat{P}_n$ and $\thetah_1^{(n)}\klk
\thetah_{\widehat{m}_n}^{(n)}$ belong to the integral closure of
$\Q[\hmd_n]$ in $\Q(\hmd_n)$, we conclude that
$\check{\theta}_1^{(n)}\klk \check{\theta}_{\widehat{m}_n}^{(n)}$
and the coefficients of $\check{P}_n$ belong to the integral
closure of the domain $\Q[\md_n]$ in its fraction field
$\Q(\md_n)$. In the same way one sees that $\check{\theta}_n$
represents an essentially division--free algorithm which computes
the polynomial $\check{P}_n$ and which can be uniquely extended to
the limit data structure $\md_n$ of $\mau_0$. For any $s\in\md_n$
we infer therefore that
$\check{\theta}_n(s):=(\check{\theta}_1^{(n)}\klk
\check{\theta}_{\widehat{m}_n}^{(n)})$ is a well defined point of
$\md_n^*$ and that $\check{P}_n^{(s)}$ is a well defined, monic,
univariate polynomial of degree $2^n$ satisfying the conditions
$\check{\theta}_n(s)=\widehat{\theta}_n(s)$ and
$\check{P}_n^{(s)}=\widehat{P}_n^{(s)}$.\medskip

Consider now the \Q--definable, holomorphic encoding
$\omega_n:\A^1\times\A^n\to\mo_n$ of the object class $\mo_n$ by
the data structure $\A^1\times\A^n$, defined for
$(t,u)\in\A^1\times\A^n$ by $\omega_n(t,u):=F_n^{(t,u)}$ (see
Section \ref{section4.2}).

Observe that
$\widehat{\sigma}_n\circ\omega_n:\A^1\times\A^n\to\md_n$ is a
dominant morphism of \Q--definable, irreducible varieties.
Therefore $\theta_1^{(n)}:=\check{\theta}_1^{(n)}\circ
\widehat{\sigma}_n\circ\omega_n\klk \theta_{\widehat{m}_n}^{(n)}:=
\check{\theta}_{\widehat{m}_n}^{(n)}\circ
\widehat{\sigma}_n\circ\omega_n$ are well--defined rational
functions belonging to $\Q(T,U)$. Observe that
$\widehat{\sigma}_n\circ\omega_n$ induces a \Q--algebra
isomorphism which maps the coordinate ring $\Q[\md_n]$ onto the
subdomain $$\mathcal{A}_n:=\Q[F_n(T,U,\gamma_1^{(n)})\klk
F_n(T,U,\gamma_{m(n)}^{(n)})]$$ of the polynomial ring
$\Q[T,U]$.

Since the rational functions $\check{\theta}_1^{(n)} \klk
\check{\theta}_{\widehat{m}_n}^{(n)}$ belong to the integral
closure of $\Q[\md_n]$ in $\Q(\md_n)$, we conclude that
$\theta_1^{(n)} \klk \theta_{\widehat{m}_n}^{(n)}$ belong to the
integral closure of $\mathcal{A}_n$ in $\Q(T,U)$. But $\Q[T,U]$ is
integrally closed in its fraction field. This implies that
$\theta_1^{(n)} \klk \theta_{\widehat{m}_n}^{(n)}$ are {\em
polynomials} belonging to $\Q[T,U]$. Thus
$\theta_n:=(\theta_1^{(n)} \klk \theta_{\widehat{m}_n}^{(n)})$
defines a morphism of algebraic varieties $\theta_n:\A^1\times\A^n
\to\md_n^*$ which satisfies for any point $(t,u)\in\A^1\times\A^n$
the identities $$\begin{array}{rcl}
\omega_n^*\big(\theta_n(t,u)\big)&=&\omega_n^*
\bigg(\check{\theta}_n\Big(\widehat{\sigma}_n
\big(\omega_n(t,u)\big)\Big)\bigg)\\
\\
&=& \omega_n^*\bigg(\thetah_n\Big(\widehat{\sigma_n}
\big(\omega_n(t,u)\big)\Big)\bigg)
\\ \\
&=&
\widehat{P}_n^{(\widehat{\sigma}_n\circ\omega_n)(t,u)}\\
\\&=& \displaystyle\prod_{(\varepsilon_1
\klk\varepsilon_n)\in\{0,1\}^n}\left(Y-
R_n^{(\widehat{\sigma}_n\circ\omega_n)(t,u)} (\varepsilon_1
\klk\varepsilon_n)\right).\end{array}$$

Let $t\in\A^1$ and $u=(u_1\klk u_n)\in\A^n$
be fixed for the
 moment. Observe that
$R_n^{(\widehat{\sigma}_n\circ\omega_n)(t,u)}$ is the {\em
unique} polynomial of the object class
 $\widehat{\mo}_n$
which satisfies the condition
 $$\left(R_n^{(\widehat{\sigma}_n\circ\omega_n)(t,u)}
(\gamma_1^{(n)})\klk R_n^{(\widehat{\sigma}_n\circ\omega_n)(t,u)}
(\gamma_{{m}(n)}^{(n)})\right)=
\widehat{\sigma}_n\circ\omega_n(t,u).$$ On the other hand we have
$$\widehat{\sigma}_n\circ\omega_n(t,u)=\left(F_n^{(t,u)}
(\gamma_1^{(n)})\klk F_n^{(t,u)}(\gamma_{{m}(n)}^{n})\right)$$
and $F_n^{(t,u)}\in\widehat{\mo}_n$. This implies
$R_n^{(\widehat{\sigma}_n\circ\omega_n)(t,u)}=F_n^{(t,u)}$ and
therefore we have
$$\begin{array}{rcl}
\displaystyle\omega_n^*\left(\theta_n(t,u)\right)&=&
\displaystyle\prod_{(\varepsilon_1
\klk\varepsilon_n)\in\{0,1\}^n}\left(Y-R_n^{
(\widehat{\sigma}_n\circ\omega_n)(t,u)} (\varepsilon_1
\klk\varepsilon_n)\right)\\ \\ &=&\displaystyle
\prod_{(\varepsilon_1 \klk\varepsilon_n)\in\{0,1\}^n}
\left(Y-F_n^{(t,u)}(\varepsilon_1
\klk\varepsilon_n)\right)\\ \\
&=& \displaystyle
\prod_{j=1}^{2^n-1}\left(Y-(j+t\prod_{i=1}^n
u_i^{[j]_i})\right).\end{array}$$

 Let
$P_n:=\prod_{j=1}^{2^n-1}\left(Y-(j+T
\prod_{i=1}^nU_i^{[j]_i})\right)\in\Q[T,U,Y]$ be the elimination
polynomial introduced in Sections \ref{section3.2} and
 \ref{section4.2}. Then we have \linebreak
$\omega_n^*\big(\theta_n(t,u)\big)
 =P_n^{(t,u)}$ for any point
$(t,u)\in\A^1\times\A^n$. Taking now
 $\A^1\times\A^n$ as input data
structure,
 $\theta_n(\A^1\times\A^n)$ as output data structure,
\linebreak $\{P_n^{(t,u)};t\in\A^1,u\in\A^n\}$ as output object
class encoded by the restriction of $\omega_n^*$ to the
\Q--definable subset $\theta_n(\A^1\times\A^n)$ of $\md_n^*$, we
see now that these data structures are \Q--constructible, that the
encoding is \Q--definable and holomorphic and that $\theta_n$
represents a {\em totally division--free} algorithm which computes
for each input code $(t,u)$ of $\A^1\times\A^n$ an output code
$\theta_n(t,u)$ which encodes the output object $P_n^{(t,u)}$.
Thus $\theta_n$ is a totally division--free elimination procedure
which computes the general solution $P_n$ of the flat family of
zero--dimensional elimination problems given by the equations
$X_1^2-X_1=0\klk X_n^2-X_n=0$ and the polynomial $F_n$ (see
Section \ref{section4.2}).\medskip

Recall that the polynomials
$\theta_1^{(n)}\klk\theta_{\widehat{m}_n}^{(n)}$ belong to the
integral closure of $\mathcal{A}_n=\Q[F_n(T,U,\gamma_1^{(n)})\klk
F_n(T,U,\gamma_{{m}(n)}^{(n)})]$ in $\Q[T,U]$ and that the
\Q--algebra $\mathcal{A}_n$ is canonically isomorphic to the
coordinate ring $\Q[\omo_n]$ of the Zariski closure of the object
class $\mo_n$. Therefore $\theta_n$ is a {\em robust} elimination
procedure in the sense of Definition \ref{robustel}. From Theorem
\ref{lemma:section4.2} we deduce now the
 estimate $\widehat{m}_n\ge
2^n$.\medskip

Since by assumption the sequential time complexity measure
$\mathcal{T}$ is monotone, we conclude now that $$2^n\le
\widehat{m}_n=|\widehat{P}_n|\le
\mathcal{T}(|\mathcal{B}_n|)\le\mathcal{T}^2(|\Psi_n|)\le
\mathcal{T}^3(|\Phi_n|)\le\mathcal{T}^3(cn^2)$$ holds for a
suitable universal constant $c>0$.

Therefore $\mathcal{T}$ cannot be a polynomial function. Finally
we remark that the same conclusion is valid if we replace in our
argumentation the formula $\Phi_n$ by the formula
$\widetilde{\Phi}_n$. We may now summarize these considerations by
the following general result:

\begin{theorem}
\label{theorem:section5.3} Let $\mathcal{P}$ be a universal
elimination procedure for the theory of algebraically closed
fields of characteristic zero with constants $0$, $1$ and let
$\mathcal{T}$ a sequential time complexity measure for
${\mathcal{P}}$. Suppose that $\mathcal{P}$ accepts as inputs
prenex existential formulas whose polynomial terms are given in
arithmetic circuit or sparse representation. Suppose that
${\mathcal{P}}$ and $\mathcal{T}$ satisfy conditions (1), (2), (3)
above. Then $\mathcal{T}$ is not a polynomial function.
\end{theorem}

%-----------------------------------------------------------------
%-----------------------------------------------------------------
%-----------------------------------------------------------------
%-----------------------------------------------------------------
%-----------------------------------------------------------------
%-----------------------------------------------------------------
%-----------------------------------------------------------------
%-----------------------------------------------------------------

\section{Conclusions.} There exists a general opinion between computer
scientists that proving lower complexity bounds for specific
problems defined by existential prenex formulas (see
\cite{Borodin93}) is an extremely difficult task which requires
tricky methods or deep mathematical insight. Simple minded
algorithmic models and the absence of operative notions of
uniformity make in our opinion excessively intricate or impossible
to prove striking complexity results for many fundamental
algorithmic problems of practical interest. A way out of this
dilemma consists in the restriction of the computational model
under consideration. Thus one may for example think to consider
only unbounded fan--in and fan--out arithmetic circuits of {\em
bounded depth} for the computation of polynomials of interest, as
e.g. the resultant of two generic univariate polynomials or more
generally, the general solution of a flat family of
zero--dimensional elimination problems.

Asymptotically optimal lower sequential time complexity bounds
become then easy to prove. However the restriction to bounded
depth circuits represents a highly artificial limitation of the
computational model (this restriction excludes for example the
evaluation of monomials of high degree by means of iterated
squaring) and the complexity result obtained in this way becomes
irrelevant as a guide for future software developers.\medskip

The ultimate aim of this paper was not a theoretical but a
practical one. We tried to give a partial answer to the following
fundamental question:

\smallskip
\noindent what has to be changed in elimination theory in order
to obtain practically efficient algorithms?
\smallskip

We established a list of implicit or explicit requirements
satisfied by all known (symbolic or seminumeric) elimination
algorithms. These requirements are: universality, no branchings
and robustness for certain simple elimination problems, capacity
of computing certain closures (as e.g. equations for the Zariski
closure of a given constructible set or the greatest common
divisor of two polynomials). Moreover, by means of a suitable
preparation
 of the input equation, all known universal
elimination
procedures may be transformed easily into Kronecker--like
procedures which are able to evaluate the corresponding canonical
elimination polynomial in any given argument or to compute its
coefficients. In this sense the known elimination procedures are
all able to ``compute canonical elimination polynomials".

The fulfillment of these requirements and the capacity of
computing canonical elimination polynomials implies the
experimentally certified non--polynomial complexity character of
these elimination procedures and explains their practical
inefficiency. The results of this paper demonstrate that the
complexity problem we focus on is not a question of optimization
of algorithms and data structures. There is no way out of the
dilemma by changing for example from dense to sparse elimination
or to fewnomial theory. Hybridization of symbolic and numeric
algorithms leads us again back to the same complexity problems we
started from.

In this sense the paper is devoted to the elaboration and
discussion of a series of ``uniformity'' notions which restrict
the (mostly implicit) computational models relevant for the
present (and probably also the future) design of implementable
elimination procedures in algebraic geometry. Emphasis was put on
the motivation of these algorithmic restrictions and not on the
mathematical depth of the techniques used in this paper in order
to prove lower complexity bounds. In fact, it turns out that
elementary methods of classical algebraic geometry are sufficient
to answer the complexity questions addressed in this paper. It is
not the first time that a refined analysis of the complexity model
produces not only elementary and simpler proofs of lower bound
results in algebraic complexity theory, but also stronger
complexity statements. Examples are the ``elementarizations'' of
Strassen's degree method \cite{Strassen73A}, due to Sch\"onhage
\cite{Schonhage76} and Baur \cite[Theorem 8.5]{BuClSh97}, and the
combinatorial method of Aldaz and Monta\~na   for the
certification of the hardness of univariate polynomials (compare
\cite[Chapter
 9]{BuClSh97} with \cite{AlHeMaMoPa98b} and
\cite{AlMaMoPa01}).\medskip

Nevertheless there are two points addressed in this paper,
which call for the development of deep new tools in mathematics
and computer science: the problem of algorithmic modeling
addressed in Section \ref{section3.3} calls for the search of
mathematical statements which generalize Hilbert's Irreducibility
Theorem to (not necessarily unirational) algebraic varieties
containing ``many'' integer or rational points and to the
characterization of unirational varieties (in the sense of
\cite{Kollar99}) by means of arithmetic properties.\medskip

On the other hand our discussion of the notion of robustness
of elimination procedures in Section \ref{section4.1} leads to the
question in which sense the concept of programmable function can
be distinguished from the notion of elementarily recursive
function (here the concepts of specification and data type make
the main difference). A programmable function appears always
together with a certificate (``correctness proof'') that it meets
its specification. The existence of such a proof necessarily
restricts the syntactical form of the underlying program and hence
the complexity model in which the running time of the program is
measured.

%-----------------------------------------------------------------
%-----------------------------------------------------------------
%-----------------------------------------------------------------
%-----------------------------------------------------------------
%-----------------------------------------------------------------
%-----------------------------------------------------------------
%-----------------------------------------------------------------
%-----------------------------------------------------------------
\appendix
\section{Appendix.} \label{appendix}
\subsection{Universal correct test and identification sequences.}
\label{unisequence}

In this section we are going to formulate a slight generalization
of the main results of Section \ref{correct} and
\ref{circuitcorrect} namely Lemma \ref{ts}, Corollary
\ref{corollary:tcs2} and Theorem \ref{theorem:encode1}. These
generalizations are based on  Baire's Theorem and lead to the
concept of {\em universal} correct test and identification
sequence.

\begin{corollary}
\label{corollary:tcs1} Let $k:=\Q$, $\overline{k}:=\C$ and let
$L$, $m$, $t$ be given natural numbers with $m>L$. Then there
exists a subset $S\subset\R^{mt}$ satisfying the following
conditions:
\begin{itemize}
\item[$(i)$] $S$ is dense in the strong topology of $\R^{mt}$.
\item[$(ii)$] any element $\gamma=(\om\gamma m)\in S$ with
$\om\gamma m\in\R^t$ is a \cts\ for the $\Q$--Zariski closure of
any $\Q$--constructible object class $\mo$ of $t$--variate
polynomial functions over $\C$ such that for $\mo$ there exists a
$\Q$--definable holomorphic encoding by a data structure of size
$L$.
\end{itemize}
\end{corollary}

A \cts\ as in Corollary \ref{corollary:tcs1}, $(ii)$ is called
{\em universal} for the corresponding set of object classes.

\begin{prf}
Observe that there are only countably many $\Q$--definable
holomorphic encodings of $\Q$--constructible object classes of
polynomial functions in $t$ variables over $\C$ by data structures
of size $L$. Therefore we may think these encodings enumerated as
$\omega_1,\omega_2,\dots$. From the second part of the proof of
Lemma \ref{ts} of Section \ref{holocorrect} we conclude that there
exists for any $i\in\N$ a $\Q$--definable, Zariski open, dense
subset $\mathcal{U}_i\subset\C^{mt}$ such that any element
$\gamma= (\om\gamma m)$ of $\mathcal{U}_i$ with $\om\gamma
m\in\C^t$ is a \cts\ for the $\Q$--Zariski closure of the object
class of $t$--variate polynomial functions over $\C$ encoded by
$\omega_i$. Observe now that
$\mathcal{U}_i^*:=\mathcal{U}_i\cap\R^{mt}$ is open and dense in
the strong topology of $\R^{mt}$. Let $S:=
\displaystyle\cap_{i\in\N}\mathcal{U}_i^*$. From Baire's Theorem
we deduce that the set $S$ is still dense in the strong topology
of $\R^{mt}$.

Let $\gamma=(\om\gamma m)$ be an arbitrary element of $S$ with
$\om\gamma m\in\R^t$ and let $\mo$ be an arbitrary
$\Q$--constructible object class of $t$--variate polynomials over
$\C$ such that for $\mo$ there exists a $\Q$--definable
holomorphic encoding by a data structure of size $L$. Then there
exist an index $i\in\N$ such that $\omega_i$ encodes $\mo$. From
$S\subset\mathcal{U}_i$ we deduce that $\gamma$ is a \cts\ for
the object class $\omo$. In conclusion $\gamma$ is a universal
\cts\ of length $m$ for the set of object classes under
consideration.
\end{prf}\cqfd

\begin{corollary}
\label{corollary:tcs3} Let $k:=\Q$, $\overline{k}:=\C$ and let
$L$, $m$, $t$ be given natural numbers with $m>2L$. Then there
exists a subset $S\subset\R^{mt}$ satisfying the following
conditions:
\begin{itemize}
\item[$(i)$] $S$ is dense in the strong topology of $\R^{mt}$.
\item[$(ii)$] Any element $\gamma=(\om\gamma m)\in S$ with
$\om\gamma m\in\R^t$ is an identification sequence for the
$\Q$--Zariski closure of any $\Q$--constructible object class
$\mo$ of $t$--variate polynomial functions over $\C$ such that for
$\mo$ there exists a $\Q$--definable holomorphic encoding by a
data structure of size $L$.
\end{itemize}
\end{corollary}

An identification sequence as in Corollary \ref{corollary:tcs3},
$(ii)$ is called {\em universal} for the corresponding set of
object classes.

The proof of Corollary \ref{corollary:tcs3} combines the statement
of Corollary \ref{corollary:tcs1} with the same arguments employed
in the proof of Corollary \ref{corollary:tcs2} of Section
\ref{holocorrect} and is omitted here.\medskip

In a similar way one may combine Corollary \ref{corollary:tcs3}
and Lemma \ref{lemma4:encode} of Section \ref{valencode} in order
to prove the following statement:

\begin{corollary}
\label{coro:encode2} Let $k:=\Q$, $\ok:=\C$ and let $L$, $m$, $t$
be given natural numbers with $m>2L$. Then there exists a subset
$S\subset\R^{mt}$ satisfying the following conditions:
\begin{itemize}
\item[$(i)$] $S$ is dense in the strong topology of $\R^{mt}$,
\item[$(ii)$] any element $\gamma=(\gammam)\in S$ with $\gammam\in
\R^t$ has the following property:
\smallskip

\noindent let $\mo$ be an arbitrary $\Q$--constructible object
class of $t$--variate polynomial functions over $\C$ such that for
$\mo$ there exists a $\Q$--definable holomorphic encoding by a
data structure of size $L$ and suppose that $\mo$ is a cone. Let
$\sigma:\omo\to\A^m(\C)$ be the map defined by
$\sigma(F):=\big(F(\gamma_1)\klk F(\gamma_m)\big)$ for $F\in\omo$
and let $\md^*:=\sigma(\omo)$. Then $\md^*$ is a cone of
$\A^m(\C)$  which is closed in the $\C$--Zariski topology of
$\A^m(\C)$ (and hence also in the strong topology) and $\sigma$
defines a bijective finite morphism of $\mo$ onto $\md^*$. For any
$\C$--irreducible component $\mc$ of $\omo$ the restriction map
$\sigma:\mc\to\sigma (\mc)$ is a birational (finite and bijective)
morphism of $\mc$ onto the $\C$--irreducible Zariski closed set
$\sigma(\mc)$. The encoding of the object class $\omo$ by the data
structure $\md^*$ defined by $\omega^*:= \sigma^{-1}$ is
continuous with respect to the $\C$--Zariski topologies of $\omo$
and $\md^*$. Moreover $\omega^*$ is holomorphic if and only if
$\omega^*$ allows to answer holomorphically the value question
about the object class $\omo$. Finally $\omega^*$ induces an
encoding of the projective variety associated to the cone $\omo$
by the projective variety associated to the cone $\md^*$ which is
continuous with respect to the strong topology.
\end{itemize}
\end{corollary}

%-----------------------------------------------------------------
%-----------------------------------------------------------------
%-----------------------------------------------------------------
%-----------------------------------------------------------------

\subsection{The VC--dimension of a holomorphically encoded object class.}
\label{VCdimension} Let $\mo$ be a $k$--constructible object
class of polynomial functions. We say that a finite set
$A\subset\A^t$ can be {\em shattered} by the object class $\mo$
if for each subset $A'\subset A$ there exists an object $F\in\mo$
such that any element $a\in A$ belongs to $A'$ if and only if
$F(a)=0$ holds. We define the {\em Vapnik--Chervonenkis (VC)
dimension} $dim_{VC}\mo$ of $\mo$ as infinite if there exist
subsets $A$ of $\A^t$ of arbitrary cardinality which can be
shattered by $\mo$. Otherwise we define $dim_{VC}\mo$ as the
maximal cardinality  of such a set (see \cite[Chapter 3,
3.6]{Vapnik00} and \cite[Chapter 3, 3.5]{BuClSh97} for details).
The following statement implies that the VC--dimension of the
object class $\mo$ is finite.

\begin{lemma}
\label{vc} Let notations and assumptions be as in Lemma \ref{ts}
of Section \ref{holocorrect}. Then $dim_{VC}\mo$ satisfies the
following estimate:
$$(dim_{VC}\mo)^{\frac{1}{2}}\le\frac{dim_{VC}\mo}{\log
dim_{VC}\mo} \le L(1+\log \Delta_2)$$ (here $\log$ denotes the
logarithm to the base 2).
\end{lemma}
\begin{prf}
We shall freely use the notations of the proof of Lemma \ref{ts}.
Let $s\in\N$ with $s\le dim_{VC}\mo$. Then there exists a finite
set $A\subset \A^t$ of cardinality $s$ which can be shattered by
$\mo$. Let $A=\{a_1\klk a_s\}$ with $a_1\klk a_s\in\A^t$. From the
construction of the ambient space $\A^N$ of $\mo$ we deduce that
there exists a $k$--definable (evaluation) map
$eval:\A^N\times\A^t\to\A^1$ which satisfies the condition
$eval(F,y)=F(y)$ for any polynomial $F\in\mo$ and any point
$y\in\A^t$. This implies that for any $a\in A$ there exists a
polynomial $\Omega_a\in \ok[Z_1\klk Z_L]$ of degree at most
$\Delta_2$ such that for any $D\in\md$ the identity
$\Omega_a(D)=eval(\Omega(D),a)=\omega(D)(a)$ holds.

Let $A'$ be an arbitrary subset of $A$. By hypothesis there exists
a polynomial $F\in\mo$ with $A'=\{a\in A;F(a)=0\}$. Consider
$$\md_{A'}:=\{D\in\md;\Omega_a(D)=0\mbox{ for }a\in A',\,
\Omega_a(D)\not=0\mbox{ for }a\in A\setminus A'\}.$$ Any code
$D\in\md$ with $\omega(D)=F$ belongs to $\md_{A'}$. Therefore
$\md_{A'}$ is nonempty. Thus $\md_{A'}$ is a
$\Omega_{a_1}\klk\Omega_{a_s}$--cell in the sense of
\cite{Heintz83}. From \cite[Theorem 2]{JeSa00} or \cite[Corollary
1]{Heintz83} one deduces that the number of
$\Omega_{a_1}\klk\Omega_{a_s}$--cells is bounded by
$(1+s\Delta_2)^L$. Since the set $A$  can be shattered by $\mo$
and different subsets of $A$ define disjoint
$\Omega_{a_1}\klk\Omega_{a_s}$--cells we conclude $$2^s\le
(1+s\Delta_2)^L.$$ This implies $\frac{s}{\log s}\le
L(1+\log\Delta_2)$. From $s\le dim_{VC} \mo$ we deduce now
$$(dim_{VC}\mo)^{\frac{1}{2}}\le\frac{dim_{VC}\mo}{\log
dim_{VC}\mo} \le L(1+\log \Delta_2).$$
\end{prf}\cqfd

Let $k:=\Q$ and $\ok:=\C$. We are going to consider the data
structure $\md_{real}:=\md\cap\R^L$ and the object class
$\mo_{real}:=\omega(\md_{real})$. Observe that $\mo_{real}$ is a
$\Q$--definable semialgebraic subset of $\R^N$. The standard
definition of the VC--dimension of $\mo_{real}$ is slightly
different from our definition of the VC--dimension of $\mo$ (see
\cite[Chapter 3, 3.6]{Vapnik00}). Taking into account the number
of different real cells of a system of $s$ real polynomials of
degree at most $\Delta_2$ in $L+1$ variables is of order
$O\left(\left(\frac{s\Delta_2}{L+1}\right)^{ L+1}\right)$ (see
\cite{PoRo93}), one concludes in the same way as in the proof of
Lemma \ref{vc} that
\begin{equation}
\label{vc1}
(dim_{VC}\mo_{real})^{\frac{1}{2}}\le\frac{dim_{VC}\mo_{real}}{\log
dim_{VC}\mo_{real}} \le (L+1)\log \Delta_2+
O\left(\frac{1}{\log dim_{VC}\mo_{real}}\right)
\end{equation}
holds.
\smallskip

Let $\overline{W}_{L,t}$ be the set of all polynomials
$F\in\okyot$ which have approximative nonscalar (sequential)
complexity over $\ok$ at most $L$. From Corollary
\ref{corollary:tcs4} and its proof we conclude that
$\overline{W}_{L,t}$ is a $k$--constructible object class which
has a $k$--definable, holomorphic encoding of size $4(L+t+1)^2+2$
by means of polynomials over $k$ of degree at most $L2^{L+1}+2$.
Thus Lemma \ref{vc} implies the estimate
$$dim_{VC}\overline{W}_{L,t}\le 8(L+t+1)^{3+\varepsilon}$$ for any
$\varepsilon>0$. From \cite[Chapter
 9, Proposition 9.1]{BuClSh97} we infer that
any univariate polynomial of
 $\okyot$ of degree at most $\frac{L^2}{4}$
belongs to $W_{L,t}$
 and hence to $\overline{W}_{L,t}$.

Let $A\subset\A^t$ be a subset of $s:=\lfloor\frac{L^2}{4}\rfloor$
elements of the form $A:=\{(a_i,0\klk 0);a_i\in k, 1\le i\le s\}$
(here $\lfloor\frac{L^2}{4}\rfloor$ denotes the largest integer
below $\frac{L^2}{4}$). Then for any subset $A'$ of $A$ there
exists a polynomial $F\in k[Y_1]$ of degree $\# A'$ such that
$A':=\{a\in A;F(a)=0\}$ holds. From $\deg F=\# A'\le\frac{L^2}{4}$
we deduce $F\in\overline{W}_{L,t}$. This consideration implies
finally $$\frac{L^2}{4}-1<dim_{VC}\overline{W}_{L,t}\le
8(L+t+1)^{3+\varepsilon}$$ for any $\varepsilon>0$.

Let $k:=\Q$ and $\ok:=\C$. We consider the set ${W}^{real}_{L,t}$
of all polynomials of $\R[Y_1\klk Y_t]$ which can be evaluated by
a totally division--free arithmetic circuit of nonscalar size at
most $L$ using only scalars from $\R$. Thus we have
${W}^{real}_{L,t}=({W}_{L,t})_{real}$. Taking into account the
estimate (\ref{vc1}),
%and \cite[Theorem 2.4]{KaWe93} we conclude
we conclude in a similar way as before that
$$\frac{L^2}{4}-1<dim_{VC}\overline{W}^{real}_{L,t}\le
8(L+t+1)^{3+\varepsilon}+ O\left(\frac{1}{\log
dim_{VC}\overline{W}_{L,t}^{real}}\right)$$ holds for any 
$\varepsilon>0$.

Analogous considerations lead to an upper bound for the set of
polynomials of $\R[Y_1\klk Y_t]$ which have approximative
complexity at most $L$ in terms of essentially division--free
arithmetic circuits using only parameters from $\R(\varepsilon)$.
\bigskip

%-----------------------------------------------------------------
%-----------------------------------------------------------------
%-----------------------------------------------------------------
%-----------------------------------------------------------------
%-----------------------------------------------------------------
%-----------------------------------------------------------------
%-----------------------------------------------------------------
%-----------------------------------------------------------------

\noindent {\bf Acknowledgment.} The authors wish to thank Rosa
Wachenchauzer who pointed to us the relation between complexity
theory and program specification and generation. They are
especially grateful to the anonymous referees for many useful
suggestions which helped to improve considerably the presentation
of the results of this paper. J. Heintz and G. Matera thank to the
Facultad de Ingenier\'\i a, Ciencias Exactas y Naturales,
Universidad Favaloro, where they did part of this work.

{\footnotesize
%\bibliography{refs,newref}
%\bibliographystyle{alpha}

\newcommand{\etalchar}[1]{$^{#1}$}

}

\end{document}